\documentclass[11pt]{amsart}
\usepackage[a4paper,margin=24mm]{geometry}
\usepackage{amsmath,amssymb,amsthm,mathtools}
\usepackage{xcolor}
\usepackage{graphicx}
\usepackage{booktabs,tikz,longtable,enumitem}
\usepackage{pgfplots}
\usepackage[section]{placeins}
\AddToHook{cmd/subsection/before}{\FloatBarrier}
\pgfplotsset{compat=1.18}
\usepgfplotslibrary{groupplots}
\usepackage[T1]{fontenc}\usepackage{lmodern,microtype,hyperref}
\usetikzlibrary{arrows.meta,calc}
\allowdisplaybreaks
\numberwithin{equation}{section}
\newtheorem{theorem}{Theorem}[section]
\newtheorem{proposition}[theorem]{Proposition}
\newtheorem{lemma}[theorem]{Lemma}
\newtheorem{corollary}[theorem]{Corollary}
\theoremstyle{definition}\newtheorem{remark}[theorem]{Remark}
\newcommand{\C}{\mathbb C}
\newcommand{\ten}{\otimes_A}
\newcommand{\id}{\mathrm{id}}
\newcommand{\Ric}{\operatorname{Ric}}
\newcommand{\Scal}{\operatorname{Scal}}
\newcommand{\Om}{\Omega}
\newcommand{\extd}{{\rm d}}

\newcommand{\Aut}{{\rm{Aut}}}

\newcommand{\R}{{\mathbb{R}}}

\newcommand{\Z}{{\mathbb{Z}}}

\newcommand{\F}{{\mathbb{F}}}

\newcommand{\CC}{{\mathcal{C}}}

\newcommand{\cu}{{\mathfrak u}}

\renewcommand{\ker}{{\rm{ker}}}

\newcommand{\tens}{\otimes}

\renewcommand{\imath}{\mathrm{i}}
\newcommand{\emath}{\mathrm{e}}

\allowdisplaybreaks
\newtheorem{conjecture}[theorem]{Conjecture}
\theoremstyle{definition}\newtheorem{definition}[theorem]{Definition}
\newcommand{\im}{\mathrm{im}}
\newcommand{\Span}{\operatorname{span}}

\newcommand{\rank}{\operatorname{rank}}

\title{Quantum Riemannian geometries on graphs over $\Bbb F_2$}
\author{Shahn Majid}
\address{School of Mathematical Sciences\\ Queen Mary University of London \\ Mile End Rd, London E1 4NS }
\email{s.majid@qmul.ac.uk}

\begin{document}
\maketitle

\begin{abstract}
A graph is Eulerian when every vertex has even degree and we similarly
call a graph co-Eulerian when any two distinct vertices have an even
codegree. We show that a connected graph $G$
admits a quantum Riemannian geometry (QRG) with calculus $\Omega_{\min}$
over the field $\F_2=\{0,1\}$ if and only if it is Eulerian or
co-Eulerian or both. Each type has a base quantum Levi-Civita
connection (QLC) $\nabla$ which is inner and $*$-preserving in an
appropriate sense. Every other $*$-preserving QLC of the chosen type
is of the form $\nabla+\alpha$, where $\alpha\in V_*(G)$, a linear
space with basis given by the connected components of the 
anti-Gallai graph $D_*(G)$. We show for large order $n$ that
almost all QRGs are Eulerian and that almost every connected Eulerian graph
has $V_*(G)$ one-dimensional. We show that the probability of a uniformly chosen connected graph admitting a QLC is asymptotic to $2^{1-n}$. Our results also suggest that asymptotically almost
all QLCs are $*$-preserving and that at even order the co-Eulerian graphs are asymptotically equally divided between Eulerian and odd-degree graphs, related by a de Morgan type duality.  We prove that the fraction of \(*\)-preserving QLCs which are flat decays exponentially in \(n^2\).  We also conjecture that for odd $n>3$ the only connected flat QRGs over $\F_2$ are the canonical geometries on the $n$-gons, and prove this for $n<17$. We display all maps between graph quantum Riemannian geometries over $\F_2$  up to six vertices in the form of a quiver.
\end{abstract}

\section{Introduction}

Quantum Riemannian geometry (QRG)\cite{BegMa} is a more explicit approach to noncommutative geometry more generally\cite{Con} that works over any field, starting with a choice of exterior algebra $(\Omega,\extd)$, metric $g\in \Omega^1\tens_A\Omega^1$, and then proceeding where possible to a bimodule quantum Levi-Civita connection (QLC) $\nabla:\Omega^1\to \Omega^1\tens_A\Omega^1$. If this exists it need not in general be unique. Here $A$ can be noncommutative but can also be commutative but not commuting with differential forms. It can, for example, be the algebra of functions on the vertices of a graph and the arrows (where every edge is viewed as bidirected) providing $\Omega^1$. There is a canonical choice $\Omega_{min}$ that extends this to a DGA for any graph, with further quotients of interest for specific graphs\cite{Ma:gra,BegMa}. Recently we have shown over $\C$ that the connected trees that admit QRGs with positive-ratio metrics are exactly the Dynkin graphs of type  A,D,E and affine $\tilde D,\tilde E$, see\cite{Ma:tre}. The problem for a general connected graph, however, remains formidable.

 In the present work, we attack this general graph case but from the opposite end of working over the field $\F_2$. In this case there is only one metric (of value 1 on every arrow), so the QRG amounts to finding all possible QLCs for it. There is also a natural $*$ operation which now just amounts to reversal of arrows. For the calculus $\Omega_{\min}$, we give a complete solution to this problem in graph-theoretic terms. We call a graph {\em Eulerian} or {\em even-degree} if all its vertices have even degree, without requiring connectedness. On the other hand, on any graph it is natural to define codegree as
\[ {\rm codeg}(x,y)=| N(x)\cap N(y)|\]
where $N(x)$ is the set of neighbours of a vertex $x$. If vertices $x,y$ coincide then this is just the usual degree. We define a graph as {\em co-Eulerian} if its codegree is even for distinct vertices $x,y$. Our main classification result for connected graphs $G$ with at least two vertices (Theorems~\ref{thm:general-F2-classification} and~\ref{thm:general-alpha-kernel}) is:

\begin{enumerate}
\item Exactly Eulerian or co-Eulerian graphs $G$ admit a QRG over $\F_2$.

\item Each case has a unique braiding $\sigma:\Omega^1\tens_A\Omega^1\to \Omega^1\tens_A\Omega^1$ that provides an inner, $*$-preserving QLC $\nabla$. 

\item In each case, all QLCs are of the form $\nabla+\alpha$ where the moduli of $\alpha$  is a  linear space that can be identified in the $*$-preserving case with the vector space spanned by the $d_*$ connected components  of the anti-Gallai graph with isolated vertices removed, denoted $D_*(G)$~\cite{Le}. Its vertices are the edges of $G$ belonging to triangles, joined when they lie in a common triangle. 
\end{enumerate}

The two cases correspond to $q=0,1$ in a single construction $\sigma_q$, with QLCs of the form $\nabla_q+\alpha$. Both types contribute when a graph is both Eulerian and co-Eulerian. Also, if we are interested in QLCs without the $*$-preserving condition then $\alpha$ lies in a space spanned by the $d$ connected components of an undirected graph $D(G)$ whose vertices are the arrows in triangles of $G$, joined when they lie in a common directed triangle. The map of inner Eulerian QRGs
\[ D(G)\to D_*(G)\]
 is a two-sheeted cover. Its split components
measure the difference between the two moduli dimensions (Theorem~\ref{thm:derived-cover}).

This classification also yields an asymptotic result independent of
the computations: almost every connected graph admitting a QLC is
Eulerian. Among all connected graph isomorphism classes of order
\(n\), the proportion admitting QLCs is asymptotic to \(2^{1-n}\)
(Proposition~\ref{prop:qlc-asymptotic}). Moreover, almost every
connected Eulerian graph has \(d=d^*=1\), and almost every graph
admitting QLCs has exactly two, both star-preserving
(Theorem~\ref{thm:typical-eulerian}). Co-Eulerian graphs contribute
\(o(e_n)\) connections even when weighted by their full QLC families,
where \(e_n\) counts connected Eulerian graphs
(Proposition~\ref{prop:coeu-counts}). The remaining issue
for the conjectured total star-preserving count \(2e_n\)
(Conjecture~\ref{conj:star-asymptotic}) is the
exponentially weighted tail of \(d^*(G)\) on Eulerian graphs.

Flatness imposes further restrictions. Apart from two non-star-preserving
connections on \(K_3=C_3\), every flat QLC is inner and star-preserving.
For Eulerian geometries, flatness characterises cycles. For co-Eulerian geometries,
we give an exact graph criterion, prove regularity and a parity
restriction on triple intersections of neighbourhoods, and bound
the degree by \(2\lfloor n/4\rfloor\) when the underlying graph is Eulerian
(Theorem~\ref{thm:flat-alpha} and
Lemma~\ref{lem:flat-sigma1} and~\ref{prop:flat-regular}).
For co-Eulerian geometries with all distinct-vertex codegrees at most two, flatness gives
a covering by a hypercube and forces the graph order to be a power
of two. We show that the only connected graphs with every vertex of degree four admitting flat co-Eulerian geometries are
\(Q_4\) and \(K_4\square K_2\)
(Proposition~\ref{prop:flat-cube} and Corollary~\ref{cor:flat-four}). We also show that, while both types of braiding are involutive and Eulerian braidings always obey the braid relation, co-Eulerian braidings do so if and only if the corresponding inner connection $\nabla$ is flat (Proposition~\ref{prop:braid}).   We conjecture that cycles are the only connected odd-order graphs
admitting flat star-preserving QLCs (Conjecture~\ref{conj:odd-flat}),
and prove this for every odd order below seventeen, as well as for
Cayley graphs of odd-order groups with ad-invariant generating sets
(Propositions~\ref{prop:flat-regular} and~\ref{prop:cayley-odd}). Proposition~\ref{prop:flat-frequency} shows that 
the fraction of star-preserving QLCs which are flat decays exponentially
in \(n^2\), without assuming the conjecture. Complete
graphs and cycles are classified explicitly in
Theorem~\ref{thm:Knunique} and Proposition~\ref{prop:Cn}. The graph criteria allow
counts through thirteen vertices without generating individual
connections; these provide evidence for the remaining asymptotic
questions in Section~\ref{ssec:graph-counts}.

Section~\ref{sec:QRG-maps} looks at general maps between QRGs over $\F_2$. We find that maps preserve the Eulerian/co-Eulerian type and that maps between Eulerian QRGs must be graph coverings (Theorem~\ref{thm:binary-QRG-map}). This is illustrated by a quiver of all maps between connected star-preserving QRGs on 2--6 vertices in Figures~\ref{fig:QRG-quiver-co-eulerian} and~\ref{fig:QRG-quiver-eulerian}. This is in the same style as our earlier quiver and classification of all Hopf algebras over $\F_2$ to dimension 4 in \cite{MaPac:cla}.  Our results in the paper also make contact with \cite{Ma:boo} which provided a power-set picture of this case and a de Morgan type duality. Of relevance is that a finite graph calculus over any field is inner with an element $\theta$ given by the sum of all arrows. Over $\F_2$, any 1-form $\omega$ corresponds to a subset of arrows and has a natural de Morgan dual $\theta+\omega$.  By working over $\F_2$ we lose immediate contact with quantum gravity, where it is widely believed that some form of noncommutative geometry could be a better model of geometry at the Planck scale where quantum gravity effects are strong, but on the plus side we are able to consider all QRGs together which would be a first step towards a full quantum gravity theory, in which we should not just integrate over all possible metrics on a fixed manifold but also sum over all topologies. For an action we need a Ricci curvature which, in the current `working definition' involves a choice of lift $i:\Omega^2_{min}\to \Omega^1\tens_A\Omega^1$. Over $\F_2$, the canonical `averaging' choice as in \cite{Ma:tre,Ma:alg4} already fails on return blocks at vertices of even degree, but the space of bimodule lifts is always nonempty. Some consideration of these matters as well as a candidate Einstein tensor specific to $\F_2$ are considered in Appendix~\ref{app:phys}.  Meanwhile, $\F_2$ or `digital' QRGs have the merit that they could be hard-wired as digital circuits as previously proposed in \cite{MaPac:dig} in a classification of parallelisable QRGs over $\F_2$ (not just graph ones) for coordinate algebras of dimension $n\leq3$. The construction of $\nabla$ as a gate and its curvature as a digital process are developed in Appendix~\ref{app:digital} and illustrated on triangle $C_3$.

Appendix~\ref{app:computational-classification} overviews the python code used to study the larger $n$ behaviour. Counts of QRGs through $n=13$ make use of our theoretical results, reducing the problem to properties of graphs. Both theory and coding were checked against an independent direct enumeration of all graph QRGs over $\F_2$ through $n=9$.  To put this in context, there are about 50 trillion connected $n=13$ graphs, about 13 billion of them admit QRGs and only one has a flat QLC. Finally, although not in our present scope, by extending the methods of the paper to $\F_p$ for other primes $p$ we may eventually get insight through their common solutions to some QRGs over $\C$ (these will not be all, as \cite{Ma:tre} already shows due to the need for metric ratios that are often irrational).

\subsection*{Acknowledgements} This work was supported by a Leverhulme Trust project grant RPG-2024-177 and by access provided by OpenAI through its ChatGPT for Academic Researchers programme.

\noindent{\bf AI declaration:}  ChatGPT(6 astra) was used as an extremely competent research assistant for the preparation of proofs and for advanced python coding. The author takes final responsibility for all proofs and reported calculations. 

\noindent{\bf Data declaration:} All data created about graphs up to 13 vertices and the associated Python codes are available at  \cite{github}.

\noindent{\bf Conflicts declaration:} The author is not aware of any conflicts of interest.

\section{Preliminaries on graph QRGs}
% ============================================================
Here we collect some basic elements of QRG and its application to graphs\cite{BegMa,Ma:gra}. We also  summarise some elements of the graph theory  for reference. 

\subsection{The graph calculus}
All graphs \(G=(V,E)\) are taken to be finite and simple: they have no loops or multiple
edges (at most one arrow for each ordered pair in the directed case).
Undirected graphs are also viewed as bidirected, replacing each edge
\(\{x,y\}\) in the edge set $E$ by the two arrows \(x\to y\) and \(y\to x\), abbreviated to
$xy$ and $yx$ respectively. The order of a graph is $|G|=|V|$. We let \(A=k(V)\)
be the algebra of functions on $V$ with pointwise operations,
where \(k\) is a field. 

We use \(N(x)\) for the neighbours of \(x\in V\), and \(\delta_x\) for
the characteristic function of \(\{x\}\). The space \(\Omega^1\)
has a \(k\)-basis \(\omega_{xy}\), with
\[
 f\omega_{xy}=f(x)\omega_{xy},\quad
 \omega_{xy}f=f(y)\omega_{xy},\quad
 \extd f=\sum_{x\to y}(f(y)-f(x))\omega_{xy}.
\]
Sums over directed edges include both orientations. Put
\(\theta=\sum_{x\to y}\omega_{xy}\), so \(\extd f=[\theta,f]\).
All tensor products of forms are balanced over \(A\):
\(\omega f\ten\eta=\omega\ten f\eta\). Thus a tensor of two edge
forms vanishes unless the edges compose. The composable paths
\(\omega_{xyz}:=\omega_{xy}\ten\omega_{yz}\) form a basis of $\Omega^1\tens_A\Omega^1$, which decomposes into \emph{endpoint blocks}
\[
 \Omega^1\ten\Omega^1=\bigoplus_{x,z\in V}P_{x,z},\qquad
 P_{x,z}=\Span_k\{\omega_{xyz}:
                         x\to y\to z\}.
\]
For a bidirected graph their dimensions are $m_{x,z}=|N(x)\cap N(z)|$, with $m_{x,x}=\deg(x)$. 

For the calculus $\Omega_{min}$, the wedge map \(\wedge:\Omega^1\ten\Omega^1\to\Omega^2_{\min}\)
is the quotient imposing only the relations
\[
\sum_{y:x\to y\to z}\omega_{xy}\wedge\omega_{yz}=0
\]
for every pair of endpoints, including \(x=z\). We write
\[
 \operatorname{red}_2:\Omega^1\ten\Omega^1\longrightarrow\Omega^2_{\min},
 \qquad \omega\ten\eta\longmapsto\omega\wedge\eta
\]
for this quotient map. In a nonempty
endpoint block its kernel is spanned by the single path sum
\[
 P^{\mathrm{sum}}_{x,z}
 =\sum_{y:x\to y\to z}\omega_{x,y}\ten\omega_{y,z}.
\]
Thus reduction identifies coefficient vectors differing by a scalar multiple of the
all-ones vector. The quotient of a nonempty block has dimension \(m_{x,z}-1\).
A block containing just one path vanishes after reduction.

Higher forms are
obtained by the differential ideal generated by the degree two relations.
On one-forms,
\(\extd\omega=\theta\wedge\omega+\omega\wedge\theta\).
For a bidirected graph over a field with involution, \(\omega_{xy}^*=-\omega_{yx}\),
and \(\dagger\) reverses tensor order and applies star to both factors.
Over \(\F_2\), the scalar involution is the identity.

\subsection{Metrics and bimodule connections}
A central nondegenerate metric has the form
\[
 g=\sum_{x\to y}g_{xy}\omega_{xy}\ten\omega_{yx},\qquad
 g_{xy}\in k^\times.
\]
Its inverse pairing is
\((\omega_{xy},\omega_{yx})=g_{yx}^{-1}\delta_x\), and is zero
on all other composable pairs.
We call the choice $g_{xy}=1$ on every arrow the \emph{Euclidean metric}, which over $\F_2$ is the unique metric.
A bimodule connection \((\nabla,\sigma)\) satisfies
\[
 \nabla(f\omega)=\extd f\ten\omega+f\nabla\omega,\qquad
 \nabla(\omega f)=(\nabla\omega)f+\sigma(\omega\ten\extd f),
\]
where \(\nabla:\Omega^1\to\Omega^1\ten\Omega^1\) is \(k\)-linear
and \(\sigma:\Omega^1\ten\Omega^1\to\Omega^1\ten\Omega^1\)
is an \(A\)-bimodule map. It has the inner presentation\cite{Ma:gra}
\[
 \nabla\omega=\theta\ten\omega-\sigma(\omega\ten\theta)+\alpha(\omega),
\]
with \(\alpha\) an \(A\)-bimodule map. We call the connection inner
when \(\alpha=0\). Its extension to tensors is
\(\nabla=\nabla\otimes_A\id+(\sigma\otimes_A\id)(\id\otimes_A\nabla)\). Bimodule maps preserve both endpoints, so \(\sigma\) acts blockwise, with restriction \(\sigma_{x,z}\).

A quantum Levi--Civita connection (QLC) is torsion-free and
metric-compatible:
\[
 \wedge\nabla=\extd,\qquad \nabla g=0.
\]
Metric compatibility is an equality in the unreduced three-tensor
space, whose basis consists of paths \(x\to y\to z\to w\).
The terms independent of \(\alpha\) have distinct overall endpoints,
while the terms involving \(\alpha\) have equal overall endpoints.
The metric equations therefore separate into equations for \(\sigma\)
and, for each solution, homogeneous linear equations for \(\alpha\). Similarly, in the inner presentation, torsion-freeness is
\begin{equation}\label{eq:tor} \wedge(\id+\sigma)=0,\quad \wedge\alpha=0\end{equation}
while $*$-preservation requires
\begin{equation}\label{eq:star} 
(\dagger\sigma)^2=\id,\quad \sigma\dagger\alpha=\alpha*.\end{equation}
A connection is $*$-compatible if we only require the first of these. Then \(\sigma\) is invertible, with 
\(\sigma^{-1}=\dagger\sigma\dagger\). The curvature is
\[
 R_\nabla=(\extd\otimes_A\id-\id\wedge\nabla)\nabla:
 \Om^1\longrightarrow\Om^2_{\min}\ten\Om^1.
\]
A connection is flat when \(R_\nabla=0\). By a graph quantum Riemannian geometry (QRG) we mean a graph with its calculus, a metric and a QLC. Star preservation will be specified separately.

%For curvature, the first two tensor legs are reduced: for each fixed
%\(x,z,w\), the relation is
%\(\sum_{y\in N(x)\cap N(z)}
%\omega_{xy}\wedge\omega_{yz}\ten\omega_{zw}=0\).

\subsection{Degree and codegree parity}
\label{ssec:graph-parity}
For distinct vertices \(x,y\), define their \emph{codegree} by
\[
 \operatorname{codeg}(x,y)=|N(x)\cap N(y)|.
\]
A graph is \emph{Eulerian}, or \emph{even-degree}, if every vertex
has even degree. Connectedness is not part of this convention.
A graph is \emph{odd-degree} if every vertex has odd degree.

We call a graph \emph{co-Eulerian} if every pair of distinct vertices
has even codegree, including zero. If \(A\) is its adjacency matrix over \(\F_2\), then
\[
 (A^2)_{xy}=\operatorname{codeg}(x,y)\pmod2\quad(x\ne y),
 \qquad (A^2)_{xx}=\deg(x)\pmod2.
\]
Thus an even-degree graph is co-Eulerian precisely when \(A^2=0\),
and an odd-degree graph is co-Eulerian precisely when \(A^2=I\).
For connected co-Eulerian graphs these are the only possibilities,
as proved in Proposition~\ref{prop:complement}.

\subsection{Cayley graphs}
Let \(H\) be a finite group and let \(\CC\subset H\setminus\{e\}\)
be an inverse-closed generating set. The Cayley graph has edges
\(x\to xs\), \(s\in \CC\). Its displacement one-forms are
\[
 e_s=\sum_{x\in H}\omega_{x,xs},\qquad
 e_s f=R_s(f)e_s,\qquad R_s(f)(x)=f(xs).
\]
They form a left \(k(H)\)-module basis of \(\Omega^1\), with
\(\theta=\sum_{s\in \CC}e_s\),
\(\extd f=\sum_{s\in \CC}(R_s(f)-f)e_s\), and, when a star structure
is specified, \(e_s^*=-e_{s^{-1}}\). The Euclidean metric is
\[
 g=\sum_{s\in \CC}e_s\otimes_A e_{s^{-1}}.
\]
The path-sum relations take the form
\[
 \sum_{\substack{s,t\in \CC\\st=h}}e_s\wedge e_t=0
 \qquad(h\in H).
\]
These relations are left-module identities: multiplication by
\(\delta_x\) recovers the relation for the endpoints \((x,xh)\). A connection is invariant under left translations precisely when its coefficients in this basis lie in $k$.

For presentation purposes it is convenient if we can take \(H=\Z_n\). We then need an inverse-closed set
\(\CC=-\CC\subset\Z_n\setminus\{0\}\) generating \(\Z_n\).
Indices are then added modulo \(n\), and
\[
 e_i=\sum_x\omega_{x,x+i},\qquad
 e_i f=R_i(f)e_i,\qquad
 \sum_{\substack{i,j\in \CC\\i+j=r}}e_i\wedge e_j=0.
\]
Thus \(e_i\otimes_Ae_j\) has displacement \(i+j\).

\begin{table}[tbp]
\centering
\small
\setlength{\tabcolsep}{4pt}
\begin{tabular}{r r r r r r r}
\toprule
\(n\)&\shortstack{Connected\\graphs}&Eulerian&Co-Eulerian&Both&Cayley&\shortstack{Cyclic\\Cayley}\\
\midrule
2&1&0&1&0&1&1\\
3&2&1&0&0&1&1\\
4&6&1&2&1&2&2\\
5&21&4&0&0&2&2\\
6&112&8&4&2&5&5\\
7&853&37&0&0&3&3\\
8&11117&184&14&7&10&8\\
9&261080&1782&1&1&7&6\\
10&11716571&31026&49&24&16&16\\
11&1006700565&1148626&7&7&7&7\\
12&164059830476&86539128&325&159&64&38\\
% Order 13: connected/Eulerian counts from OEIS A001349/A003049;
% Cayley counts from qlc-cayley-counts-through13.json.
% Co-Eulerian/intersection counts from qlc-eulerian-counts13-user-result.json;
% at odd order every connected co-Eulerian graph is Eulerian.
13&50335907869219&12798435868&108&108&13&13\\
\bottomrule
\end{tabular}
\caption{Connected graph isomorphism classes. Co-Eulerian and intersection counts come from 
Python enumeration through \(n=13\).}
\label{tab:graph-classes}
\end{table}

The summary Table~\ref{tab:graph-classes} counts the discussed graph classes up to
isomorphism, with the later columns obtained by Python-based enumeration. Cayley presentations over different groups count once
when their underlying graphs are isomorphic. The only Cayley graphs through $n=13$ that need a non-Abelian group are 18 of them at $n=12$.

\section{Classification of QRGs over \(\mathbb F_2\)}
\label{sec:classification}

This section contains the main theoretical results of the paper. Graphs are taken to be finite and bidirected with their unique metric over \(\F_2\) and we use \(\Omega_{\min}\) for the exterior algebra.

\subsection{QLC equations over $\F_2$}

To find QLCs,  we first parametrise torsion-free bimodule connections, then impose metric
compatibility. Invertibility and star compatibility of \(\sigma\)
will follow from these equations.

Torsion-freeness requires every column of \(\id+\sigma_{x,z}\)
to lie in \(\ker\operatorname{red}_2\), giving
\begin{equation}
 \sigma_{x,z}=I_{m_{x,z}}+\boldsymbol1_{x,z}\lambda_{x,z}^T,
 \qquad \lambda_{x,z}\in\F_2^{m_{x,z}}.
 \label{eq:F2-sigma-block}
\end{equation}
Choose an ordering of the intermediate vertices in each nonempty
block. Here \(\boldsymbol1_{x,z}\) is the all-ones column in this
two-path basis, and columns of \(\sigma_{x,z}\) are images of basis paths. For a basis path \(p_j\), this means
\(\sigma(p_j)=p_j+\lambda_jP^{\mathrm{sum}}_{x,z}\).
The bimodule property gives \(\alpha(\omega_{xy})\in P_{x,y}\),
and \(\operatorname{red}_2\alpha=0\) therefore gives
\begin{equation}
 \alpha(\omega_{xy})=a_{xy}P^{\mathrm{sum}}_{x,y},
 \qquad a_{xy}\in\F_2,
 \label{eq:F2-alpha-block}
\end{equation}
when \(P_{x,y}\ne0\), and zero otherwise. Such a two-path exists
precisely when the edge \(xy\) belongs to a triangle.

\subsection{$*$-compatibility over $\F_2$}
\label{ssec:star}

Over \(\F_2\), \(\omega_{xy}{}^*=\omega_{yx}\), and its order-reversing extension to two-tensors is
\(\dagger(\omega_{xyz})=\omega_{zyx}\).

\begin{proposition}
\label{prop:det-equiv}
For a torsion-compatible map on any finite graph
with the stated calculus over \(\F_2\), the following are equivalent:
\begin{enumerate}
\item \(\sigma\) is invertible;
\item \(\sigma^2=\id\);
\item \(\sum_j(\lambda_{x,z})_j=0\) for every nonempty endpoint block.
\end{enumerate}
In this case $\sigma$ is $*$-compatible precisely when
\(\lambda_{x,z}=\lambda_{z,x}\), identifying intermediate vertices
by path reversal.
\end{proposition}

\begin{proof}
Write \(S=I+\boldsymbol1\lambda^T\) on a nonempty block. Then
\[
\det S=1+\lambda^T\boldsymbol1,\qquad
S^2=I+(\lambda^T\boldsymbol1)\boldsymbol1\lambda^T.
\]
If the sum is zero, both assertions follow. If it is one, \(\lambda\)
is nonzero, so \(S\) is singular and \(S^2\ne I\).
Under these equivalent conditions, \(\sigma^{-1}=\sigma\), so
\((\dagger\sigma)^2=\id\) is equivalent to
\(\dagger\sigma\dagger=\sigma\). Path reversal exchanges the blocks
\((x,z)\) and \((z,x)\) without changing their middle vertices,
giving precisely \(\lambda_{x,z}=\lambda_{z,x}\).
\end{proof}

\subsection{Moduli of QLCs}

We first classify $\sigma$, then the compatible $\alpha$. Integer counts in scalar expressions are reduced modulo two.

\subsubsection{Moduli of braidings}
\label{ssec:general-F2-classification}

\begin{theorem}[General braiding classification and existence]
\label{thm:general-F2-classification}
Let \(G\) be a connected finite graph with at least two
vertices and 
\(\Omega_{\min}\) calculus. For \(q\in\F_2\), define
\begin{equation}
\label{eq:general-F2-sigma}
\sigma_q(\omega_{x,y}\otimes_A\omega_{y,z})
=\omega_{x,y}\otimes_A\omega_{y,z}
+(q+\delta_{x,z})P^{\mathrm{sum}}_{x,z}.
\end{equation}
\begin{enumerate}[label=\textup{(\roman*)}]
\item \(G\) admits a QLC braiding \(\sigma_0\) if and only if it is Eulerian.
\item \(G\) admits a QLC braiding \(\sigma_1\) if and only if it is co-Eulerian.
\end{enumerate}
These are the only QLC braidings, each realised by the inner
connection \((\sigma_q,0)\).
\end{theorem}

\begin{proof}
Torsion-freeness gives a parameter \(L_{xyz}\in\F_2\) on each
composable two-path, with
\[
\sigma(\omega_{x,y}\otimes_A\omega_{y,z})
=\omega_{x,y}\otimes_A\omega_{y,z}
+L_{xyz}P^{\mathrm{sum}}_{x,z}.
\]
The metric equation separates into its inner and \(\alpha\) parts:
\[
\theta\otimes_Ag+\sigma_{12}\sigma_{23}(g\otimes_A\theta)=0,
\qquad
(\alpha\otimes_A\id)g+\sigma_{12}(\id\otimes_A\alpha)g=0.
\]
Here \(\sigma_{12}=\sigma\ten\id\) and
\(\sigma_{23}=\id\ten\sigma\), acting on three-tensors. Set
\(H_a=\sum_{j\in N(a)}L_{aja}\).
For a three-path \((a,b,c,d)\) with \(d\ne a\) and \(a\sim d\) (i.e. forming an edge),
the coefficient in the first expression is
\begin{equation}
\label{eq:general-F2-metric-component}
0=\delta_{b,d}+\delta_{a,c}(1+H_a)+L_{bad}
 +\sum_{j\in N(a)\cap N(c)}L_{jad}L_{ajc}.
\end{equation}
If \(a\not\sim d\), that coefficient is automatically zero.
Equation~\eqref{eq:general-F2-metric-component} follows by applying
\(\sigma_{23}\), then \(\sigma_{12}\), to the return-path metric
in \(g\otimes_A\theta\); the intermediate vertex \(j\) must lie
in \(N(a)\cap N(c)\).

Take \(c=a\), fix an edge \(a\sim d\), and compare all choices
\(b\in N(a)\). The quantity \(L_{bad}+\delta_{b,d}\) is independent
of \(b\). Hence there is a value \(q_{ad}\) on each directed edge,
such that, on every two-path,
\[
L_{xyz}=q_{yz}+\delta_{x,z}.
\]
Substitution into the same equations gives
\begin{equation}
\label{eq:general-F2-edge-values}
q_{da}=H_a+(1+H_a)q_{ad},\qquad
H_a=\sum_{j\in N(a)}(q_{ja}+1).
\end{equation}
If \(H_a=1\), the first equation forces \(q_{da}=1\) for all
\(d\in N(a)\), making the second equation give \(H_a=0\), a
contradiction. Consequently \(H_a=0\) for every vertex,
\(q_{ad}=q_{da}\), and
\[
\sum_{j\in N(a)}q_{ja}=\deg(a)\pmod2.
\]

For \(a\ne c\) with a common neighbour, put
\(T_{ac}=\sum_{j\in N(a)\cap N(c)}q_{jc}\).
For every \(d\in N(a)\cap N(c)\), choose a three-path
\((a,b,c,d)\) and use \eqref{eq:general-F2-metric-component}.
It reduces to
\[
q_{dc}=q_{ad}(1+T_{ac}).
\]
If \(T_{ac}=1\), all \(q_{dc}\) in its defining sum vanish,
a contradiction. Therefore \(T_{ac}=0\), and \(q_{dc}=q_{ad}\)
for every common neighbour \(d\). Thus any two distinct edges
incident at \(d\) have the same value. Connectedness propagates this
equality along edges: every \(q_{ad}\) equals a single \(q\).
The conditions already proved now read
\[
(q+1)\deg(a)=0,\qquad q\,m_{a,c}=0\quad(a\ne c).
\]
These are exactly the two parity criteria in the theorem.

Conversely, suppose one of those criteria holds and use
\(L_{xyz}=q+\delta_{x,z}\). Then \(H_a=0\), and the coefficient
\eqref{eq:general-F2-metric-component} vanishes: for \(c=a\), use
\eqref{eq:general-F2-edge-values}; for \(c\ne a\), use
\(T_{ac}=q\,m_{a,c}=0\). All other inner coefficients are zero
by their endpoints. Taking \(\alpha=0\) proves metric compatibility;
torsion-freeness follows from the defining path-sum relations.
If both braidings occur, every vertex has positive even degree,
so their return blocks \(I+J\) and \(I\) are different.
\end{proof}

We call QRG where the QLC has braiding \(\sigma_0\) an \emph{Eulerian geometry},
and one with braiding \(\sigma_1\) a \emph{co-Eulerian geometry}.
These names distinguish connections: a graph that is both Eulerian
and co-Eulerian supports both types of geometry.

\begin{corollary}[Automatic invertibility and the first star condition]
\label{conj:automatic-invertibility}
Every QLC on a finite graph over \(\F_2\)
is $*$-compatible and has \(\sigma^2=\id\). 
\end{corollary}

\begin{proof}
On a connected component with an edge, the block torsion vector is
\(\lambda_{x,z}=(q+\delta_{x,z})\boldsymbol1\).
Its sum vanishes: on return blocks this is
\((q+1)\deg(x)=0\); on non-return blocks it is \(q\,m_{x,z}=0\).
Proposition~\ref{prop:det-equiv} gives \(\sigma^2=\id\).
Path reversal identifies equal vectors \(\lambda_{x,z}\) and
\(\lambda_{z,x}\), proving the first star condition as well.
For a disconnected graph, all spaces and equations split by connected
components. Isolated vertices contribute no one-forms, so the
conclusion holds there vacuously.
\end{proof}

\subsubsection{Moduli of \(\alpha\)}
\label{ssec:general-alpha-parameters}

Let \(G\) be a finite graph and write
\(E_\triangle^{\mathrm{dir}}\) for its directed edges in triangles.
For coefficients \(a_{xy}\in\F_2\), define the bimodule map
\[
 \alpha_a(\omega_{xy})=a_{xy}P^{\mathrm{sum}}_{xy},
 \qquad (x,y)\in E_\triangle^{\mathrm{dir}},
\]
and set it to zero on all other edges. These are precisely the
bimodule maps satisfying \(\wedge\alpha=0\).

\begin{definition}
\label{def:alpha-spaces}
Let \(D(G)\) and \(D_*(G)\) be the following undirected graphs:
\begin{enumerate}
\item The vertices of \(D(G)\) are \(E_\triangle^{\mathrm{dir}}\).
For every ordered triangle \((x,u,w)\), join \((x,w)\) to \((u,x)\).
\item The vertices of \(D_*(G)\) are the unoriented edges in triangles.
Join two whenever they belong to a common triangle.
\end{enumerate}
Define \(V(G)\) to be the space of maps \(\alpha_a\) whose
coefficients are constant on each connected component of \(D(G)\).
Define \(V_*(G)\) similarly using \(D_*(G)\), assigning the
same coefficient to both orientations of an edge. Thus
\[
 V_*(G)\subseteq V(G),\qquad
 V(G)\cong\F_2^d,\qquad V_*(G)\cong\F_2^{d^*},
\]
where \(d\) and \(d^*\) count the respective connected components.
A basis map assigns one to the edges in one component and zero elsewhere.
\end{definition}

The graph \(D_*(G)\) is the anti-Gallai graph of \(G\), with its
isolated vertices removed \cite{Le}. Thus \(d^*(G)\) counts the
components of edges linked through triangles, excluding edges in no triangle.

\begin{theorem}[Classification of the \(\alpha\)-families]
\label{thm:general-alpha-kernel}
Let \(G\) be connected. For either braiding \(\sigma_q\) that occurs
in Theorem~\ref{thm:general-F2-classification}, all QLCs with that
braiding are precisely
\[
 (\sigma_q,\alpha),\qquad \alpha\in V(G).
\]
They are star-preserving precisely when \(\alpha\in V_*(G)\).
In particular, the same spaces of maps apply to both braidings whenever
both occur.
\end{theorem}

\begin{proof}
The remaining metric equation is homogeneous and linear in \(\alpha\).
Its only possible three-paths have equal overall endpoints and hence
are \((x,u,w,x)\), with \((x,u,w)\) a triangle. The coefficient
from \((\alpha\ten\id)g\) is \(a_{xw}\), while
\(\sigma_{12}(\id\ten\alpha)g\) contributes
\(a_{ux}+q\sum_{y\in N(x)\cap N(w)}a_{yx}\). Thus the equations are
\begin{equation}
\label{eq:general-alpha-matrix}
 a_{xw}+a_{ux}+qT_{xw}=0,\qquad
 T_{xw}=\sum_{y\in N(x)\cap N(w)}a_{yx}.
\end{equation}
For \(q=0\), these are exactly the component equalities defining
\(V(G)\). For \(q=1\), the number \(m_{x,w}\) of common
neighbours is even. Summing over \(u\in N(x)\cap N(w)\) therefore
gives \(T_{xw}=0\), leaving the same equalities. Conversely, those
equalities imply \(T_{xw}=m_{x,w}a_{xw}=0\), so they suffice.

For an edge in a triangle,
\[
 \sigma_q(P^{\mathrm{sum}}_{xy})
 =(1+q\,m_{x,y})P^{\mathrm{sum}}_{xy}=P^{\mathrm{sum}}_{xy},
 \qquad \dagger P^{\mathrm{sum}}_{xy}=P^{\mathrm{sum}}_{yx}.
\]
The second star condition is consequently \(a_{xy}=a_{yx}\).
Identifying reversed edges in the triangle equations gives exactly
\(V_*(G)\). The first star condition holds by
Corollary~\ref{conj:automatic-invertibility}.
\end{proof}

For example, on \(K_3\) the auxiliary graph has two components,
one for each cyclic orientation, giving two independent parameters and four
QLCs. Reversal exchanges the two components, so the star-preserving
family has one parameter and two connections. For two triangles sharing
exactly one vertex and no edge, the original graph has only even
degrees. Its auxiliary graph has four components, giving four
independent parameters, sixteen QLCs, and four star-preserving connections.

\begin{remark}[Counts, symmetry and curvature]
For a connected graph admitting QLCs, the full family is a disjoint
union of its Eulerian and co-Eulerian families, when they occur;
both have the same \(\alpha\)-space.
Graph automorphisms preserve each \(\sigma_q\) and act
on this space by permuting the auxiliary components of
Definition~\ref{def:alpha-spaces}.
The flat subfamily need not be a vector space: on \(K_3\), for
example, three of the four \(\alpha\)-values are flat. Curvature can
contain quadratic terms in the independent parameters.
\end{remark}

\subsection{Triangle components and derived geometries}
\label{ssec:derived}

The undirected graphs \(D(G)\) and \(D_*(G)\) need not be connected.
We equip each with its own \(\Omega_{\min}\) calculus and unique
\(\F_2\)-metric.

\begin{proposition}
\label{prop:derived}
If an edge \(e\) of \(G\) belongs to \(t_e\) triangles, its
corresponding vertices in \(D(G)\) and \(D_*(G)\) have degree
\(2t_e\). Hence both derived graphs are Eulerian, and every connected
component admits the inner star-preserving Eulerian geometry.
No component of \(D(G)\) is co-Eulerian.
\end{proposition}

\begin{proof}
Each triangle containing \(e\) contributes two distinct neighbours
in either derived graph. Neighbours contributed by different triangles
are distinct, proving the degree formula. Apply
Theorem~\ref{thm:general-F2-classification} and
Corollary~\ref{conj:automatic-invertibility} componentwise.
Adjacent vertices \(xy,yz\) of \(D(G)\) have exactly one common
neighbour, \(zx\), so the \(\sigma_1\) criterion fails.
\end{proof}

\begin{theorem}[The derived double cover]
\label{thm:derived-cover}
Forgetting orientation defines a two-sheeted graph covering
\[
 \pi:D(G)\longrightarrow D_*(G),\qquad xy\longmapsto\{x,y\},
\]
with involution \(\tau(xy)=yx\), acting on coefficients by
\((\tau a)_{xy}=a_{yx}\). Let \(s(G)\) count the connected
components of \(D_*(G)\) whose inverse images under \(\pi\) are
disconnected. Then
\[
 d=d^*+s(G),\qquad
 V_*(G)=\{\alpha\in V(G):\tau\alpha=\alpha\}.
\]
An inverse image splits precisely when its original edges can be
oriented so that every triangle is a directed cycle. If \(G\) is
connected and admits QLCs, the fraction which are star-preserving is
\(2^{-s(G)}\), before quotienting by graph automorphisms.
\end{theorem}

\begin{proof}
A graph covering is a graph map inducing a bijection on every neighbour
set. The neighbours of \(xy\) in \(D(G)\) are \(yz\) and \(zx\),
for \(z\in N(x)\cap N(y)\). These map bijectively to the neighbours
\(\{y,z\}\) and \(\{z,x\}\) of \(\{x,y\}\). Each vertex has exactly
two lifts, exchanged by \(\tau\).

Over a connected component downstairs, a two-sheeted cover is either
connected or consists of two copies of that component: lifting paths
from a fixed vertex reaches either both sheets or only one. Thus each
component contributes one to \(d^*\) and either one or two to \(d\).
Reversal fixes a connected inverse image and exchanges the two
components of a split inverse image. Invariance of its component
coefficients is exactly \(a_{xy}=a_{yx}\), proving the assertion
about \(V_*(G)\).

Choosing one component of a split inverse image chooses one orientation
of every original edge. The triangle downstairs associated to a triangle of \(G\) lifts to
its two oriented cycles, so the chosen orientations make each original
triangle cyclic.
Conversely, such orientations select a union of components containing
exactly one lift of every vertex, and hence split the cover.
For each braiding there are \(2^d\) QLCs and \(2^{d^*}\)
star-preserving ones; their ratio is \(2^{-s(G)}\).
\end{proof}

With the inner Eulerian geometry on every component, \(\pi\)
is also a map of QRGs in the sense of Section~\ref{sec:QRG-maps}:
the covering condition and \(\alpha=0\) satisfy
Theorem~\ref{thm:binary-QRG-map} componentwise.
For example,
\[
 D(K_3)=K_3\sqcup K_3,\qquad D_*(K_3)=K_3,
\]
so \(s(K_3)=1\) and half its QLCs are star-preserving.
For \(K_4\), the graph \(D_*(K_4)\) is octahedral and \(D(K_4)\)
is connected. Indeed, an orientation of \(K_4\) has a vertex with
two outgoing edges, so cannot make every triangle cyclic.
Thus \(s(K_4)=0\) and all its QLCs are star-preserving.
More generally, if \(D_*(G)\) is connected and \(G\) contains a
\(K_4\), its cover cannot split, so \(d=d^*=1\). Each braiding
that occurs then gives exactly two QLCs, both star-preserving.

\medskip\paragraph{\bf Open problem: anti-Gallai components of Eulerian graphs.}
For a uniformly chosen connected Eulerian graph isomorphism class of
order \(n\), determine the tail of the distribution of \(d^*(G)\).
Theorem~\ref{thm:typical-eulerian} below proves \(d=d^*=1\) with
probability tending to one. Does the stronger assertion
\[
 \mathbb E\bigl[2^{d^*(G)}\bigr]\longrightarrow2
\]
hold? This requires control of the exponentially weighted exceptional
components. For the full QLC count one must also control \(s(G)\),
since \(d(G)=d^*(G)+s(G)\).
Section~\ref{ssec:graph-counts} formulates the corresponding counting
conjecture and gives the current evidence.

The derived graphs can be iterated, taking the inner Eulerian
connection on every component at each stage. Iteration through
\(D_*(G)\) retains information about higher cliques: edges linked
through triangles are followed by triangles linked through tetrahedra.
Write \(k_r(G)\) for the number of copies of \(K_r\) in \(G\),
\(t(G)=k_3(G)\), and \(D_*^j\) for \(j\) successive applications
of \(D_*\). Two triangles are called equivalent if they can be
joined by a sequence of triangles in which consecutive triangles
belong to a common \(K_4\).

\begin{proposition}[Iteration of the derived graphs]
\label{prop:derived-iteration}
Let \(G\) be a finite graph.
\begin{enumerate}[label=\textup{(\roman*)}]
\item
\[
 D_*(D(G))\cong\bigsqcup_{1}^{2t(G)}K_3,
 \qquad D(D(G))\cong\bigsqcup_{1}^{4t(G)}K_3.
\]
\item The components of \(D_*^2(G)\) correspond to the equivalence
classes of triangles of \(G\). If \(G\) contains no \(K_4\), then
\[
 D_*^2(G)\cong\bigsqcup_{1}^{t(G)}K_3,
 \qquad D(D_*(G))\cong\bigsqcup_{1}^{2t(G)}K_3.
\]
\item The clique counts satisfy
\[
 t(D_*(G))=t(G)+4k_4(G),\qquad
 k_r(D_*(G))=(r+1)k_{r+1}(G)\quad(r\geq4).
\]
\item If \(G\) has clique number \(w\geq3\), then
\[
 D_*^{w-1}(G)\cong\bigsqcup_{1}^{N(G)}K_3,
 \qquad N(G)=\sum_{r=3}^{w}\frac{r!}{6}k_r(G).
\]
Further applications of \(D_*\) leave this graph unchanged up to
isomorphism. If \(G\) is triangle-free, both derived graphs are empty.
\end{enumerate}
\end{proposition}

\begin{proof}
Each edge of \(D(G)\) belongs to exactly one triangle, arising
from an oriented triangle of \(G\). These \(2t(G)\) triangles
have disjoint edge sets. Applying \(D_*\) separates them into
copies of \(K_3\); applying \(D\) gives two oriented copies
of each, proving (i).

A triangle in \(D_*(G)\) is a set of three pairwise incident
edges of \(G\). They either form a triangle of \(G\), or share
a vertex. In the latter case their other endpoints must be pairwise
adjacent, giving a \(K_4\) with a distinguished vertex. These
are the two types of triangle, proving the first formula in (iii).
Four or more pairwise incident edges must share a vertex, so an
\(r\)-clique in \(D_*(G)\), \(r\geq4\), is exactly an
\((r+1)\)-clique of \(G\) with a distinguished vertex. This
proves the remaining formula in (iii).

Every edge of \(D_*(G)\) is a pair of edges of a unique triangle
\(T\) of \(G\), meeting at a distinguished vertex of \(T\).
Thus the vertices of \(D_*^2(G)\) are the three incidences
\((T,v)\) for each triangle \(T\). The first type of triangle
above joins the three incidences of \(T\); the second joins
incidences from faces of a \(K_4\) at a common vertex. Within
each \(K_4\), these joins connect all four faces. They generate
exactly the equivalence relation in (ii). When there is no \(K_4\),
each edge of \(D_*(G)\) lies in its unique triangle of the first
type. The same argument as in (i) then gives both decompositions
in (ii).

By (iii), applying \(D_*\) reduces a clique number greater than
three by one. After \(w-3\) iterations the graph contains no
\(K_4\), so two more applications give a disjoint union of
triangles by (ii). Their number is \(t(D_*^{w-3}(G))\), since
the triangle count is unchanged once there is no \(K_4\).
For \(0\leq j\leq w-4\), repeated use of (iii) gives
\[
 k_4(D_*^j(G))=\frac{(j+4)!}{4!}k_{j+4}(G).
\]
Summing the increments \(4k_4(D_*^j(G))\) of the triangle count
proves the formula for \(N(G)\). Finally \(D_*(K_3)=K_3\),
and \(D_*\) preserves disjoint unions.
\end{proof}

Theorem~\ref{thm:derived-cover} applied to \(D_*(G)\) gives the
two-sheeted cover
\[
 D(D_*(G))\longrightarrow D_*^2(G).
\]
The components downstairs parametrise the star-preserving
\(\alpha\)-moduli of the derived geometry; splitting of this cover
measures the additional parameters without the star condition.

For example, successive applications of \(D_*\) give
\[
 K_4\longmapsto\text{the octahedral graph}
 \longmapsto\text{the cuboctahedral graph}
 \longmapsto8K_3.
\]
The second derived graph has vertices the vertex--face incidences of a
tetrahedron, joined when they share the vertex or the face. It is
the line graph of the tetrahedron's bipartite vertex--face incidence
graph, which is \(K_{4,4}\) minus a perfect matching, hence a cube. Moreover,
\(D(D_*(K_4))\) is two copies of the cuboctahedral graph: the
octahedral graph is \(K_{2,2,2}\), and orienting its edges cyclically
between the three parts makes every triangle cyclic, splitting the
cover by Theorem~\ref{thm:derived-cover}.

\subsection{Graph complementation and de Morgan duality}
\label{ssec:complement}

Write \(\overline G\) for the complement of a graph \(G\)
on the same vertex set. Its adjacency matrix over \(\F_2\) is
\(\overline A=J+I+A\), where every entry of \(J\) is one.

\begin{proposition}
\label{prop:complement}
Let \(G\) be connected and co-Eulerian. Its vertex degrees
have the same parity. If they are even, \(G\) is Eulerian and admits
both braidings; if they are odd, it admits only \(\sigma_1\).
If \(G\) has even order, \(\overline G\) is also co-Eulerian,
with opposite degree parity. When \(\overline G\) is connected,
complementation therefore exchanges these two classes.
\end{proposition}

\begin{proof}
The common-neighbour criterion says that \(A^2\) is diagonal,
with entries the degree parities. Since \(A\) commutes with \(A^2\),
adjacent vertices have the same degree parity. Connectedness therefore
gives \(A^2=dI\), where \(d\in\F_2\), and \(AJ=JA=dJ\).
A connected graph is Eulerian exactly when all its vertex degrees are
even, which is also the criterion for \(\sigma_0\). This proves the
dichotomy. If the order \(n\) is even, \(J^2=0\), so
\[
 \overline A^2=(J+I+A)^2=(d+1)I.
\]
Thus \(\overline G\) satisfies the \(\sigma_1\) criterion, with
opposite degree parity.
\end{proof}

The co-Eulerian hypothesis is essential:
odd degrees alone do not imply its common-neighbour condition.
For example, \(K_{1,3}\) has all degrees odd, but two leaves have
one common neighbour. Likewise, an Eulerian graph always admits
\(\sigma_0\), but need not admit \(\sigma_1\).

\begin{remark}
Complementation is de Morgan duality on the Boolean algebra of edge
subsets of the complete graph. There is also an operator interpretation.
When \(n\) is even, under the hypotheses of
Proposition~\ref{prop:complement}, both adjacency matrices act on
\[
 W=\{v\in\F_2^n:\textstyle\sum_xv_x=0\}/\langle\boldsymbol1\rangle.
\]
Indeed, \(A\boldsymbol1=d\boldsymbol1\) and \(A\) is symmetric;
evenness of \(n\) puts \(\boldsymbol1\) in the numerator. On \(W\),
\(J=0\), so \(\overline A=I+A\). This is the substitution
\(a\mapsto1+a\) of \cite[Lemma~5.2]{Ma:boo}, exchanging the relations
\(a^2=0\) and \(a^2=1\).

The differential-geometric de Morgan duality in
\cite[Section~3.2 and Proposition~5.3]{Ma:boo} instead complements
subsets of a fixed arrow set and transports the algebra and calculus.
Graph complementation changes the arrow set, generally changing
\(\dim\Omega^1\). The proposition concerns braiding existence;
it does not identify the full QLC families or their curvatures. \end{remark}

We can also represent every graph \(G\) on \(n\geq3\) vertices
inside the fixed calculus \(\Omega_{\min}(K_n)\). Writing
\(A_{xy}\in\F_2\) for its adjacency coefficients, put
\[
 \omega_G=\sum_{x\ne y}A_{xy}\omega_{xy},\qquad
 H_G=\sum_{x\ne y}A_{xy}\omega_{xy}\ten\omega_{yx}.
\]
Then \(\omega_G^*=\omega_G\), while \(H_G\) is central and edge-symmetric. It is immediate that
\[
 \omega_{\overline G}=\theta+\omega_G,\qquad
 H_{\overline G}=g+H_G,\quad  
 (\ ,\ )H_G=0\ \Longleftrightarrow\ G\text{ is Eulerian}.
\]
where $g$ is the unique metric and $(\ ,\ )$ is inverse. 

\begin{corollary}
\label{cor:graph-form}
Let $\sigma$ be the unique QLC braiding on $K_n$ over $\F_2$,
where $n\geq3$. For any graph $G$ on the same vertex set,

\[
 \sigma(\omega_G\ten\omega_G)=\omega_G\ten\omega_G
 \quad\Longleftrightarrow\quad
 \begin{cases}
 G\text{ is Eulerian},&n\text{ odd},\\
 G\text{ is co-Eulerian},&n\text{ even}.
 \end{cases}
\]
\end{corollary}
\begin{proof}
By Theorem~\ref{thm:general-F2-classification}, $\sigma=\sigma_0$
if $n$ is odd and $\sigma=\sigma_1$ if $n$ is even. Then

 \[
 \sigma(\omega_G\ten\omega_G)+\omega_G\ten\omega_G
 =\begin{cases}
 \displaystyle\sum_x\deg_G(x)P^{\mathrm{sum}}_{x,x},&n\text{ odd},\\[4pt]
 \displaystyle\sum_{x\ne z}\operatorname{codeg}_G(x,z)
                  P^{\mathrm{sum}}_{x,z},&n\text{ even},
 \end{cases}
\]
where the path sums are those of \(K_n\) and counts are reduced
modulo two. Indeed, the coefficient of \(\omega_{xyz}\) in
\(\omega_G\ten\omega_G\) is \(A_{xy}A_{yz}\), whose sum over
\(y\) counts degrees on return blocks and codegrees otherwise.
The displayed path sums are nonzero and lie in distinct endpoint blocks,
so their independence gives the stated equivalence.
\end{proof}

\section{Flat QLCs}
\label{ssec:flat}

Recall that \(\alpha(\omega_{xy})=a_{xy}P^{\mathrm{sum}}_{x,y}\),
with \(a_{xy}\in\F_2\), extended by zero off edges in triangles. The triangle equations give
\(a_{xy}=a_{yz}=a_{zx}\) on each oriented triangle.

\begin{theorem}
\label{thm:flat-alpha}
Let \(G\) be connected with \(n\geq2\) vertices.
\begin{enumerate}[label=\textup{(\roman*)}]
\item Every flat QLC is inner and star-preserving, except for the
      two non-inner, non-star-preserving flat QLCs on \(C_3\).
\item A flat Eulerian geometry exists if and only if
      \(G=C_n\), \(n\geq3\). It is unique for \(n\geq4\).
\end{enumerate}
\end{theorem}

\begin{proof}
Write the part of \(\nabla\omega_{xy}\) starting at \(x\) as
\[
 \sum_{v\sim x,\ z\sim v}B_{xy}(v,z)\omega_{xv}\ten\omega_{vz},
 \qquad B_{xy}(v,z)=\delta_{vy}+(q+\delta_{xz})A_{yz}
                         +a_{xy}\delta_{zy}
\]
where $A$ is the adjacency matrix of $G$. The remaining part is \(\theta\ten\omega_{xy}\). In the curvature,
terms with one incoming \(\theta\) cancel, and those with two vanish
by the path-sum relations. Thus
\begin{equation}
\label{eq:flat-curvature-coeff}
 R_\nabla(\omega_{xy})=
 \sum_{\substack{v\sim x,\ b\sim v,\ c\sim b}}
 \left(\sum_{z\sim v}B_{xy}(v,z)B_{vz}(b,c)\right)
 \omega_{xv}\wedge\omega_{vb}\ten\omega_{bc}.
\end{equation}
For fixed \(x,b,c\), a coefficient vector vanishes in the wedge
quotient precisely when it is constant on \(N(x)\cap N(b)\).

Taking \(b=x,c=y\), and using the parity and triangle equations,
the corresponding curvature component is
\[
 a_{xy}a_{yx}\sum_{v\in N(x)\cap N(y)}
       \omega_{xv}\wedge\omega_{vx}\ten\omega_{xy}.
\]
For an edge in a triangle this vector is nonzero: its support is
nonempty and omits \(v=y\) from the return block \(N(x)\).
Flatness therefore implies
\begin{equation}
\label{eq:flat-alpha-product}
 a_{xy}a_{yx}=0.
\end{equation}
For star-preserving connections \(a_{xy}=a_{yx}\), so \(\alpha=0\).

For distinct \(y,c\in N(x)\), set
\[ T_{vyc}=|N(v)\cap N(y)\cap N(c)|.\]
In the return block \(b=x\), a representative coefficient vector
from \eqref{eq:flat-curvature-coeff} is
\begin{align*}
 f_{y,c}(v)={}&\delta_{vy}+(q+\delta_{vc})(qT_{vyc}+1)
    +\delta_{vy}a_{vc}A_{vc}
    +qA_{yc}a_{vc}A_{vc}+a_{xy}(q+\delta_{vc})A_{yc}A_{vy}.
\end{align*}
If \(q=0\), this is supported on \(\{y,c\}\), with values
\[
 f_{y,c}(y)=1+a_{yc}A_{yc},\qquad
 f_{y,c}(c)=1+a_{xy}A_{yc}.
\]
If \(\deg(x)>2\), constancy forces \(y\sim c\) and \(a_{xy}=1\).
Interchanging \(y,c\) also gives \(a_{xc}=1\), whereas the triangle
equations give \(a_{cx}=a_{xy}=1\), contradicting
\eqref{eq:flat-alpha-product}. Every degree is positive and even,
so \(G\) is a cycle. The cycle formulas below and
Remark~\ref{rem:K3twoparam} give the stated classification for
\(\sigma_0\).

If \(q=1\), even common-neighbour counts give
\(f_{y,c}(y)=f_{y,c}(c)=0\). For \(v\notin\{y,c\}\),
\[
 f_{y,c}(v)=1+T_{vyc}
       +A_{yc}(a_{vc}A_{vc}+a_{xy}A_{vy}).
\]
Suppose \(a_{xy}=1\). Then \(a_{yx}=0\), and there are distinct
\(c,d\in N(x)\cap N(y)\), since this set is nonempty and even.
The triangle equations give \(a_{xc}=a_{xd}=a_{dy}=0\), and also
\(a_{dc}A_{dc}=0\). Flatness would require
\(f_{y,c}(d)=f_{c,y}(d)=0\), but their sum is
\[
 (a_{dc}+a_{xc})A_{dc}
       +(a_{dy}+a_{xy})A_{dy}=1.
\]
This contradiction gives \(\alpha=0\) for every flat co-Eulerian
geometry. Together with the Eulerian case, this proves that the only
non-inner exceptions occur on \(K_3\). Inner QLCs are star-preserving
by the general classification.
\end{proof}

\begin{lemma}[Graph criterion for flatness of co-Eulerian geometries]
\label{lem:flat-sigma1}
Let \(G\) be co-Eulerian. Its inner QRG is flat if
and only if, for every edge \(xy\), vertex \(b\) and \(c\in N(b)\),
\begin{equation}
\label{eq:flat-triple-test}
 v\longmapsto\delta_{vy}+(1+\delta_{vc})
       \bigl(|N(v)\cap N(y)\cap N(c)|+A_{xc}\bigr)
\end{equation}
is constant on \(N(x)\cap N(b)\), with values in \(\F_2\).
\end{lemma}

\begin{proof}
Set \(q=1\) and \(a=0\) in \eqref{eq:flat-curvature-coeff}.
The coefficient at \((x,v,b,c)\) is
\[
 \delta_{vy}+(1+\delta_{xb})A_{yb}
 +(1+\delta_{vc})(T_{vyc}+A_{xc}).
\]
The middle term is independent of \(v\) and vanishes in the wedge
quotient. The remaining coefficient is \eqref{eq:flat-triple-test}.
\end{proof}

\begin{proposition}[The braid relation]
\label{prop:braid}
Let \(\sigma_q\) be a QLC braiding on a connected graph over \(\F_2\).
\begin{enumerate}[label=\textup{(\roman*)}]
\item The Eulerian braiding \(\sigma_0\) always satisfies
\[
 \sigma_{12}\sigma_{23}\sigma_{12}
 =\sigma_{23}\sigma_{12}\sigma_{23}.
\]
\item The co-Eulerian braiding \(\sigma_1\) satisfies this relation
if and only if its inner QLC \((\sigma_1,0)\) is flat.
\end{enumerate}
In these cases the adjacent braidings on
\((\Omega^1)^{\otimes_A r}\) define an action of the symmetric group
\(S_r\).
\end{proposition}

\begin{proof}
For \(q=0\), put \(U=\sigma_{12}+\id\) and
\(V=\sigma_{23}+\id\). These maps replace a return two-path by
the sum of all return paths at its base, and vanish on non-return
paths. Even vertex degrees give \(U^2=V^2=0\).
In \(UVU\), the middle operation selects the intermediate vertex
equal to the final vertex of the three-path, and the last operation
selects the original base. Hence \(UVU=U\), and similarly
\(VUV=V\). Expanding the braid relation leaves
\(U+V+UVU+VUV=0\).

For \(q=1\), write
\(T_{abc}=|N(a)\cap N(b)\cap N(c)|\) modulo two.
The coefficient of \((x,u,v,w)\) in the braid defect applied to
\((x,y,z,w)\) is
\begin{align*}
 &\delta_{vz}(1+\delta_{xz})+\delta_{uy}(1+\delta_{yw})\\
 &\quad +(1+\delta_{xz})(1+\delta_{xv})(T_{xzv}+A_{xw})
 +(1+\delta_{yw})(1+\delta_{uw})(T_{yuw}+A_{xw}).
\end{align*}
This follows by expanding with \(\sigma=\id+(\sigma+\id)\);
the squared terms vanish by involutivity.
For adjacent endpoints \(x,w\), these coefficients vanish precisely
when every three distinct vertices with a common neighbour have an
odd number of common neighbours. Necessity follows by comparing
\((x,y,x,w)\) with \((x,u,x,w)\), for distinct
\(u,y,w\in N(x)\); its coefficient is \(T_{uyw}+1\).
Conversely this condition and even distinct-vertex codegrees make
each of the two groups of terms above vanish.
The same condition is equivalent to the adjacent-endpoint part of
Lemma~\ref{lem:flat-sigma1}: its coefficient is constantly
zero if \(y\ne w\), and constantly one if \(y=w\).

For nonadjacent endpoints \(x,w\), including \(x=w\), let \(C\) be the binary
matrix with rows indexed by \(N(x)\), columns by \(N(w)\), and
\(C_{yz}=A_{yz}\). Treat the row and column sets as disjoint
copies when they overlap. Put
\[
 M=I+CC^T,\qquad K=I+C^TC.
\]
Both have diagonal entries one, by even codegrees. The braid-defect
coefficient reduces to
\[
 M_{uy}+K_{vz},\qquad C_{yz}=C_{uv}=1.
\]
The flatness criterion says that, for each column \(v\) and row
\(y\), the entry \(M_{uy}\) is independent of \(u\) among rows
with \(C_{uv}=1\). The braid relation implies this by choosing any
\(z\) with \(C_{yz}=1\). If row \(y\) is isolated, the required
entries are all zero.

Conversely, flatness makes the rows of \(M\) equal whenever their
vertices belong to the same nontrivial connected component of the
bipartite graph of \(C\). Its diagonal entries are one, so
\(M_{uy}\) is one within such a component and zero between components.
Applying flatness with the endpoints reversed gives the same statement
for \(K\). Since \(yz\) and \(uv\) are edges of this bipartite
graph, \(M_{uy}=K_{vz}\). This proves equivalence for all remaining endpoints.

Finally, involutivity and commutation on disjoint tensor legs give
the remaining defining relations of \(S_r\).
\end{proof}

\begin{proposition}[Restrictions on flat co-Eulerian geometries]
\label{prop:flat-regular}
A connected graph of order \(n\geq2\) admitting a flat co-Eulerian geometry
is regular, of degree \(r\), and any three distinct vertices with a
common neighbour have an odd number of common neighbours. If
\(r\geq3\), no two distinct vertices have the same neighbourhood.
If the graph is also Eulerian, its adjacency matrix satisfies
\[
 r\leq\rank_{\F_2}A\leq n/2,\qquad r\leq2\lfloor n/4\rfloor.
\]
If \(n\) is odd, the graph is necessarily Eulerian and \(4\mid r\).
\end{proposition}

\begin{proof}
By Theorem~\ref{thm:flat-alpha} and its proof the connection is inner.
For an edge \(xv\), let
\(P_{xv}:\F_2^{N(x)}\to\F_2^{N(v)}\) have rows indexed by
\(z\in N(v)\), columns indexed by \(y\in N(x)\), and entries
\[
 (P_{xv})_{zy}=B_{xy}(v,z)
 =\delta_{vy}+(1+\delta_{xz})A_{yz}.
\]
The return-block curvature formula
\eqref{eq:flat-curvature-coeff} says that \(P_{vx}P_{xv}\) is
independent of \(v\in N(x)\). Write this common matrix as \(Q_x\).
If \(e_v\) is the coordinate vector indexed by \(v\), then
\(P_{xv}e_v=e_x\) and \(P_{vx}e_x=e_v\). Thus \(Q_x\) fixes
every coordinate vector, so \(Q_x=I\). Interchanging \(x,v\)
also gives \(P_{xv}P_{vx}=I\). The two neighbour spaces have equal
dimension, and connectedness gives regularity.

For distinct \(v,y,c\in N(x)\), take \(b=x\) in
\eqref{eq:flat-triple-test}. Its value at \(y\) is zero, by the
even common-neighbour condition, while its value at \(v\) is
\(|N(v)\cap N(y)\cap N(c)|+1\). Constancy gives the required
oddness.

If \(N(u)=N(v)\), choose a common neighbour \(x\) and, when
\(r\geq3\), a third vertex \(w\in N(x)\setminus\{u,v\}\).
The common-neighbour count of \(u,v,w\) equals
\(\operatorname{codeg}(u,w)\), which is even, a contradiction.

Now suppose \(r\) is even. Write \(a_v=Ae_v\) for the adjacency
columns and fix \(x\). In a relation
\(\sum_{v\in N(x)}c_va_v=0\), the \(x\)-coordinate gives
\(\sum_vc_v=0\). For distinct \(y,z\in N(x)\), take its scalar
product with the coordinatewise product \(a_y\odot a_z\), where \((a_y\odot a_z)_u=A_{uy}A_{uz}\).
The terms \(v=y,z\) vanish by even codegree, and every remaining
term is \(c_v\) by the triple-neighbour condition. Hence
\(0=\sum_{v\ne y,z}c_v=c_y+c_z\), so all coefficients agree.
Conversely, \(\sum_{v\in N(x)}a_v=A^2e_x=0\). These columns
therefore have rank \(r-1\). Since a symmetric zero-diagonal matrix
over \(\F_2\) has even rank, \(\rank A\geq r\). The relation
\(A^2=0\) gives \(\rank A\leq n/2\), and \(r\) is even.

More generally, every Eulerian co-Eulerian graph has an even number
of edges. Indeed, \(A^2=0\) gives \(\det(tI-A)=t^n\), whose
coefficient of \(t^{n-2}\) is \(|E|\) modulo two. For a regular
graph of odd order, \(2|E|=nr\) therefore implies \(4\mid r\).
\end{proof}

\begin{proposition}[A cubical covering]
\label{prop:flat-cube}
Let \(G\) be connected with a flat co-Eulerian geometry and degree
\(r\geq1\). If every distinct-vertex codegree is at most two,
there is a graph covering from the \(r\)-dimensional cube \(Q_r\)
onto \(G\). In particular, the order of \(G\) divides \(2^r\).
\end{proposition}

\begin{proof}
Every nonzero codegree is two. The transport matrix \(P_{xv}\) in
Proposition~\ref{prop:flat-regular} is therefore a permutation: it
sends \(v\in N(x)\) to \(x\in N(v)\), and sends
\(y\in N(x)\setminus\{v\}\) to the unique common neighbour of
\(v,y\) other than \(x\).
Consider states \((x,L)\), where \(L:\{1,\ldots,r\}\to N(x)\)
is a bijection. Define \(T_i\) by moving to \(v=L(i)\) and
transporting the numbering by \(P_{xv}\). The inverse transport
identity gives \(T_i^2=\id\).

For \(i\ne j\), put \(v=L(i)\), \(w=L(j)\), and let \(b\)
be their common neighbour other than \(x\). Both \(T_jT_i\) and
\(T_iT_j\) end at \(b\). The endpoints \(x,b\) have precisely
the two intermediates \(v,w\), and flatness in
\eqref{eq:flat-curvature-coeff} gives
\[
 P_{vb}P_{xv}=P_{wb}P_{xw}.
\]
The resulting numberings agree, so the \(T_i\) commute.
Starting from one state, applying
\(T_1^{a_1}\cdots T_r^{a_r}\) and forgetting the numbering defines
a graph map \(Q_r\to G\), for \((a_1,\ldots,a_r)\in\F_2^r\).
At each cube vertex its \(r\) neighbours map bijectively to the
neighbours of its image. Connectedness makes this a covering onto
\(G\). A finite covering of a connected graph has constant fibre
size, so \(|G|\) divides \(2^r\).
\end{proof}

\begin{corollary}
\label{cor:flat-four}
The connected four-regular graphs admitting a flat co-Eulerian geometry
are precisely \(Q_4\) and the Cartesian product \(K_4\square K_2\).
\end{corollary}

\begin{proof}
By Proposition~\ref{prop:flat-regular}, distinct vertices cannot have
equal neighbourhoods, so every nonzero codegree is two. The induced
graph on any four-element neighbourhood therefore has all degrees
zero or two: it is \(4K_1\), \(K_3\sqcup K_1\), or \(C_4\).
In the last case, opposite vertices have at least three common
neighbours in the original graph, contradicting the codegree bound.
If \(G\) contains a triangle, a vertex in it consequently belongs
to a \(K_4\), say \(B\).
Each \(u\in B\) has exactly one neighbour \(f(u)\) outside \(B\).
These neighbours are distinct, since \(f(u)=f(v)\) would give
\(\operatorname{codeg}(u,v)=3\).
For \(u\ne v\) in \(B\), the vertices \(u,f(v)\) share \(v\),
and their only possible additional common neighbour is \(f(u)\).
Even codegree forces \(f(u)\sim f(v)\). Thus \(f(B)\) is another
\(K_4\), and the matching exhausts every degree. Connectedness gives
\(G=K_4\square K_2\).

If \(G\) is triangle-free, Proposition~\ref{prop:flat-cube} gives
a covering \(Q_4\to G\). Four-regularity requires at least five
vertices, so the order is eight or sixteen. In the eight-vertex case,
the neighbourhoods of the endpoints of any edge are disjoint
independent sets of size four. They partition the vertex set, and
four-regularity forces \(G=K_{4,4}\), contradicting the absence of
equal neighbourhoods. Thus \(G\) has sixteen vertices and the
covering is an isomorphism. Both graphs have flat connections by
Proposition~\ref{prop:cayley-involutions} below.
\end{proof}

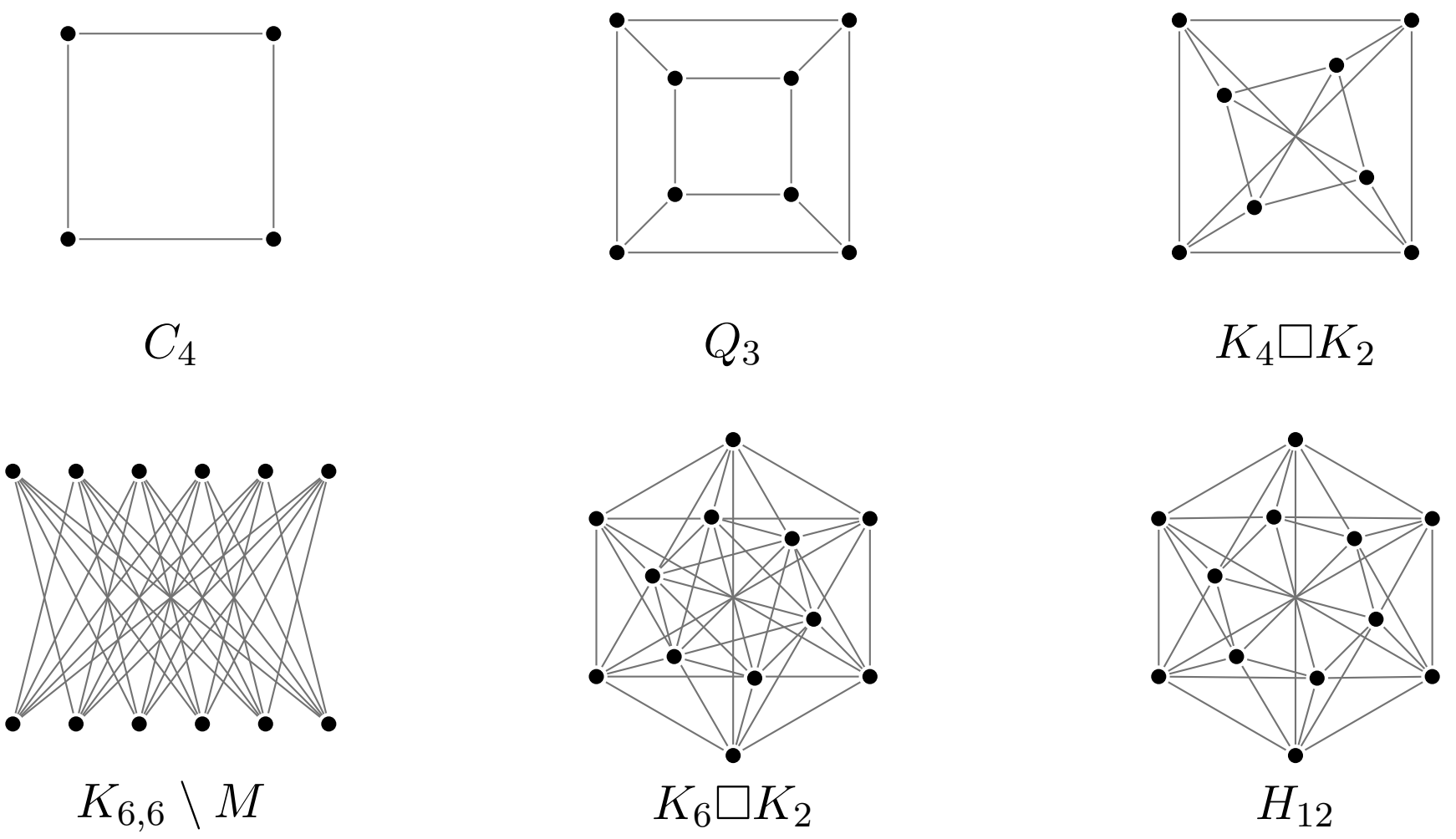
\begin{figure}[tbp]
\centering
\begin{tikzpicture}[x=1cm,y=1cm,line width=.4pt]
\begin{scope}[shift={(0.000,-0.000)}]
\draw[black!55] (0.8132,0.8132) -- (-0.8132,0.8132);
\draw[black!55] (0.8132,0.8132) -- (0.8132,-0.8132);
\draw[black!55] (-0.8132,0.8132) -- (-0.8132,-0.8132);
\draw[black!55] (-0.8132,-0.8132) -- (0.8132,-0.8132);
\fill[white] (0.8132,0.8132) circle[radius=2.5pt];
\fill (0.8132,0.8132) circle[radius=1.7pt];
\fill[white] (-0.8132,0.8132) circle[radius=2.5pt];
\fill (-0.8132,0.8132) circle[radius=1.7pt];
\fill[white] (-0.8132,-0.8132) circle[radius=2.5pt];
\fill (-0.8132,-0.8132) circle[radius=1.7pt];
\fill[white] (0.8132,-0.8132) circle[radius=2.5pt];
\fill (0.8132,-0.8132) circle[radius=1.7pt];
\node at (0,-1.65) {$C_4$};
\end{scope}
\begin{scope}[shift={(4.450,-0.000)}]
\draw[black!55] (0.9192,0.9192) -- (-0.9192,0.9192);
\draw[black!55] (0.9192,0.9192) -- (0.9192,-0.9192);
\draw[black!55] (0.9192,0.9192) -- (0.4596,0.4596);
\draw[black!55] (-0.9192,0.9192) -- (-0.9192,-0.9192);
\draw[black!55] (-0.9192,0.9192) -- (-0.4596,0.4596);
\draw[black!55] (-0.9192,-0.9192) -- (0.9192,-0.9192);
\draw[black!55] (-0.9192,-0.9192) -- (-0.4596,-0.4596);
\draw[black!55] (0.9192,-0.9192) -- (0.4596,-0.4596);
\draw[black!55] (0.4596,0.4596) -- (-0.4596,0.4596);
\draw[black!55] (0.4596,0.4596) -- (0.4596,-0.4596);
\draw[black!55] (-0.4596,0.4596) -- (-0.4596,-0.4596);
\draw[black!55] (-0.4596,-0.4596) -- (0.4596,-0.4596);
\fill[white] (0.9192,0.9192) circle[radius=2.5pt];
\fill (0.9192,0.9192) circle[radius=1.7pt];
\fill[white] (-0.9192,0.9192) circle[radius=2.5pt];
\fill (-0.9192,0.9192) circle[radius=1.7pt];
\fill[white] (-0.9192,-0.9192) circle[radius=2.5pt];
\fill (-0.9192,-0.9192) circle[radius=1.7pt];
\fill[white] (0.9192,-0.9192) circle[radius=2.5pt];
\fill (0.9192,-0.9192) circle[radius=1.7pt];
\fill[white] (0.4596,0.4596) circle[radius=2.5pt];
\fill (0.4596,0.4596) circle[radius=1.7pt];
\fill[white] (-0.4596,0.4596) circle[radius=2.5pt];
\fill (-0.4596,0.4596) circle[radius=1.7pt];
\fill[white] (-0.4596,-0.4596) circle[radius=2.5pt];
\fill (-0.4596,-0.4596) circle[radius=1.7pt];
\fill[white] (0.4596,-0.4596) circle[radius=2.5pt];
\fill (0.4596,-0.4596) circle[radius=1.7pt];
\node at (0,-1.65) {$Q_3$};
\end{scope}
\begin{scope}[shift={(8.900,-0.000)}]
\draw[black!55] (0.9192,0.9192) -- (-0.9192,0.9192);
\draw[black!55] (0.9192,0.9192) -- (-0.9192,-0.9192);
\draw[black!55] (0.9192,0.9192) -- (0.9192,-0.9192);
\draw[black!55] (0.9192,0.9192) -- (0.3250,0.5629);
\draw[black!55] (-0.9192,0.9192) -- (-0.9192,-0.9192);
\draw[black!55] (-0.9192,0.9192) -- (0.9192,-0.9192);
\draw[black!55] (-0.9192,0.9192) -- (-0.5629,0.3250);
\draw[black!55] (-0.9192,-0.9192) -- (0.9192,-0.9192);
\draw[black!55] (-0.9192,-0.9192) -- (-0.3250,-0.5629);
\draw[black!55] (0.9192,-0.9192) -- (0.5629,-0.3250);
\draw[black!55] (0.3250,0.5629) -- (-0.5629,0.3250);
\draw[black!55] (0.3250,0.5629) -- (-0.3250,-0.5629);
\draw[black!55] (0.3250,0.5629) -- (0.5629,-0.3250);
\draw[black!55] (-0.5629,0.3250) -- (-0.3250,-0.5629);
\draw[black!55] (-0.5629,0.3250) -- (0.5629,-0.3250);
\draw[black!55] (-0.3250,-0.5629) -- (0.5629,-0.3250);
\fill[white] (0.9192,0.9192) circle[radius=2.5pt];
\fill (0.9192,0.9192) circle[radius=1.7pt];
\fill[white] (-0.9192,0.9192) circle[radius=2.5pt];
\fill (-0.9192,0.9192) circle[radius=1.7pt];
\fill[white] (-0.9192,-0.9192) circle[radius=2.5pt];
\fill (-0.9192,-0.9192) circle[radius=1.7pt];
\fill[white] (0.9192,-0.9192) circle[radius=2.5pt];
\fill (0.9192,-0.9192) circle[radius=1.7pt];
\fill[white] (0.3250,0.5629) circle[radius=2.5pt];
\fill (0.3250,0.5629) circle[radius=1.7pt];
\fill[white] (-0.5629,0.3250) circle[radius=2.5pt];
\fill (-0.5629,0.3250) circle[radius=1.7pt];
\fill[white] (-0.3250,-0.5629) circle[radius=2.5pt];
\fill (-0.3250,-0.5629) circle[radius=1.7pt];
\fill[white] (0.5629,-0.3250) circle[radius=2.5pt];
\fill (0.5629,-0.3250) circle[radius=1.7pt];
\node at (0,-1.65) {$K_4\square K_2$};
\end{scope}
\begin{scope}[shift={(0.000,-3.650)}]
\draw[black!55] (-1.2500,1.0000) -- (-0.7500,-1.0000);
\draw[black!55] (-1.2500,1.0000) -- (-0.2500,-1.0000);
\draw[black!55] (-1.2500,1.0000) -- (0.2500,-1.0000);
\draw[black!55] (-1.2500,1.0000) -- (0.7500,-1.0000);
\draw[black!55] (-1.2500,1.0000) -- (1.2500,-1.0000);
\draw[black!55] (-0.7500,1.0000) -- (-1.2500,-1.0000);
\draw[black!55] (-0.7500,1.0000) -- (-0.2500,-1.0000);
\draw[black!55] (-0.7500,1.0000) -- (0.2500,-1.0000);
\draw[black!55] (-0.7500,1.0000) -- (0.7500,-1.0000);
\draw[black!55] (-0.7500,1.0000) -- (1.2500,-1.0000);
\draw[black!55] (-0.2500,1.0000) -- (-1.2500,-1.0000);
\draw[black!55] (-0.2500,1.0000) -- (-0.7500,-1.0000);
\draw[black!55] (-0.2500,1.0000) -- (0.2500,-1.0000);
\draw[black!55] (-0.2500,1.0000) -- (0.7500,-1.0000);
\draw[black!55] (-0.2500,1.0000) -- (1.2500,-1.0000);
\draw[black!55] (0.2500,1.0000) -- (-1.2500,-1.0000);
\draw[black!55] (0.2500,1.0000) -- (-0.7500,-1.0000);
\draw[black!55] (0.2500,1.0000) -- (-0.2500,-1.0000);
\draw[black!55] (0.2500,1.0000) -- (0.7500,-1.0000);
\draw[black!55] (0.2500,1.0000) -- (1.2500,-1.0000);
\draw[black!55] (0.7500,1.0000) -- (-1.2500,-1.0000);
\draw[black!55] (0.7500,1.0000) -- (-0.7500,-1.0000);
\draw[black!55] (0.7500,1.0000) -- (-0.2500,-1.0000);
\draw[black!55] (0.7500,1.0000) -- (0.2500,-1.0000);
\draw[black!55] (0.7500,1.0000) -- (1.2500,-1.0000);
\draw[black!55] (1.2500,1.0000) -- (-1.2500,-1.0000);
\draw[black!55] (1.2500,1.0000) -- (-0.7500,-1.0000);
\draw[black!55] (1.2500,1.0000) -- (-0.2500,-1.0000);
\draw[black!55] (1.2500,1.0000) -- (0.2500,-1.0000);
\draw[black!55] (1.2500,1.0000) -- (0.7500,-1.0000);
\fill[white] (-1.2500,1.0000) circle[radius=2.5pt];
\fill (-1.2500,1.0000) circle[radius=1.7pt];
\fill[white] (-0.7500,1.0000) circle[radius=2.5pt];
\fill (-0.7500,1.0000) circle[radius=1.7pt];
\fill[white] (-0.2500,1.0000) circle[radius=2.5pt];
\fill (-0.2500,1.0000) circle[radius=1.7pt];
\fill[white] (0.2500,1.0000) circle[radius=2.5pt];
\fill (0.2500,1.0000) circle[radius=1.7pt];
\fill[white] (0.7500,1.0000) circle[radius=2.5pt];
\fill (0.7500,1.0000) circle[radius=1.7pt];
\fill[white] (1.2500,1.0000) circle[radius=2.5pt];
\fill (1.2500,1.0000) circle[radius=1.7pt];
\fill[white] (-1.2500,-1.0000) circle[radius=2.5pt];
\fill (-1.2500,-1.0000) circle[radius=1.7pt];
\fill[white] (-0.7500,-1.0000) circle[radius=2.5pt];
\fill (-0.7500,-1.0000) circle[radius=1.7pt];
\fill[white] (-0.2500,-1.0000) circle[radius=2.5pt];
\fill (-0.2500,-1.0000) circle[radius=1.7pt];
\fill[white] (0.2500,-1.0000) circle[radius=2.5pt];
\fill (0.2500,-1.0000) circle[radius=1.7pt];
\fill[white] (0.7500,-1.0000) circle[radius=2.5pt];
\fill (0.7500,-1.0000) circle[radius=1.7pt];
\fill[white] (1.2500,-1.0000) circle[radius=2.5pt];
\fill (1.2500,-1.0000) circle[radius=1.7pt];
\node at (0,-1.65) {$K_{6,6}\setminus M$};
\end{scope}
\begin{scope}[shift={(4.450,-3.650)}]
\draw[black!55] (0.0000,1.2500) -- (-1.0825,0.6250);
\draw[black!55] (0.0000,1.2500) -- (-1.0825,-0.6250);
\draw[black!55] (0.0000,1.2500) -- (-0.0000,-1.2500);
\draw[black!55] (0.0000,1.2500) -- (1.0825,-0.6250);
\draw[black!55] (0.0000,1.2500) -- (1.0825,0.6250);
\draw[black!55] (0.0000,1.2500) -- (-0.1708,0.6375);
\draw[black!55] (-1.0825,0.6250) -- (-1.0825,-0.6250);
\draw[black!55] (-1.0825,0.6250) -- (-0.0000,-1.2500);
\draw[black!55] (-1.0825,0.6250) -- (1.0825,-0.6250);
\draw[black!55] (-1.0825,0.6250) -- (1.0825,0.6250);
\draw[black!55] (-1.0825,0.6250) -- (-0.6375,0.1708);
\draw[black!55] (-1.0825,-0.6250) -- (-0.0000,-1.2500);
\draw[black!55] (-1.0825,-0.6250) -- (1.0825,-0.6250);
\draw[black!55] (-1.0825,-0.6250) -- (1.0825,0.6250);
\draw[black!55] (-1.0825,-0.6250) -- (-0.4667,-0.4667);
\draw[black!55] (-0.0000,-1.2500) -- (1.0825,-0.6250);
\draw[black!55] (-0.0000,-1.2500) -- (1.0825,0.6250);
\draw[black!55] (-0.0000,-1.2500) -- (0.1708,-0.6375);
\draw[black!55] (1.0825,-0.6250) -- (1.0825,0.6250);
\draw[black!55] (1.0825,-0.6250) -- (0.6375,-0.1708);
\draw[black!55] (1.0825,0.6250) -- (0.4667,0.4667);
\draw[black!55] (-0.1708,0.6375) -- (-0.6375,0.1708);
\draw[black!55] (-0.1708,0.6375) -- (-0.4667,-0.4667);
\draw[black!55] (-0.1708,0.6375) -- (0.1708,-0.6375);
\draw[black!55] (-0.1708,0.6375) -- (0.6375,-0.1708);
\draw[black!55] (-0.1708,0.6375) -- (0.4667,0.4667);
\draw[black!55] (-0.6375,0.1708) -- (-0.4667,-0.4667);
\draw[black!55] (-0.6375,0.1708) -- (0.1708,-0.6375);
\draw[black!55] (-0.6375,0.1708) -- (0.6375,-0.1708);
\draw[black!55] (-0.6375,0.1708) -- (0.4667,0.4667);
\draw[black!55] (-0.4667,-0.4667) -- (0.1708,-0.6375);
\draw[black!55] (-0.4667,-0.4667) -- (0.6375,-0.1708);
\draw[black!55] (-0.4667,-0.4667) -- (0.4667,0.4667);
\draw[black!55] (0.1708,-0.6375) -- (0.6375,-0.1708);
\draw[black!55] (0.1708,-0.6375) -- (0.4667,0.4667);
\draw[black!55] (0.6375,-0.1708) -- (0.4667,0.4667);
\fill[white] (0.0000,1.2500) circle[radius=2.5pt];
\fill (0.0000,1.2500) circle[radius=1.7pt];
\fill[white] (-1.0825,0.6250) circle[radius=2.5pt];
\fill (-1.0825,0.6250) circle[radius=1.7pt];
\fill[white] (-1.0825,-0.6250) circle[radius=2.5pt];
\fill (-1.0825,-0.6250) circle[radius=1.7pt];
\fill[white] (-0.0000,-1.2500) circle[radius=2.5pt];
\fill (-0.0000,-1.2500) circle[radius=1.7pt];
\fill[white] (1.0825,-0.6250) circle[radius=2.5pt];
\fill (1.0825,-0.6250) circle[radius=1.7pt];
\fill[white] (1.0825,0.6250) circle[radius=2.5pt];
\fill (1.0825,0.6250) circle[radius=1.7pt];
\fill[white] (-0.1708,0.6375) circle[radius=2.5pt];
\fill (-0.1708,0.6375) circle[radius=1.7pt];
\fill[white] (-0.6375,0.1708) circle[radius=2.5pt];
\fill (-0.6375,0.1708) circle[radius=1.7pt];
\fill[white] (-0.4667,-0.4667) circle[radius=2.5pt];
\fill (-0.4667,-0.4667) circle[radius=1.7pt];
\fill[white] (0.1708,-0.6375) circle[radius=2.5pt];
\fill (0.1708,-0.6375) circle[radius=1.7pt];
\fill[white] (0.6375,-0.1708) circle[radius=2.5pt];
\fill (0.6375,-0.1708) circle[radius=1.7pt];
\fill[white] (0.4667,0.4667) circle[radius=2.5pt];
\fill (0.4667,0.4667) circle[radius=1.7pt];
\node at (0,-1.65) {$K_6\square K_2$};
\end{scope}
\begin{scope}[shift={(8.900,-3.650)}]
\draw[black!55] (0.0000,1.2500) -- (-0.1708,0.6375);
\draw[black!55] (0.0000,1.2500) -- (-1.0825,0.6250);
\draw[black!55] (0.0000,1.2500) -- (-0.6375,0.1708);
\draw[black!55] (0.0000,1.2500) -- (-0.0000,-1.2500);
\draw[black!55] (0.0000,1.2500) -- (1.0825,0.6250);
\draw[black!55] (0.0000,1.2500) -- (0.4667,0.4667);
\draw[black!55] (-0.1708,0.6375) -- (-1.0825,0.6250);
\draw[black!55] (-0.1708,0.6375) -- (-0.6375,0.1708);
\draw[black!55] (-0.1708,0.6375) -- (0.1708,-0.6375);
\draw[black!55] (-0.1708,0.6375) -- (1.0825,0.6250);
\draw[black!55] (-0.1708,0.6375) -- (0.4667,0.4667);
\draw[black!55] (-1.0825,0.6250) -- (-0.6375,0.1708);
\draw[black!55] (-1.0825,0.6250) -- (-1.0825,-0.6250);
\draw[black!55] (-1.0825,0.6250) -- (-0.4667,-0.4667);
\draw[black!55] (-1.0825,0.6250) -- (1.0825,-0.6250);
\draw[black!55] (-0.6375,0.1708) -- (-1.0825,-0.6250);
\draw[black!55] (-0.6375,0.1708) -- (-0.4667,-0.4667);
\draw[black!55] (-0.6375,0.1708) -- (0.6375,-0.1708);
\draw[black!55] (-1.0825,-0.6250) -- (-0.4667,-0.4667);
\draw[black!55] (-1.0825,-0.6250) -- (-0.0000,-1.2500);
\draw[black!55] (-1.0825,-0.6250) -- (0.1708,-0.6375);
\draw[black!55] (-1.0825,-0.6250) -- (1.0825,0.6250);
\draw[black!55] (-0.4667,-0.4667) -- (-0.0000,-1.2500);
\draw[black!55] (-0.4667,-0.4667) -- (0.1708,-0.6375);
\draw[black!55] (-0.4667,-0.4667) -- (0.4667,0.4667);
\draw[black!55] (-0.0000,-1.2500) -- (0.1708,-0.6375);
\draw[black!55] (-0.0000,-1.2500) -- (1.0825,-0.6250);
\draw[black!55] (-0.0000,-1.2500) -- (0.6375,-0.1708);
\draw[black!55] (0.1708,-0.6375) -- (1.0825,-0.6250);
\draw[black!55] (0.1708,-0.6375) -- (0.6375,-0.1708);
\draw[black!55] (1.0825,-0.6250) -- (0.6375,-0.1708);
\draw[black!55] (1.0825,-0.6250) -- (1.0825,0.6250);
\draw[black!55] (1.0825,-0.6250) -- (0.4667,0.4667);
\draw[black!55] (0.6375,-0.1708) -- (1.0825,0.6250);
\draw[black!55] (0.6375,-0.1708) -- (0.4667,0.4667);
\draw[black!55] (1.0825,0.6250) -- (0.4667,0.4667);
\fill[white] (0.0000,1.2500) circle[radius=2.5pt];
\fill (0.0000,1.2500) circle[radius=1.7pt];
\fill[white] (-0.1708,0.6375) circle[radius=2.5pt];
\fill (-0.1708,0.6375) circle[radius=1.7pt];
\fill[white] (-1.0825,0.6250) circle[radius=2.5pt];
\fill (-1.0825,0.6250) circle[radius=1.7pt];
\fill[white] (-0.6375,0.1708) circle[radius=2.5pt];
\fill (-0.6375,0.1708) circle[radius=1.7pt];
\fill[white] (-1.0825,-0.6250) circle[radius=2.5pt];
\fill (-1.0825,-0.6250) circle[radius=1.7pt];
\fill[white] (-0.4667,-0.4667) circle[radius=2.5pt];
\fill (-0.4667,-0.4667) circle[radius=1.7pt];
\fill[white] (-0.0000,-1.2500) circle[radius=2.5pt];
\fill (-0.0000,-1.2500) circle[radius=1.7pt];
\fill[white] (0.1708,-0.6375) circle[radius=2.5pt];
\fill (0.1708,-0.6375) circle[radius=1.7pt];
\fill[white] (1.0825,-0.6250) circle[radius=2.5pt];
\fill (1.0825,-0.6250) circle[radius=1.7pt];
\fill[white] (0.6375,-0.1708) circle[radius=2.5pt];
\fill (0.6375,-0.1708) circle[radius=1.7pt];
\fill[white] (1.0825,0.6250) circle[radius=2.5pt];
\fill (1.0825,0.6250) circle[radius=1.7pt];
\fill[white] (0.4667,0.4667) circle[radius=2.5pt];
\fill (0.4667,0.4667) circle[radius=1.7pt];
\node at (0,-1.65) {$H_{12}$};
\end{scope}
\end{tikzpicture}
\caption{Additional flat co-Eulerian QRGs through twelve vertices, all with $\alpha=0$.}
\label{fig:flat-extra}
\end{figure}

For \(n>3\), every flat QLC is inner and star-preserving by
Theorem~\ref{thm:flat-alpha}. Apart from the Eulerian geometries on
\(C_n\) and the co-Eulerian geometries on \(K_n\) for even \(n\),
the enumeration through twelve vertices gives precisely the additional
flat QRGs in Figure~\ref{fig:flat-extra}. Here \(M\) is a perfect
matching and
\[
 H_{12}=\operatorname{Cay}\bigl(\Z_6\times\Z_2,
 \{(0,1),(3,0),(\pm1,0),(\pm1,1)\}\bigr).
\]
All six geometries are co-Eulerian, as required by
Theorem~\ref{thm:flat-alpha}(ii) for non-cycle graphs.

\begin{conjecture}
\label{conj:odd-flat}
A connected graph of odd order \(n>3\) admits a flat QLC if and
only if it is \(C_n\), whose unique flat connection is the inner 
\(*\)-preserving one.
\end{conjecture}

By Theorem~\ref{thm:flat-alpha}, it remains to exclude flat
co-Eulerian geometries at odd order. Proposition~\ref{prop:flat-regular}
forces a possible counterexample to be regular of degree divisible
by four, with \(r\leq2\lfloor n/4\rfloor\).
Corollary~\ref{cor:flat-four} excludes degree four, so
\(r\geq8\) and \(n\geq17\). Thus our theoretical results imply the conjecture for all odd \(3<n<17\). 
Proposition~\ref{prop:flat-cube} also requires some codegree to be at
least four. At these codegrees the transport matrices need not be
permutations, so the cubical argument does not apply.
Section~\ref{ssec:cayley} proves the conjecture for all odd-order
Cayley graphs with ad-invariant generators.

\section{Cayley graphs}
\label{ssec:cayley}

Let \(H\) be a finite group with inverse-closed generating set
\(\CC\subset H\setminus\{e\}\), and write
\(u_{\CC}=\sum_{s\in \CC}s\in\F_2H\).

\begin{corollary}
\label{cor:cayley-parity}
The graph \(\operatorname{Cay}(H,\CC)\) is Eulerian if and only if \(|\CC|\) is even, and is co-Eulerian if and only if
\[
 u_{\CC}^2=\begin{cases}e &{\rm if}\  |\CC|\ {\rm odd}\\ 0 &{\rm if}\  |\CC|\ {\rm even}.\end{cases}
\]
\end{corollary}

\begin{proof}
The degree is \(|\CC|\). The coefficient of \(h\) in \(u_{\CC}^2\)
counts factorizations \(h=st\) with \(s,t\in \CC\), modulo two.
These correspond to common neighbours \(s\) of \(e\) and \(h\),
since \(s=ht^{-1}\). The coefficient at \(e\) is
\(|\CC|\bmod2\), giving the co-Eulerian criterion.
\end{proof}

If \(\CC\) is conjugation-invariant, \(u_{\CC}\) belongs to the centre
of \(\F_2H\). If \(H\) is abelian, \(u_{\CC}^2=\sum_{s\in \CC}s^2\),
so the second criterion says that every nonidentity element has an
even number of square roots in \(\CC\). 

\begin{proposition}[A family of flat Cayley geometries]
\label{prop:cayley-involutions}
Let \(H=(\Z_2)^k\) and let \(\CC\subset H\setminus\{0\}\) generate
\(H\). If the sums \(s+t\), over unordered pairs of distinct
elements of \(\CC\), are all distinct, then \(\operatorname{Cay}(H,\CC)\) is co-Eulerian and its inner co-Eulerian connection is 
 flat, with
\[
 \sigma_1(e_s\ten e_t)=e_t\ten e_s,\qquad \nabla e_s=0.
\]
\end{proposition}

\begin{proof}
The group-algebra criterion holds since
\(u_{\CC}^2=(|\CC|\bmod2)e\). A non-returning two-path of displacement
\(s+t\) has precisely the two intermediates corresponding to
\((s,t)\) and \((t,s)\); \(\sigma_1=I+J\) interchanges them.
Return paths have \(s=t\) and are fixed. Hence \(\sigma_1\) is the
displayed flip, and
\(\nabla e_s=\theta\ten e_s+\sigma_1(e_s\ten\theta)=0\), proving
flatness.
\end{proof}

Taking \(\CC\) to be a basis gives the hypercube in every dimension.
For independent \(a,b,c\), the set \(\CC=\{a,b,a+b,c\}\) gives
\(K_4\square K_2\), establishing the existence assertion in
Corollary~\ref{cor:flat-four}.

\begin{proposition}[Odd-order groups with ad-invariant generators]
\label{prop:cayley-odd}
Let \(H\) have odd order $n\geq3$ and let \(\CC\) be ad-invariant (so the graph calculus is bicovariant). Then
\(\operatorname{Cay}(H,\CC)\) is Eulerian but not co-Eulerian. 
It admits a flat star-preserving QLC if and only if
\(\operatorname{Cay}(H,\CC)\cong C_n\), with its unique inner flat connection.
\end{proposition}

\begin{proof}
Inversion has no fixed points in \(\CC\), so \(|\CC|\) is even and
\(\sigma_0\) exists. Ad-invariance makes \(u_{\CC}\) central in
\(\F_2H\). By Maschke's theorem this algebra is semisimple;
its centre is a product of fields and has no nonzero nilpotents.
Since \(u_{\CC}\ne0\), the equation \(u_{\CC}^2=0\) required for
\(\sigma_1\) cannot hold.
By Theorem~\ref{thm:flat-alpha}, flatness of an Eulerian geometry forces
the graph to be a cycle, equivalently \(|\CC|=2\). In odd order,
an inverse-closed two-element generating set has the form
\(\{s,s^{-1}\}\), with \(H=\langle s\rangle\). The cycle has
exactly one flat star-preserving QLC, including the case \(|H|=3\).
\end{proof}

\subsubsection{Complete graphs}

Use the cyclic presentation \(K_n=\operatorname{Cay}(\Z_n,I_n)\),
where \(I_n=\Z_n\setminus\{0\}\), and put
\[
 D_s=\sum_{\substack{a,b\in I_n\\a+b=s}}e_a\ten e_b,
 \qquad \bar n=n\bmod2.
\]
Thus  the metric is \(g=D_0\) and \(\theta=\sum_{i\in I_n}e_i\).

\begin{theorem}[QLCs on complete graphs]
\label{thm:Knexistence}
\label{thm:Knunique}
\label{prop:Kn-sigma-unique}
\label{cor:nabla-explicit}
For \(n\geq3\), the graph \(K_n\) has the unique QLC braiding
\[
 \sigma(e_i\ten e_j)
 =e_i\ten e_j+(1+\bar n+\delta_{i+j,0})D_{i+j}.
\]
For \(t\in\F_2\), the maps \(\alpha_t(e_i)=tD_i\) give
star-preserving QLCs
\begin{equation}
\label{eq:nabla-explicit}
 \nabla_t e_i=\theta\ten e_i+e_i\ten\theta+
 \begin{cases}
 \displaystyle\sum_{s\ne0,i}D_s+tD_i,&n\text{ even},\\[3pt]
 D_0+tD_i,&n\text{ odd}.
 \end{cases}
\end{equation}
For \(n\geq4\), these are all the QLCs.
\end{theorem}

\begin{proof}
The vertex degrees are \(n-1\), and distinct vertices have \(n-2\)
common neighbours. Theorem~\ref{thm:general-F2-classification}
therefore gives precisely \(\sigma_{1+\bar n}\).
The coefficients \(a_{xy}=t\) satisfy the triangle and star equations
of Theorem~\ref{thm:general-alpha-kernel}. Substitution into the
inner connection formula gives \eqref{eq:nabla-explicit}.
For \(n\geq4\), the triangle equations \(a_{xw}=a_{ux}\)
force all incoming coefficients at each vertex to agree: given two
neighbours \(u,v\) of \(x\), choose \(w\notin\{x,u,v\}\).
Writing their common value as \(b_x\), the same equations give
\(a_{xw}=b_x=b_w\). Thus all coefficients have a common value \(t\in\F_2\).
\end{proof}

Put
\[
C_i=\sum_{v\ne0}e_i\wedge e_v\ten e_{-v},\qquad
C'_i=e_i\wedge\theta\ten\theta
       +\sum_{v\ne0}\theta\wedge e_v\ten e_{-v}.
\]
\begin{theorem}[Curvature of the complete-graph connections]
\label{thm:Knflatness}
For \(n\geq3\), the connections \eqref{eq:nabla-explicit} have curvature
\[
R_{\nabla_t}(e_i)=tC_i+(1+t)\bar n C'_i.
\]
\begin{enumerate}[label=\textup{(\roman*)}]
\item For even \(n\geq4\), \(\nabla_0\) is flat and \(\nabla_1\) is not.
\item For odd \(n\geq5\), neither connection is flat.
\item For \(n=3\), \(C'_i=0\), so \(\nabla_0\) is flat and
      \(\nabla_1\) is not.
\end{enumerate}
\end{theorem}

\begin{proof}
All displacement indices below are read modulo \(n\), and tensor indices
are nonzero. Put \(h=\bar n\), \(k=1+h\), and \(l=k+t\). The
coefficient of \(e_a\ten e_b\) in \(\nabla_t(e_i)\) is
\[
N_i(a,b)=\delta_{a,i}+\delta_{b,i}+k+\delta_{a+b,0}
                        +l\delta_{a+b,i}.
\]
Since coefficients are constant, expansion of
\((\extd\ten\id+\id\wedge\nabla_t)\nabla_t\) gives the raw
coefficient of \(e_u\wedge e_v\ten e_w\) as
\[
N_i(u,w)+N_i(v,w)+\sum_{b\ne0}N_i(u,b)N_b(v,w).
\]
The row sum is \(\sum_{b\ne0}N_i(u,b)=t(1+\delta_{u,i})\).
Substitution and collection, writing \(r=u+v\), gives
\[
K_i(r,w)+h(1+t)(\delta_{u,i}+\delta_{v+w,0})
              +t\delta_{u,i}\delta_{v+w,0},
\]
where the term independent of the intermediate index is
\[
K_i(r,w)=\delta_{w,i}+\delta_{r,0}+l\delta_{r,i}
          +k+l(\delta_{r+w,0}+\delta_{r+w,i}).
\]
For fixed \(r,w\), its contribution is
\(K_i(r,w)\,\wedge D_r\ten e_w=0\), by the path-sum
relations. The other two terms are precisely
\(h(1+t)C'_i+tC_i\), proving the formula.

The connection \(\nabla_1\) is star-preserving with nonzero
\(\alpha\), so is not flat by Theorem~\ref{thm:flat-alpha}.
For even \(n\), the formula gives \(R_{\nabla_0}=0\).
For odd \(n\geq5\), inspect the return-displacement block \(u+v=0\)
in \(C'_i\), with third factor \(e_w\), \(w\ne i\). Its coefficient vector has ones
at \(u=i\) and \(u=w\), and zeros elsewhere. Since \(n-1\geq4\),
this is nonzero in the quotient.
For \(n=3\), the wedge relations are
\(e_1\wedge e_1=e_2\wedge e_2=0\) and
\(e_1\wedge e_2=e_2\wedge e_1\). Thus
\(e_i\wedge\theta=e_1\wedge e_2=\theta\wedge e_v\), so
\(C'_i=e_1\wedge e_2\ten\theta+e_1\wedge e_2\ten\theta=0\).
\end{proof}

\begin{remark}[The exceptional graph \(K_3\)]
\label{rem:K3twoparam}
On \(K_3\), the two oriented triangle components give the full family
\[
 \nabla_{t_1,t_2}e_1=t_1e_2\ten e_2,\qquad
 \nabla_{t_1,t_2}e_2=t_2e_1\ten e_1,
 \qquad (t_1,t_2)\in\F_2^2.
\]
Its curvature is
\[
 R_{\nabla_{t_1,t_2}}(e_i)=t_1t_2 e_1\wedge e_2\ten e_i\qquad(i=1,2).
\]
It is flat precisely when \(t_1t_2=0\), and star-preserving
precisely when \(t_1=t_2\). The latter is the family
\eqref{eq:nabla-explicit}. On \(K_2\), the unique QLC has
\(\sigma=\id\) and \(\alpha=0\).
\end{remark}

\subsubsection{Cycle graphs}

Use \(C_n=\operatorname{Cay}(\Z_n,\{1,-1\})\) and its basis
\(e_1,e_{-1}\).

\begin{proposition}[QLCs on cycles]
\label{prop:Cn}
\label{conj:Cnunique}
For \(n\geq5\), the unique QLC on \(C_n\) is the Eulerian flat,
star-preserving inner one with
\[
 \sigma(e_a\ten e_b)=e_b\ten e_a\quad(a,b\in\{1,-1\}),
 \qquad \nabla e_{\pm1}=0.
\]
\end{proposition}

\begin{proof}
Every vertex has degree two, whereas vertices at distance two have
one common neighbour. Theorem~\ref{thm:general-F2-classification}
therefore gives only \(\sigma_0\), which is the displayed flip.
There are no triangles, so \(\alpha=0\). The connection formula
gives \(\nabla e_{\pm1}=0\), hence flatness; star preservation
follows from Corollary~\ref{conj:automatic-invertibility}.
\end{proof}

For \(C_4\), both braidings occur and \(\alpha=0\). The braiding
\(\sigma_0\) is again the flip, with \(\nabla e_{\pm1}=0\).
The other connection is
\[
 \sigma_1(e_{\pm1}\ten e_{\pm1})=e_{\mp1}\ten e_{\mp1},
 \qquad
 \sigma_1(e_{\pm1}\ten e_{\mp1})=e_{\pm1}\ten e_{\mp1},
 \qquad \nabla e_{\pm1}=\theta\ten\theta.
\]
Both are star-preserving and flat; for the second connection,
\(\extd\theta=0=\nabla\theta\).

\section{Asymptotic prevalence of Eulerian graphs}
\label{ssec:asymptotic}

For fixed \(n\), choose uniformly from the connected graph
isomorphism classes. The probability of admitting a QLC is
\[
 P_n=\frac{\#\{G:\text{\(G\) admits a QLC}\}}{\#\{G\}}.
\]
This is also the probability of admitting a star-preserving QLC,
since each QLC braiding gives the star-preserving connection
\((\sigma_q,0)\). Let \(g_n\) and \(e_n\) count connected graphs and connected Eulerian
graphs, respectively, up to isomorphism. The latter are precisely the
graphs admitting \(\sigma_0\). Their exact enumeration can be obtained
by Burnside's lemma and the component decomposition
\cite{Rob}. For example, if \(u_n\) counts all Eulerian graphs,
including disconnected ones, then the formal generating function is
\[
 U(x)=1+\sum_{n\geq1}u_nx^n
     =\prod_{m\geq1}(1-x^m)^{-e_m}.
\]
The analogous identity recovers \(g_n\) from the counts of all graphs.

\begin{proposition}[Asymptotic QLC existence]
\label{prop:qlc-asymptotic}
As \(n\to\infty\),
\[
 g_n\sim\frac{2^{\binom n2}}{n!},\qquad
 e_n\sim\frac{2^{\binom{n-1}2}}{n!},\qquad P_n\sim2^{1-n}.
\]
The number of inner QLCs is also asymptotic to \(e_n\).
Writing \(b_n\) for the number of connected co-Eulerian graph
isomorphism classes of order \(n\), we have
\[
 b_n\leq2^{n^2/4+O(n)}.
\]
\end{proposition}

\begin{proof}
There are \(2^{\binom n2}\) labelled graphs and
\(2^{\binom{n-1}2}\) labelled even-degree graphs. For the latter,
choose the edges among \(n-1\) vertices freely; their degree parities
force the edges to the last vertex, whose degree is then even.
Thus the even-degree probability is exactly \(2^{1-n}\) before
imposing connectedness or forgetting labels.

We now pass from labelled graphs to isomorphism classes and then restrict to connected graphs.
For a permutation moving \(k\geq2\) vertices, the deficit in its
number of orbits on unordered pairs is at least
\[
 \Delta_k=\frac{k(2n-k-2)}4.
\]
Indeed, moved--fixed pairs contribute a deficit at least
\(k(n-k)/2\), and moved--moved pairs contribute at least
\(k(k-2)/4\). An invariant graph is specified by one value in \(\F_2\) per pair orbit.
On these values, evenness imposes at least \(n-k-1\) independent
conditions, using edges between fixed vertices; when this bound is
negative it is vacuous. Relative to the identity term, the total
nonidentity contribution to the even-degree Burnside count is at most
\[
 \sum_{k=2}^n n^k2^{-\Delta_k+k}=o(1).
\]
The ordinary graph count has the stronger bound without \(+k\).
Hence forgetting labels divides the leading counts by \(n!\).
For a specified cut of sizes \(k,n-k\), the fraction of labelled
even-degree graphs having no crossing edges is
\(2^{1-k(n-k)}\). Summing over cuts gives \(o(1)\); the ordinary
graph bound is smaller. The same Burnside bounds control the
nonidentity terms for disconnected graphs, proving the stated
connected asymptotics.

It remains to bound co-Eulerian graphs. Such a connected graph
has adjacency matrix \(A\) satisfying \(A^2=\delta I\), where
\(\delta\in\F_2\) is its common degree parity. The symmetric matrix
\(B=A+\delta I\) satisfies \(B^2=0\). Its image \(W\) is totally
isotropic for the standard bilinear form, since
\((Bu)^T(Bv)=u^TB^2v=0\), so \(r=\dim W\leq n/2\).
An ordered basis of \(W\) is a tuple of independent orthogonal
vectors. After \(j\) vectors have been chosen, their orthogonal
complement has \(2^{n-j}\) elements. Dividing this upper bound on
tuples by \(|\mathrm{GL}(r,2)|\geq2^{r^2}/4\) bounds the number
of possible images by
\[
 4\,2^{rn-r(r-1)/2-r^2}.
\]
For a fixed image with basis matrix \(X\), symmetry gives
\(B=XCX^T\) with \(C\) symmetric, allowing at most
\(2^{r(r+1)/2}\) choices. Counting both parities and all ranks,
and bounding unlabelled graphs by labelled ones, gives
\[
 b_n\leq8\sum_{r=0}^{\lfloor n/2\rfloor}2^{r(n-r+1)}
      =2^{n^2/4+O(n)}=o(e_n).
\]
The number of graphs admitting QLCs lies between \(e_n\) and
\(e_n+b_n\), while the number of inner QLCs is exactly
\(e_n+b_n\). Dividing the former by \(g_n\) proves the result.
\end{proof}

In particular, asymptotically almost every connected graph admitting
a QLC is Eulerian: the proportion is \(e_n/(e_n+o(e_n))\to1\).
This is a statement about graph isomorphism classes, rather than
connection-weighted counts.
The straight-line asymptote in Figure~\ref{fig:QLC-frequencies} is
\[
 \log_{10}P_n=-(n-1)\log_{10}2+o(1).
\]
The probability asymptotically halves with each additional vertex.

\begin{proposition}[Co-Eulerian connection counts]
\label{prop:coeu-counts}
For a co-Eulerian graph, \(d^*(G)\leq d(G)\leq |E|/6\).
Consequently the total number of co-Eulerian geometries on connected
graphs of order \(n\) is at most \(2^{n^2/3+O(n)}=o(e_n)\).
\end{proposition}

\begin{proof}
Consider the subgraph formed by the original edges in one component
of \(D_*(G)\). Every such edge belongs to a positive even number of
triangles, all of whose edges remain in the component. This subgraph
has minimum degree at least three and hence at least six edges.
If its derived cover is connected, it contributes one to \(d(G)\).
If the cover splits, orient its edges so that every triangle is cyclic,
using Theorem~\ref{thm:derived-cover}. An edge \(x\to y\) belongs
to at least two triangles, giving at least two incoming edges at
\(x\) and two outgoing edges at \(y\). Every vertex has both
orientations incident to it, and applying this observation to each
shows that both its indegree and outdegree are at least two.
The subgraph therefore has minimum degree at least four. It cannot
have only five vertices, since it would be \(K_5\), whose triangles
cannot all be cyclically oriented. It has at least six vertices and
twelve edges, and contributes two to \(d(G)\).
The components have disjoint edge sets, proving \(6d(G)\leq|E|\).
The double cover also gives \(d^*(G)\leq d(G)\).

Each co-Eulerian graph contributes \(2^{d(G)}\) co-Eulerian geometries.
Proposition~\ref{prop:qlc-asymptotic} bounds
their total by
\[
 b_n\,2^{n(n-1)/12}\leq2^{n^2/3+O(n)}=o(e_n).
\]
\end{proof}

\begin{theorem}[Typical Eulerian moduli]
\label{thm:typical-eulerian}
For a uniformly chosen connected Eulerian graph isomorphism class
of order \(n\), the probability that \(d(G)=d^*(G)=1\) tends to one.
Almost every connected graph admitting QLCs has exactly two QLCs,
both star-preserving.
\end{theorem}

\begin{proof}
First choose uniformly among labelled even-degree graphs on \(n\)
vertices. For a fixed vertex \(x\), its neighbourhood \(S=N(x)\)
is uniform among the even subsets of the other \(n-1\) vertices:
conditional on \(S\), the graph \(G-x\) has prescribed degree
parities \(1_S\), with \(2^{\binom{n-2}2}\) realisations.
If \(S\) omits a vertex \(z\ne x\), all edges of \(G-\{x,z\}\)
can be chosen freely, and the prescribed parities force the edges
to \(z\). Thus, conditional on such an \(S\), the induced graph
on \(S\) is uniform among all graphs on \(|S|\) vertices.

For \(m\geq2\), a uniform labelled graph on \(m\) vertices is
disconnected with probability at most
\[
 \sum_{k=1}^{\lfloor m/2\rfloor}\binom mk2^{-k(m-k)}
 \leq 2m\,2^{-m/2}.
\]
For \(m\geq4\), bound the sum by
\((1+2^{-m/2})^m-1\leq2m2^{-m/2}\); for \(m=2,3\) the
last bound is already at least one. Put
\(\rho=(1+2^{-1/2})/2<1\). Averaging over the even subsets \(S\),
and treating \(S=\varnothing\) and \(S=V\setminus\{x\}\)
separately, bounds the probability of an empty or disconnected
neighbourhood by
\[
 2^{3-n}+4(n-1)\rho^{n-1}.
\]
A union bound shows that all neighbourhoods are nonempty and
connected with probability tending to one.

Connectedness of \(G\) also has probability tending to one, by the
proof of Proposition~\ref{prop:qlc-asymptotic}. A path in \(G[N(x)]\)
links any two edges incident at \(x\) through triangles. Hence, when
\(G\) and all its neighbourhoods are connected, \(D_*(G)\) is
connected. Moreover, the induced graph on any fixed \(n-1\) vertices
is uniform among all graphs. Disjoint sets of four vertices give
\[
 \Pr(G\text{ contains no }K_4)
 \leq(63/64)^{\lfloor(n-1)/4\rfloor}=o(1).
\]
The derived cover consequently does not split, proving
\(d=d^*=1\) with probability tending to one.

These exceptional sets are invariant under relabelling. Their identity
contribution to Burnside's count is \(o(2^{\binom{n-1}2})\), and
their nonidentity contribution is bounded by the same order of error
in Proposition~\ref{prop:qlc-asymptotic}. Thus the assertion holds
for connected graph isomorphism classes as well. Finally, only
\(o(e_n)\) classes are co-Eulerian, so almost every graph admitting
QLCs has only \(\sigma_0\), with its two star-preserving
\(\alpha\)-values.
\end{proof}

We can also bound the mean of \(s(G)=d(G)-d^*(G)\), whose observed
decay is shown in the third pane of Figure~\ref{fig:dstar-histograms}. 

\begin{corollary}
\label{cor:mean-s}
Let \(G\) be uniformly distributed over connected Eulerian graph
isomorphism classes of order \(n\). For some \(c>0\),
\begin{equation}
 (1-o(1))\binom n3 2^{-n}\leq\mathbb E[s(G)]
 \leq\frac{n(n-1)}6\Pr(s(G)>0)=O(e^{-cn})
 \qquad(n\to\infty).
 \label{eq:mean-s-bounds}
\end{equation}
\end{corollary}

\begin{proof}
The upper bound follows from \(s(G)\leq d^*(G)\leq |E(G)|/3\).
The proof of Theorem~\ref{thm:typical-eulerian} bounds the exceptional
probability exponentially in \(n\), and the polynomial factor can be
absorbed by reducing the exponential rate.

For the lower bound, a triangle sharing no edge with another triangle
is a component of \(D_*(G)\) whose cover splits, contributing one to
\(s(G)\). For \(n\geq5\), fix its three vertices in a labelled even-degree graph and
put \(m=n-3\). Each outside vertex joins at most one of the three.
Requiring each of the three attachment counts to be even leaves
\[
 B_m=\frac18\bigl(4^m+(3+(-1)^m)2^m\bigr)
\]
choices, by averaging the eight parity characters. Each choice permits
\(2^{\binom{m-1}2}\) outside graphs with the required degree parities.
The probability that a fixed triple gives such a triangle is therefore
\[
 B_m2^{-3m}=2^{-n}+(3+(-1)^{n-3})2^{-2n+3}\sim2^{-n}.
\]
The cut estimate in Proposition~\ref{prop:qlc-asymptotic} applies
uniformly to prescribed degree parities, so the outside graph is
connected with probability \(1-o(1)\). The attachment is nonempty
with probability \(1-1/B_m\). Thus imposing connectedness leaves the
leading term unchanged. Summing over triples gives the labelled
lower bound. Dividing the labelled total by \(n!\) gives a lower
bound for the total over isomorphism classes, and
\(e_n\sim2^{\binom{n-1}2}/n!\) proves the lower bound.
\end{proof}

Counting all connections involves further weights. To isolate the
Eulerian contribution, sum over connected Eulerian graphs of order \(n\):
\[
 Z_n(u,v)=\sum_{\substack{G\text{ connected Eulerian}\\  |G|=n}}
                  u^{d^*(G)}v^{s(G)},
\]
where \(s(G)\) is the number of split components in
Theorem~\ref{thm:derived-cover}. Then \(Z_n(1,1)=e_n\), while
\(Z_n(2,2)\) and \(Z_n(2,1)\) count the Eulerian geometries and
star-preserving Eulerian geometries, respectively. By
Proposition~\ref{prop:coeu-counts}, the total counts including
co-Eulerian geometries are therefore \(Z_n(2,2)+o(e_n)\) and
\(Z_n(2,1)+o(e_n)\). Theorem~\ref{thm:typical-eulerian} identifies
the typical summand, but does not control the exponentially weighted
exceptional terms.

\subsection{Computational counts and asymptotic questions}
\label{ssec:graph-counts}
\label{sec:computations}

Tables~\ref{tab:graph-classes} and~\ref{tab:QLC-statistics} use the graph
classification counts obtained by enumerating connected
graphs up to isomorphism and applying the parity tests,
triangle-component counts and flatness criterion above. Connections
are counted on one labelled representative, before quotienting by
graph automorphisms. No individual connections need be generated.
The counts through nine vertices also agree with the solution archive.
At order thirteen, Eulerian enumeration captures every co-Eulerian
graph as well, giving the inner and star-preserving counts.

\begin{table}[tbp]
\centering
\begin{tabular}{r r r r r}
\toprule
\(n\)&Inner QLCs&QLCs&Flat QLCs&Star QLCs\\
\midrule
2&1&1&1&1\\
3&1&4&3&2\\
4&3&4&3&4\\
5&4&23&1&9\\
6&12&31&2&21\\
7&37&289&1&95\\
8&198&1230&4&478\\
9&1783&17741&1&4731\\
10&31075&180470&2&77024\\
11&1148633&5283235&1&2703573\\
12&86539453&295948615&5&192306329\\
13&12798435976&34870883156&1&27274476608\\
\bottomrule
\end{tabular}
\caption{QLC assignments on connected graph representatives over
\(\F_2\), before quotienting by automorphisms. The inner count is
also the number of braidings.}
\label{tab:QLC-statistics}
\end{table}

The graphs admitting QLCs form the union of the Eulerian and
co-Eulerian classes in Table~\ref{tab:graph-classes}; their intersection
contributes twice to the inner count in Table~\ref{tab:QLC-statistics}.
By Proposition~\ref{prop:complement}, connected co-Eulerian graphs
outside the intersection are odd-degree. They occur only at even
orders, and complementation pairs them with Eulerian co-Eulerian
graphs whenever the complements are connected. Differences between
these two counts arise from disconnected complements. At odd order,
the graphs admitting QLCs are precisely the Eulerian graphs.
By Theorem~\ref{thm:flat-alpha}, the flat star-preserving count equals
the flat QLC count except at \(n=3\), where it is one rather than three.

The program \texttt{qlc-graph-counts.py} streams graph representatives
and saves counts in resumable batches. Appendix~\ref{app:computational-classification}
describes the programs; commands and output formats are given in the
separate supplementary documentation.

\begin{figure}[tbp]
\centering
\begin{tikzpicture}
\begin{groupplot}[
 group style={group size=3 by 1,horizontal sep=.055\textwidth},
 width=.36\textwidth,height=5.1cm,
 xlabel={$n$},ylabel={Probability},ymode=log,
 xmin=2,xmax=12,xtick={2,4,6,8,10,12},
 ymin=0.000000001,ymax=1.2,grid=both,
 tick label style={font=\small},label style={font=\small},
 title style={font=\small}]
\nextgroupplot[title={QLC existence},xmax=13,xtick={3,5,7,9,11,13},
 legend style={font=\scriptsize,at={(0.04,0.04)},anchor=south west,
 draw=none,fill=white,cells={anchor=west}}]
\addplot[blue,thick,mark=*] coordinates {
 (2,1) (3,0.5) (4,0.333333333333) (5,0.190476190476)
 (6,0.0892857142857) (7,0.0433763188746) (8,0.0171808941261)
 (9,0.00682549410142) (10,0.00265017811098) (11,0.00114098078409) (12,0.000527486184454) (13,0.000254260554935305)};
\addlegendentry{\(P_n\)}
\addplot[black,dashed,thick,domain=2:13,samples=60] {2^(1-x)};
\addlegendentry{\(2^{1-n}\)}
\nextgroupplot[title={Flatness among star QLCs},title style={font=\footnotesize},ylabel={},
 xmax=13,xtick={3,5,7,9,11,13},ymin=0.00000000001,
 legend style={font=\scriptsize,at={(0.04,0.04)},anchor=south west,
 draw=none,fill=white,cells={anchor=west}}]
\addplot[purple,thick,mark=*] coordinates {
 (2,1) (3,0.5) (4,0.75) (5,0.111111111111)
 (6,0.0952380952381) (7,0.0105263157895) (8,0.00836820083682)
 (9,0.000211371803001) (10,2.59659326963e-05) (11,3.69880894653e-07)
 (12,2.60001843205e-08) (13,3.666431493342334e-11)};
\addlegendentry{\(F_n\)}
% F_n^{ref} below, evaluated in logarithmic form to avoid large intermediate powers.
\addplot[orange!85!black,dashdotted,thick,domain=5:13,samples=100,forget plot]
 {2^(ln(3*sqrt(2*pi*x)/4)/ln(2)
      +x*(ln(x)-1)/ln(2)-(x-1)*(x-2)/2)};
\nextgroupplot[title={Star preservation},ylabel={},ymode=normal,
 xmax=13,xtick={3,5,7,9,11,13},
 ymin=0,ymax=1.05,ytick={0,0.25,0.5,0.75,1},
 legend style={font=\scriptsize,at={(0.96,0.04)},anchor=south east,
 draw=none,fill=white,cells={anchor=west}}]
\addplot[teal,thick,mark=*] coordinates {
 (2,1)
 (3,0.5)
 (4,1)
 (5,0.391304347826)
 (6,0.677419354839)
 (7,0.328719723183)
 (8,0.388617886179)
 (9,0.266670424441)
 (10,0.426796697512) (11,0.511726811319) (12,0.649796347248) (13,0.782156175568701)};
\addlegendentry{\(S_n\)}
\addplot[black,dashed,thick,domain=2:13,samples=2,forget plot] {1};
\end{groupplot}
\end{tikzpicture}
\caption{Left: ratio of graph isomorphism classes admitting QLCs through \(n=13\), against its theoretical
asymptote \(2^{1-n}\).  Middle: 
ratio of flat connections among star-preserving QLCs and dashed-dotted the Stirling
reference \(F_n^{\mathrm{ref}}\). Right: ratio of QLCs that are $*$-preserving, with possible asymptote to 1.}
\label{fig:QLC-frequencies}
\end{figure}
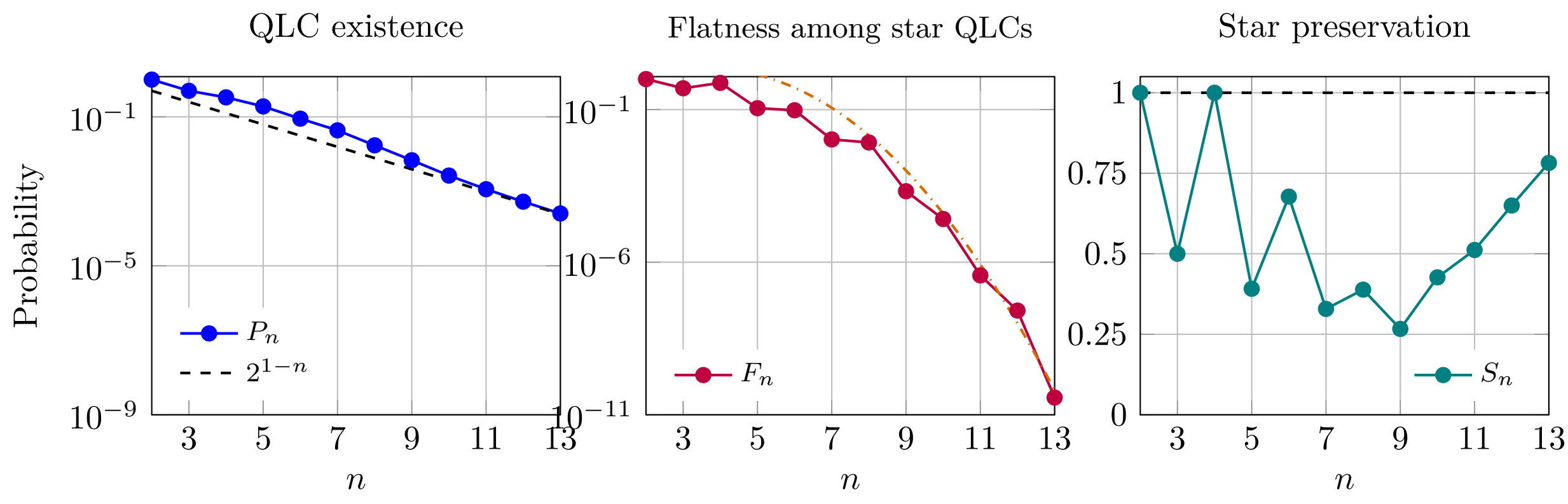

For the middle pane of Figure~\ref{fig:QLC-frequencies}, choose uniformly from all
star-preserving QLC assignments at that \(n\), and put
\[
 F_n=\frac{\#\{\text{flat star-preserving QLCs}\}}
 {\#\{\text{star-preserving QLCs}\}}.
\]
For the right pane, choose uniformly from all QLC assignments and put
\[
 S_n=\frac{\#\{\text{star-preserving QLCs}\}}{\#\{\text{QLCs}\}}.
\]
These two sampling rules weight graphs by their respective numbers of
connections, before quotienting by automorphisms. By
Theorem~\ref{thm:derived-cover}, \(S_n\) is the QLC-count-weighted
average of \(2^{-s(G)}\) over the graphs admitting QLCs. Although not clear from the limited data, the plot is consistent with a trend towards $1$ at large $n$. 

Next, approximating the star-preserving QLC count by \(2e_n\) (see below) and the
flat count by \(1.5\) gives the Stirling reference curve
\[
 F_n^{\mathrm{ref}}=
 \frac34\sqrt{2\pi n}\left(\frac{n}{\emath}\right)^n
 2^{-(n-1)(n-2)/2}.
\]
The numerator averages the tentative counts one at odd order and two
at even order, with exceptional counts at \(n=4,8,12\). This is a reference curve rather than an asymptotic result. The numerator of \(F_n\) depends
on Conjecture~\ref{conj:odd-flat} at odd order and on the outstanding classification of flat co-Eulerian geometries at even order.

\begin{proposition}[Decay of the flatness frequency]
\label{prop:flat-frequency}
The fraction \(F_n\) of star-preserving QLCs which are flat decays
exponentially in \(n^2\). More precisely, for \(n\geq4\),
\[
 \frac{1}{2(1+2^{1/3})^{\binom n2}}
 \leq F_n\leq\frac{1+b_n}{e_n}
 \leq2^{-n^2/4+O(n\log n)}.
\]
\end{proposition}

\begin{proof}
By Theorem~\ref{thm:flat-alpha}, the Eulerian geometries contribute
one flat connection, on \(C_n\), and each co-Eulerian graph contributes
at most one. Hence the flat count lies between one and \(1+b_n\).
There are at least \(e_n\) star-preserving QLCs, so
Proposition~\ref{prop:qlc-asymptotic} gives the upper bound.

Each component of \(D_*(G)\) contains at least the three edges of a
triangle, so \(d^*(G)\leq|E|/3\). A graph therefore contributes at
most \(2\cdot2^{|E|/3}\) star-preserving QLCs, counting both braidings.
Bounding graph isomorphism classes by labelled graphs gives
\[
 Q_{*,n}\leq2\sum_{G\ {\rm labelled}}2^{|E(G)|/3}
          =2(1+2^{1/3})^{\binom n2}.
\]
The existence of the flat cycle connection proves the lower bound.
\end{proof}

\begin{figure}[tbp]
\centering
\begin{tikzpicture}
\begin{groupplot}[
 group style={group size=3 by 1,horizontal sep=.065\textwidth},
 width=.36\textwidth,height=5.4cm,
 xlabel={$n$},grid=both,
 tick label style={font=\small},label style={font=\small},
 title style={font=\footnotesize}]
\nextgroupplot[title={Eulerian among co-Eulerian},
 xmin=4,xmax=12,xtick={4,6,8,10,12},
 ymin=.45,ymax=.55,ytick={.45,.5,.55},
 legend style={font=\scriptsize,at={(.04,.04)},anchor=south west,
 draw=none,fill=white,cells={anchor=west}}]
\addplot[orange!85!black,thick,mark=*,mark size=1pt] coordinates {
 (4,0.5) (6,0.5) (8,0.5) (10,0.489795918367347) (12,0.489230769230769)};
\addlegendentry{\(c_n/b_n\)}
\addplot[black,dashed,thick,domain=4:12,samples=2,forget plot] {.5};
\nextgroupplot[title={Co-Eulerian probability},
 xmin=3,xmax=13,xtick={3,5,7,9,11,13},
 ymode=log,ymin=0.000000001,ymax=1,
 legend style={font=\scriptsize,at={(.04,.04)},anchor=south west,
 draw=none,fill=white,cells={anchor=west}}]
\addlegendimage{empty legend}
\addlegendentry{\(b_n/(e_n+o_n)\)}
\addplot[blue,thick,mark=*,mark size=1pt] coordinates {
 (4,1) (6,0.4) (8,0.0732984293193717) (10,0.00157804901613475) (12,3.75551942912777e-06)};
\addlegendentry{Even \(n\)}
\addplot[purple,solid,thick,mark=square*,mark size=1pt] coordinates {
 (9,0.000561167227833895) (11,6.09423781108907e-06) (13,8.43853117004969e-09)};
\addlegendentry{Odd \(n\)}
\addplot[black,dashed,thick,domain=5:12,samples=100,forget plot]
 {2^(-x*x/3+3*x-6)};
\nextgroupplot[title={Star QLCs / Eulerian},
 xmin=3,xmax=13,xtick={3,5,7,9,11,13},
 ymin=1.8,ymax=4.2,ytick={2,2.5,3,3.5,4},
 legend style={font=\scriptsize,at={(.97,.97)},anchor=north east,
 draw=none,fill=white,cells={anchor=west}}]
\addplot[teal,thick,mark=*,mark size=1pt] coordinates {
 (3,2)
 (4,4)
 (5,2.25)
 (6,2.625)
 (7,2.56756756757)
 (8,2.59782608696)
 (9,2.65488215488)
 (10,2.48256301167)
 (11,2.35374525738)
 (12,2.22218935462)
 (13,2.13107889818)};
\addlegendentry{\(Q_{*,n}/e_n\)}
\addplot[black,dashed,thick,domain=3:13,samples=2,forget plot] {2};
\end{groupplot}
\end{tikzpicture}
\caption{Left: Eulerian among co-Eulerian graphs at even orders.
Middle: ratio of co-Eulerian graphs among those that admit a QLC and dashed reference curve
$2^{-n^2/3+3n-6}$ for even $n$.
Right: star-preserving QLCs per Eulerian graph, with possible asymptote to 2.}
\label{fig:eulerian-ratios}
\end{figure}
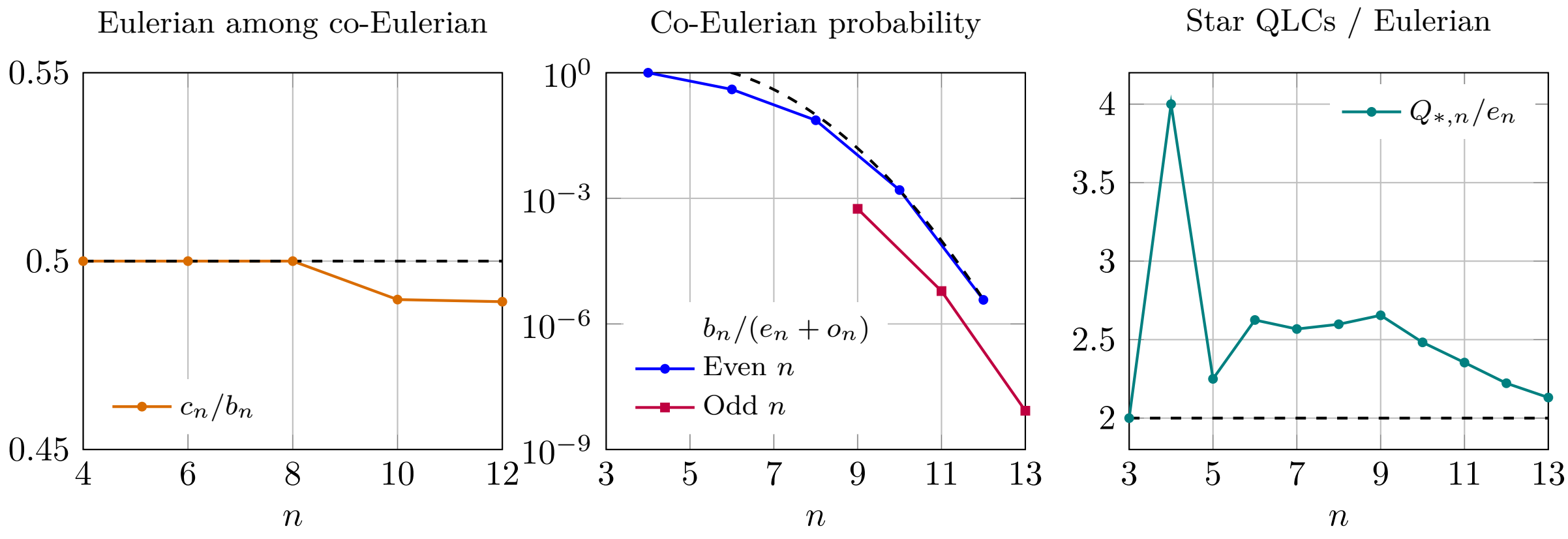

Let \(o_n\) count connected co-Eulerian graphs with all degrees odd (i.e. not Eulerian, by Proposition~\ref{prop:complement}). Let  \(c_n\) count those graphs that are both Eulerian and co-Eulerian, including \(c_1=1\), so  \(b_n=o_n+c_n\). The left pane of Figure~\ref{fig:eulerian-ratios}
shows \(c_n/b_n\) at even orders from Table~\ref{tab:graph-classes}. We see that this ratio is around \(1/2\) and conjecture that this could be its asymptotic limit along even orders. This is because 
complementation pairs the two types when both graphs are connected (Proposition~\ref{prop:complement}) which could be more likely for large $n$. There is also an exact relation
between their counts.
In formal power series over \(\mathbb Q\), define
\[
 C(x)=\prod_{m\ {\rm odd}}(1-x^m)^{-c_m}.
\]
Then
\begin{equation}
\label{eq:coeu-parity-series}
 \prod_{m\ {\rm even}}(1-x^m)^{-(o_m-c_m)}
 =\frac{C(x)+C(-x)}2.
\end{equation}
Indeed, a graph with \(A^2=0\) is a disjoint union of connected
Eulerian co-Eulerian graphs. At even order complementation bijects
these graphs with those satisfying \(A^2=I\), whose connected
components all have odd degree. Taking the even part of the first
Euler product and cancelling its even-order factors proves the
identity. Thus the differences at even order are determined by the
odd-order intersection counts. Table~\ref{tab:graph-classes} gives
\(c_3=c_5=c_7=0\), \(c_9=1\), \(c_{11}=7\), \(c_{13}=108\).
Expanding \eqref{eq:coeu-parity-series} for these values then gives
\[o_n=c_n+c_{n-1}, \quad n=2,4,\cdots, 16\]
 \[ o_{18}=c_{18}+c_{17}+1,\qquad o_{20}=c_{20}+c_{19}+6,\qquad
 o_{22}=c_{22}+c_{21}+129 \]
 For example, the extra term at $n=20$ comes from \(
c_9c_{11}-\binom{c_9+1}{2}=6\). In general, for even $n\geq4$, 
\[ {c_n\over b_n}={c_n\over 2c_n+c_{n-1}+\mathcal P_n(c_3,c_5,\ldots,c_{n-3})}.\]
Here $\mathcal P_n$ are polynomials determined from \eqref{eq:coeu-parity-series}. The limit $c_n/b_n\to1/2$ is equivalent to \(o_n-c_n=o(c_n)\) along even
orders. Meanwhile, the middle of Figure~\ref{fig:eulerian-ratios} shows the probability $b_n/(e_n+o_n)$ that a uniformly chosen connected graph admitting a QLC is co-Eulerian.
It tends to zero by Proposition~\ref{prop:qlc-asymptotic}. Finally, the right of Figure~\ref{fig:eulerian-ratios} suggests:

\begin{conjecture}
\label{conj:star-asymptotic} The total number  \(Q_{*,n}\) of star-preserving QLCs on connected graphs of order $n$ obeys \(Q_{*,n}\sim2e_n\) as $n\to \infty$. \end{conjecture}

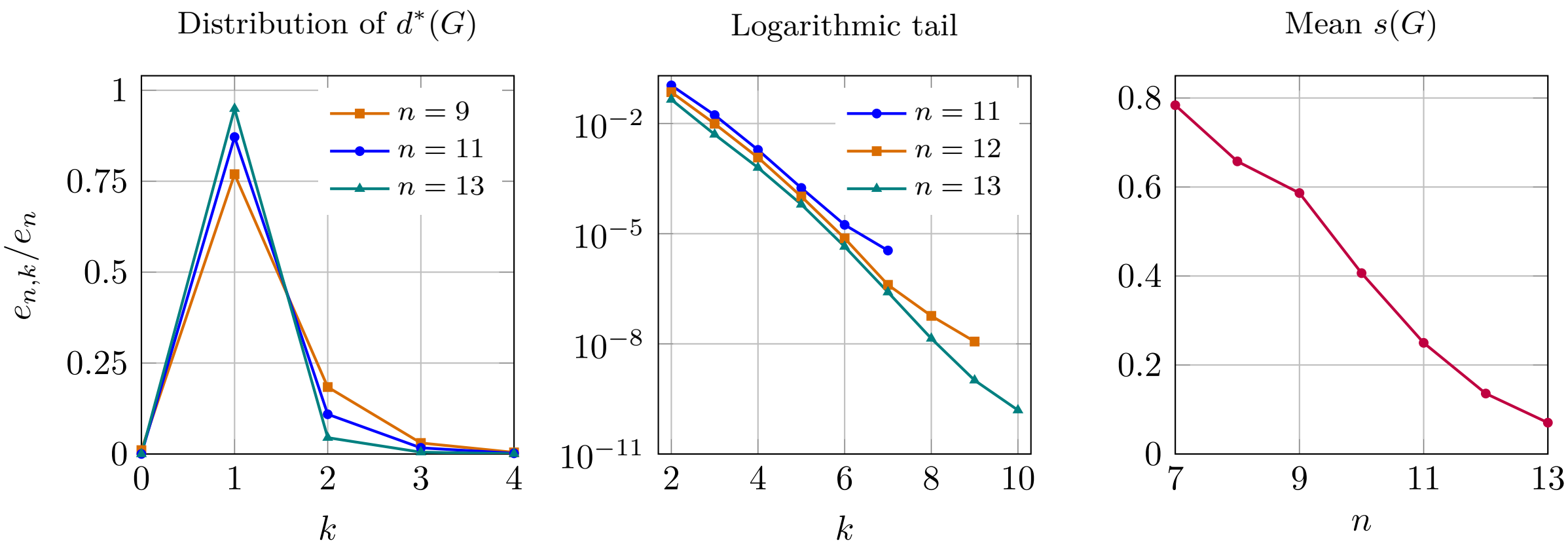
\begin{figure}[tbp]
\centering
\begin{tikzpicture}
\begin{groupplot}[
 group style={group size=3 by 1,horizontal sep=.09\textwidth},
 width=.33\textwidth,height=5.4cm,
 xlabel={$k$},grid=major,
 tick label style={font=\small},label style={font=\small},
 title style={font=\footnotesize}]
\nextgroupplot[title={Distribution of $d^*(G)$},
 ylabel={$e_{n,k}/e_n$},
 xlabel={$k$},xmin=0,xmax=4,xtick={0,1,2,3,4},
 ymin=0,ymax=1.04,ytick={0,.25,.5,.75,1},
 legend style={font=\scriptsize,at={(.97,.97)},anchor=north east,
 draw=none,fill=white,cells={anchor=west}}]
\addplot[orange!85!black,thick,solid,mark=square*,mark size=1pt] coordinates {
 (0,0.010662177328844) (1,0.769360269360269) (2,0.184062850729517) (3,0.0303030303030303) (4,0.00448933782267116)};
\addlegendentry{$n=9$}
\addplot[blue,thick,solid,mark=*,mark size=1pt] coordinates {
 (0,0.000199368636962771) (1,0.871475136380336) (2,0.109240953974575) (3,0.0169715816984815) (4,0.00191620248888672)};
\addlegendentry{$n=11$}
\addplot[teal,thick,solid,mark=triangle*,mark size=1pt] coordinates {
 (0,6.36562161503656e-07) (1,0.949336189071257) (2,0.0449294227771803) (3,0.00503743220382092) (4,0.000629004128553563)};
\addlegendentry{$n=13$}
\nextgroupplot[title={Logarithmic tail},
 xmin=1.7,xmax=10.3,xtick={2,4,6,8,10},
 ymode=log,ymin=0.00000000001,ymax=.2,
 legend style={font=\scriptsize,at={(.97,.97)},anchor=north east,
 draw=none,fill=white,cells={anchor=west}}]
\addplot[blue,thick,mark=*,mark size=1pt] coordinates {
 (2,0.109240953974575) (3,0.0169715816984815) (4,0.00191620248888672) (5,0.000175862291119999) (6,1.74121080316831e-05) (7,3.48242160633661e-06)};
\addlegendentry{$n=11$}
\addplot[orange!85!black,thick,mark=square*,mark size=1pt] coordinates {
 (2,0.0714445146708666) (3,0.00984794993543268) (4,0.00118473576484385) (5,0.000103213427341214) (6,7.48794233286011e-06) (7,4.04441329706951e-07) (8,5.77773328152787e-08) (9,1.15554665630557e-08)};
\addlegendentry{$n=12$}
\addplot[teal,thick,mark=triangle*,mark size=1pt] coordinates {
 (2,0.0449294227771803) (3,0.00503743220382092) (4,0.000629004128553563) (5,6.25985087758382e-05) (6,4.44632458113748e-06) (7,2.55265567894003e-07) (8,1.39860840688787e-08) (9,1.0157491223208e-09) (10,1.56269095741661e-10)};
\addlegendentry{$n=13$}
\nextgroupplot[title={Mean $s(G)$},
 xlabel={$n$},xmin=7,xmax=13,xtick={7,9,11,13},
 ymin=0,ymax=.85,ytick={0,.2,.4,.6,.8}]
\addplot[purple,thick,mark=*,mark size=1pt] coordinates {
 (7,0.783783783783784) (8,0.657608695652174)
 (9,0.586419753086420) (10,0.406433313994714)
 (11,0.249794101822525) (12,0.135854858625338)
 (13,0.0703697094933164)};
\end{groupplot}
\end{tikzpicture}
\caption{Connected Eulerian graphs: distribution of \(d^*\), its logarithmic tail,
and mean \(s=d-d^*\).}
\label{fig:dstar-histograms}
\end{figure}

We count each graph up to isomorphism and both braidings when they occur.  The right of Figure~\ref{fig:eulerian-ratios} shows $Q_{*,n}/e_n$ through $n=13$ and the conjecture is that this tends to 2. 
By Proposition~\ref{prop:coeu-counts}, the conjecture is equivalent to
\(Z_n(2,1)\sim2e_n\). Writing \(e_{n,k}\) for the number of connected
Eulerian graphs with \(d^*(G)=k\), the exact identity
\[
 \frac{Z_n(2,1)}{e_n}
 =2-\frac{e_{n,0}}{e_n}
   +\frac1{e_n}\sum_{k\ge2}(2^k-2)e_{n,k}
\]
and Theorem~\ref{thm:typical-eulerian} reduce this form of the conjecture  to
\(\sum_{k\ge2}(2^k-2)e_{n,k}=o(e_n)\).
Thus, the remaining task is to control the weighted contribution of
graphs with large \(d^*\). Figure~\ref{fig:dstar-histograms} shows the normalised distributions of \(d^*\) increasingly concentrated at $k=1$ as
$n$ increases in agreement with Theorem~\ref{thm:typical-eulerian}.   The stronger estimate \(Z_n(2,2)\sim2e_n\) would imply both 
Conjecture~\ref{conj:star-asymptotic} and \(S_n\to1\) but its proof would also require controlling the factor \(2^{s(G)}\) from split components
of the derived cover. The mean \(s(G)\), shown in the third pane,
decreases from \(0.784\) at \(n=7\) to \(0.0704\) at \(n=13\).
Corollary~\ref{cor:mean-s} gives exponential decay, with a
lower bound asymptotic to \(\binom n3 2^{-n}\).

\section{Maps between graph quantum Riemannian geometries over \(\F_2\)}
\label{sec:QRG-maps}

We study graph QRG maps over a field \(k\), then specialise to \(\F_2\).

\subsection{Digraph maps and pullback over a field}

Let \(k\) be a field, and let \((W,F)\) and \((V,E)\) be finite
digraphs. The first-order correspondence of \cite[Thm.~3.1]{Ma:gra}
uses vertex maps \(\psi:W\to V\) satisfying
\begin{equation}
\label{eq:digraph-map-condition}
a\to b\in F\quad\Longrightarrow\quad
\psi(a)=\psi(b)\ \text{or}\ \psi(a)\to\psi(b)\in E.
\end{equation}
The function pullback is \({}_0\phi_\psi:k(V)\to k(W)\),
\({}_0\phi_\psi(f)=f\circ\psi\), and its one-form pullback \(\phi_\psi={}_1\phi_\psi\) is
\begin{equation}
\label{eq:QRG-form-pullback}
\phi_\psi(\omega_{x,y})=
\sum_{\substack{a\to b\in F\\\psi(a)=x,\ \psi(b)=y}}
\omega_{a,b}.
\end{equation}
It obeys \(\phi_\psi(f\omega h)={}_0\phi_\psi(f)\phi_\psi(\omega){}_0\phi_\psi(h)\)
and \(\phi_\psi\extd_V=\extd_W{}_0\phi_\psi\).

Write \(\Phi_\psi=\phi_\psi\otimes\phi_\psi\) for the induced map
on balanced two-tensors. Explicitly,
\[
\Phi_\psi(\omega_{x,y}\otimes_{k(V)}\omega_{y,z})
=\sum_{\substack{a\to b\to c\text{ in }W\\
                 \psi(a)=x,\ \psi(b)=y,\ \psi(c)=z}}
  \omega_{a,b}\otimes_{k(W)}\omega_{b,c}.
\]

\subsection{Metric, connection and star compatibility}

Assume now that both graphs are bidirected and have nondegenerate graph
metrics \(g_V,g_W\). If star structures are included, equip \(k\)
with an involution \(\lambda\mapsto\overline\lambda\), which may be
the identity, and use the pointwise star on functions and
\(\omega_{x,y}^*=-\omega_{y,x}\) on one-forms. QLCs $\nabla$ are given via bimodule maps $(\sigma,\alpha)$. 

\begin{definition}[A map of graph QRGs]
\label{def:QRG-map}
For \(\mathcal Q_W=(W,g_W,\sigma_W,\alpha^W)\) and
\(\mathcal Q_V=(V,g_V,\sigma_V,\alpha^V)\), a geometric or QRG map
\(\psi:\mathcal Q_W\to\mathcal Q_V\) is a digraph map
\eqref{eq:digraph-map-condition} whose induced pullback
satisfies
\begin{equation}
\label{eq:QRG-morphism-equations}
\Phi_\psi(g_V)=g_W,\qquad
       \Phi_\psi\nabla_V=\nabla_W\phi_\psi.
\end{equation}
For star-preserving QRGs, we also require the pullbacks to respect star. The geometric arrow \(W\to V\) and the one-form arrow
\(\Omega^1(V)\to\Omega^1(W)\) have opposite directions; throughout
this section the quiver uses the geometric direction.
\end{definition}

\begin{proposition}[Automatic star and braiding compatibility]
\label{prop:QRG-map-star-sigma}
The pullbacks of a bidirected graph map respect star.
Connection compatibility in \eqref{eq:QRG-morphism-equations} implies
\begin{equation}
\label{eq:QRG-sigma-naturality}
\Phi_\psi\sigma_V=\sigma_W\Phi_\psi.
\end{equation}
Under metric preservation in \eqref{eq:QRG-morphism-equations}, connection compatibility is equivalent to
\eqref{eq:QRG-sigma-naturality} together with
\begin{equation}
\label{eq:QRG-alpha-naturality}
\Phi_\psi\alpha^V=\alpha^W\phi_\psi.
\end{equation}
\end{proposition}

\begin{proof}
The coefficients in \eqref{eq:QRG-form-pullback} are \(0\) and \(1\).
Reversing a lifted edge exchanges the fibres of its endpoints, so
\(\phi_\psi(\omega_{x,y}^*)=\phi_\psi(\omega_{x,y})^*\).
The same holds for function pullback and tensor reversal.

Apply connection compatibility to \(\omega f\) and use both right
Leibniz rules. After cancelling the terms involving \(\nabla\omega\),
we obtain
\[
\Phi_\psi\sigma_V(\omega\otimes \extd_Vf)
=\sigma_W\bigl(\phi_\psi(\omega)\otimes \extd_W{}_0\phi_\psi(f)\bigr).
\]
The tensors \(\omega\otimes \extd_Vf\) span the balanced tensor space:
one-forms are generated by \(a\,\extd b\), and \(a\) can be absorbed
into the first factor. This proves \eqref{eq:QRG-sigma-naturality}.
The next proposition shows that metric preservation forbids collapsed
edges, and consequently \(\phi_\psi(\theta_V)=\theta_W\).
Substituting the inner connection formula into
\eqref{eq:QRG-morphism-equations}, and using
\eqref{eq:QRG-sigma-naturality}, leaves precisely
\eqref{eq:QRG-alpha-naturality}. Conversely those two identities and
\(\phi_\psi(\theta_V)=\theta_W\) give connection compatibility.
\end{proof}

\begin{proposition}[The metric condition is a graph immersion condition]
\label{prop:QRG-metric-immersion}
For nonzero directed-edge metric coefficients, \(\Phi_\psi(g_V)=g_W\)
is equivalent to the following conditions:
\begin{enumerate}[label=\textup{(\roman*)}]
\item every edge of \(W\) maps to an edge of \(V\), with distinct endpoints;
\item \(\psi:N_W(b)\to N_V(\psi(b))\) is injective for every \(b\in W\);
\item \(g^W_{ab}=g^V_{\psi(a),\psi(b)}\) for every directed edge \(a\to b\).
\end{enumerate}
For Euclidean metric on both graphs, condition \textup{(iii)} is automatic.
\end{proposition}

\begin{proof}
On a lifted two-path \((a,b,c)\), the coefficient of \(\Phi_\psi(g_V)\)
is zero if \(\psi(a)=\psi(b)\); otherwise it is
\[
g^V_{\psi(a),\psi(b)}\,\delta_{\psi(a),\psi(c)}.
\]
The coefficient in \(g_W\) is \(g^W_{ab}\delta_{a,c}\).
On return paths this forces noncollapse and equality of edge weights.
On a non-return path with \(a,c\in N_W(b)\), nonzero weights then
force \(\psi(a)\ne\psi(c)\), exactly the neighbour-set injectivity.
Conversely these conditions make every two-path coefficient agree.
\end{proof}

\begin{remark}
Naturality of the inverse-metric pairing is an additional stronger
requirement: for Euclidean metrics it requires local surjectivity
on neighbour sets, giving graph coverings.
\end{remark}

\begin{proposition}[Composition and the degree-two calculus]
\label{prop:QRG-map-composition}
QRG maps compose and include identity maps. The pullbacks
\({}_0\phi_\psi\) and \(\phi_\psi={}_1\phi_\psi\) extend to a map
of differential graded algebras \(\Omega_{\min}(V)\to\Omega_{\min}(W)\),
with degree-two part
\({}_2\phi_\psi:\Omega^2_{\min}(V)\to\Omega^2_{\min}(W)\).
The curvatures obey
\[
R_{\nabla_W}\phi_\psi=({}_2\phi_\psi\otimes\phi_\psi)R_{\nabla_V}.
\]
In particular, a map \(\mathcal Q_W\to\mathcal Q_V\) with flat
\(\mathcal Q_V\) has flat \(\mathcal Q_W\).
\end{proposition}

\begin{proof}
Pullbacks compose in the opposite order to vertex maps, and preserve
both equations \eqref{eq:QRG-morphism-equations}. For each pair \(x,z\),
noncollapse gives
\[
\Phi_\psi(P^{\mathrm{sum}}_{x,z})
=\sum_{\substack{a,c\in W\\\psi(a)=x,\ \psi(c)=z}}
  P^{\mathrm{sum}}_{a,c}.
\]
Each summand is a defining relation after applying the wedge. Thus
\(\Phi_\psi\) descends to degree two. Multiplicatively extending the
one-form pullback preserves the ideal generated by these relations;
compatibility with the differential also preserves its differential
closure. Hence the extension is defined in all degrees of this
\(\Omega_{\min}\) calculus. The inner formula for the
differential and \(\phi_\psi(\theta_V)=\theta_W\) give differential
compatibility, and substitution in the curvature formula gives the
stated identity. Finally, local injectivity implies that any directed
edge form \(\omega_{a,b}\) is
\(\delta_a\phi_\psi(\omega_{\psi(a),\psi(b)})\).
The pullback image therefore generates \(\Omega^1(W)\) as a left
\(k(W)\)-module. Since curvature is left-linear, flatness of \(V\)
forces flatness of \(W\).
\end{proof}

\subsection{Morphism tests over \(\F_2\)}

For a graph with vertex set \(V\), write
\[
 \alpha^V(\omega_{xy})=a^V_{xy}P^{\mathrm{sum},V}_{xy},
\]
with \(a^V_{xy}=0\) on nonedges (including \(x=y\)) and on edges
in no triangle, and use the analogous notation for \(W\). Eulerian and co-Eulerian
QRGs refer to their respective braidings in
Theorem~\ref{thm:general-F2-classification}; a graph supporting both
has two types of geometry.

\begin{theorem}[Binary QRG morphism criterion]
\label{thm:binary-QRG-map}
Let \(\mathcal Q_W,\mathcal Q_V\) be connected star-preserving QRGs over $\F_2$
on graphs with at least two vertices. A vertex map \(\psi:W\to V\)
is a QRG map precisely when:
\begin{enumerate}[label=\textup{(\roman*)}]
\item it is a noncollapsing, locally injective graph homomorphism;
\item the two QRGs have the same type and:
      \begin{itemize}
      \item for Eulerian QRGs, every neighbour-set map is bijective;
      \item for co-Eulerian QRGs, for every distinct \(a,c\) with
            \(N_W(a)\cap N_W(c)\ne\varnothing\), the map
            \[
            \psi:N_W(a)\cap N_W(c)\longrightarrow
            N_V(\psi(a))\cap N_V(\psi(c))
            \]
            is bijective;
      \end{itemize}
\item for every two-path \((a,b,c)\) in \(W\),
      \begin{equation}
      \label{eq:binary-map-alpha-test}
      a^V_{\psi(a),\psi(c)}=a^W_{ac}.
      \end{equation}
\end{enumerate}
Thus maps between Eulerian QRGs are coverings subject to the
\(\alpha\)-equation; maps between co-Eulerian QRGs can also be proper
immersions. There are no maps between the two types.
\end{theorem}

\begin{proof}
Condition \textup{(i)} is Proposition~\ref{prop:QRG-metric-immersion}.
Write \(q_W,q_V\) for the two braiding parameters.
Fix \(a\in W\) and a neighbour \(b\), with images \(x,y\).
Apply braiding naturality to the return path in \(V\) \((x,y,x)\),
and restrict to return paths based at \(a\). The coefficient of the
return-block path sum is \(q_V+1\) on one side and \(q_W+1\) on
the other, forcing \(q_W=q_V=:q\). If a neighbour \(y\in N_V(x)\)
has no lift at \(a\), its coefficient is instead \(q_V+1\) on one
side and zero on the other. For \(q=0\), all neighbours must therefore
lift. With local injectivity, this is the covering condition.
For \(q=1\), return blocks are identities and impose no further condition.

On a non-return block in \(W\) with endpoints \(a,c\), local
injectivity identifies its distinct middle vertices with a subset of
the middle vertices for \((\psi(a),\psi(c))\) in \(V\).
The braiding adds \(q\) times the complete block path sum.
For \(q=0\) this addition vanishes. For \(q=1\), naturality for
each middle vertex in \(V\) requires that it lift in the given block:
the middle-vertex subset must be the entire set. Blocks with no
\(W\)-paths contribute zero and impose no condition.
These arguments also prove sufficiency for braiding naturality.

Finally, the coefficient of a two-path \((a,b,c)\) in
\(\Phi_\psi\alpha^V\) is the left-hand side of
\eqref{eq:binary-map-alpha-test}; in \(\alpha^W\phi_\psi\) it is
the right-hand side. Checking all edge forms in \(V\) amounts to these
coefficient equations. Proposition~\ref{prop:QRG-map-star-sigma}
now proves the equivalence with connection compatibility. 
\end{proof}

Write \(\alpha_\triangle\) for the map with \(a_{xy}=1\) on every
edge in a triangle and zero elsewhere. Thus \(\alpha_\triangle\in V_*(G)\).

\begin{proposition}[Triangle-lifting coverings]
\label{prop:triangle-cover}
Let \(\mathcal Q_V=(V,\sigma_0,\alpha_\triangle)\) be the Eulerian
QRG on a connected graph \(G=(V,E)\) with at least two vertices.
A QRG map \(\psi:\mathcal Q_W\to\mathcal Q_V\) from a connected
star-preserving graph QRG on at least two vertices is precisely a
graph covering \(\psi:W\to V\) which lifts every triangle to a triangle, with
\(\sigma_W=\sigma_0\) and \(\alpha^W=\alpha_\triangle\).
If every set of at most four vertices of \(G\) has a common neighbour,
every such map is an isomorphism.
\end{proposition}

\begin{proof}
Theorem~\ref{thm:binary-QRG-map} makes the map a covering of Eulerian
geometries. Lift two consecutive edges of a triangle of \(G\) to
\(a\to b\to c\). The target coefficient in
\eqref{eq:binary-map-alpha-test} is one, so \(a^W_{ac}=1\).
In particular \(a\) and \(c\) are adjacent, and the lifted path
closes to a triangle. Every triangle in the source graph maps to a triangle,
so its edge coefficients are one. Conversely, triangle lifting
makes the two coefficients in \eqref{eq:binary-map-alpha-test}
agree on every two-path, proving the first assertion.

For the last assertion, fix \(o\in V\) and a lift \(\widetilde o\).
For \(x\ne o\), choose a common neighbour \(y\) of \(o,x\), and
define \(s(x)\) by lifting \(o\to y\to x\) from \(\widetilde o\).
This is independent of \(y\): given another choice \(y'\), choose
\(w\) adjacent to \(o,x,y,y'\). The triangles \(oyw\) and \(xyw\)
identify the endpoints of the lifts of \(oyx\) and \(owx\), and
similarly for \(y'\). Put \(s(o)=\widetilde o\).
For adjacent \(x,x'\), a common neighbour of \(o,x,x'\) and triangle
lifting show that \(s(x)\) and \(s(x')\) are adjacent. When one
endpoint is \(o\), the triangle \(oyx\) gives the same conclusion.
Thus \(s\) is a graph section of the covering. Its image contains
all neighbours of each of its vertices, by local bijectivity.
Connectedness of the source graph makes this image all of \(W\).
\end{proof}

\begin{corollary}
\label{cor:typical-map-rigidity}
For almost every connected Eulerian graph isomorphism class
\(G=(V,E)\), every QRG map from a connected star-preserving QRG on
at least two vertices into \(\mathcal Q_V=(V,\sigma_0,\alpha_\triangle)\)
is an isomorphism.
\end{corollary}

\begin{proof}
A uniformly chosen labelled even-degree graph \(G\) is obtained by
choosing all edges of \(G-z\) independently with probability \(1/2\)
and then fixing the edges incident at \(z\) by degree parity.
For any four-element set \(S\), choose \(z\notin S\). The events
that \(w\) is adjacent to every vertex of \(S\), for
\(w\notin S\cup\{z\}\), are independent, each of probability
\(1/16\). Hence
\[
 \Pr(\text{some four vertices have no common neighbour})
 \leq\binom n4(15/16)^{n-5}=o(1).
\]
For \(n\geq5\), the property for four-element sets implies it for
all smaller sets. The same Burnside estimate used in
Theorem~\ref{thm:typical-eulerian} transfers this relabelling-invariant
exceptional set to graph isomorphism classes. The result follows
from Proposition~\ref{prop:triangle-cover}.
\end{proof}

\subsection{Isomorphism classes}

\begin{theorem}[QRG isomorphisms]
\label{thm:QRG-isomorphism-clusters}
Let \(\mathcal Q_W=(W,\sigma_{q_W},\alpha^W)\) and
\(\mathcal Q_V=(V,\sigma_{q_V},\alpha^V)\) be connected
star-preserving QRGs with at least two vertices and the unique metrics
over \(\F_2\). A vertex map \(\psi:W\to V\) is a QRG isomorphism
if and only if it is a graph isomorphism, the QRGs have the same type, and
\[
 a^W_{ac}=a^V_{\psi(a),\psi(c)}
 \qquad\text{for every directed edge \((a,c)\) in a triangle.}
\]
Every graph automorphism preserves each QLC braiding \(\sigma_q\).
Consequently, after choosing one labelled representative of each
underlying graph isomorphism class, QRG isomorphisms occur only
between geometries of the same type on the same graph.
\end{theorem}

\begin{proof}
An isomorphism and its inverse give a graph isomorphism.
Theorem~\ref{thm:binary-QRG-map} then requires \(q_W=q_V\), and
its neighbour-set conditions hold automatically. The remaining
condition is equality of the transported \(\alpha\)-coefficients.
These conditions also hold for the inverse.
\end{proof}

\begin{corollary}[The orbit description of QRG isomorphism classes]
\label{cor:QRG-alpha-orbits}
Let \(G\) be connected and admit a QLC with braiding \(\sigma_q\).
The QRG isomorphism classes with this braiding form
\[
 V_*(G)\big/\Aut(G),
\]
where graph automorphisms act by permuting directed-edge coordinates.
The number of orbits is
\[
 \frac{1}{|\Aut(G)|}
 \sum_{\pi\in\Aut(G)}
 2^{\dim_{\F_2}\{\alpha\in V_*(G): \pi\alpha=\alpha\}}.
\]
\end{corollary}

\begin{proof}
The theorem identifies isomorphisms on a fixed representative with
its graph automorphisms carrying one \(\alpha\)-value to another.
The fixed-point sets are binary vector spaces, so the displayed
count follows from the orbit-counting lemma.
\end{proof}

For disconnected graphs, isomorphisms may permute isomorphic
components; the braidings must agree on corresponding components.

\subsection{The quiver through six vertices}
\label{ssec:quiver-six}

Choose one labelled representative of each connected graph of order
\(2\leq n\leq6\), and take all its star-preserving QLCs as objects.
Arrows are the vertex maps of Definition~\ref{def:QRG-map}, including
identities and isomorphisms. There are \(37\) objects on \(18\)
underlying graphs, with \(2534\) arrows: \(2224\) isomorphisms and
\(310\) non-isomorphisms. The isomorphism count includes the \(37\)
identities. There are \(36\) QRG isomorphism classes: the two
connections on a pair of triangles sharing one vertex, with precisely
one triangle parameter nonzero, are isomorphic.

\begin{figure}[tbp]
\centering
% BEGIN GENERATED QRG CO-EULERIAN
\definecolor{qrgC16a34a}{HTML}{16A34A}
\definecolor{qrgC334155}{HTML}{334155}
\definecolor{qrgC64748b}{HTML}{64748B}
\definecolor{qrgC7c3aed}{HTML}{7C3AED}
\definecolor{qrgC94a3b8}{HTML}{94A3B8}
\begin{tikzpicture}[x=.008cm,y=-.008cm]
\path[use as bounding box] (0,0) rectangle (2000,1110);
% Q04 -> Q06: 24 maps
\draw[->,>=stealth,draw=qrgC7c3aed,line width=.55pt] (965.020,193.990) .. controls (765.670,185.930) and (568.770,403.200) .. (352.030,286.830);
\node[anchor=center,text=qrgC7c3aed,font={\fontsize{9}{10}\selectfont},inner sep=0pt] at (591.990,261.860) {24};
% Q04 -> Q30: 32 maps
\draw[->,>=stealth,draw=qrgC7c3aed,line width=.55pt] (1014.890,192.690) .. controls (1215.780,174.020) and (1424.220,295.980) .. (1625.110,277.310);
\node[anchor=center,text=qrgC7c3aed,font={\fontsize{9}{10}\selectfont},inner sep=0pt] at (1373.720,282.310) {32};
% Q20 -> Q30: 48 maps
\draw[->,>=stealth,draw=qrgC7c3aed,line width=.55pt] (1650.000,815.000) .. controls (1650.000,651.670) and (1650.000,463.330) .. (1650.000,300.000);
\node[anchor=center,text=qrgC7c3aed,font={\fontsize{9}{10}\selectfont},inner sep=0pt] at (1617.000,504.800) {48};
% Q01 -> Q06: 12 maps
\draw[->,>=stealth,draw=qrgC7c3aed,line width=.55pt] (950.640,620.340) .. controls (752.760,546.600) and (532.760,426.600) .. (350.020,289.970);
\node[anchor=center,text=qrgC7c3aed,font={\fontsize{9}{10}\selectfont},inner sep=0pt] at (562.210,480.320) {12};
% Q01 -> Q07: 12 maps
\draw[->,>=stealth,draw=qrgC7c3aed,line width=.55pt] (948.320,629.840) .. controls (748.810,605.120) and (672.960,324.910) .. (469.420,280.350);
\node[anchor=center,text=qrgC7c3aed,font={\fontsize{9}{10}\selectfont},inner sep=0pt] at (699.430,387.120) {12};
% Q01 -> Q04: 8 maps
\draw[->,>=stealth,draw=qrgC7c3aed,line width=.55pt] (990.000,593.000) .. controls (990.000,488.330) and (990.000,341.670) .. (990.000,220.000);
\node[anchor=center,text=qrgC7c3aed,font={\fontsize{9}{10}\selectfont},inner sep=0pt] at (957.000,373.710) {8};
% Q01 -> Q30: 22 maps
\draw[->,>=stealth,draw=qrgC7c3aed,line width=.55pt] (1026.870,614.890) .. controls (1210.000,515.000) and (1430.000,395.000) .. (1628.050,286.970);
\node[anchor=center,text=qrgC7c3aed,font={\fontsize{9}{10}\selectfont},inner sep=0pt] at (1363.510,383.430) {22};
% Q01 -> Q31: 22 maps
\draw[->,>=stealth,draw=qrgC7c3aed,line width=.55pt] (1029.440,620.570) .. controls (1258.440,536.770) and (1516.780,416.770) .. (1743.290,287.400);
\node[anchor=center,text=qrgC7c3aed,font={\fontsize{9}{10}\selectfont},inner sep=0pt] at (1478.460,474.210) {22};
% Q01 -> Q20: 16 maps
\draw[->,>=stealth,draw=qrgC7c3aed,line width=.55pt] (1032.000,634.930) .. controls (1231.360,634.570) and (1394.400,886.270) .. (1625.400,844.450);
\node[anchor=center,text=qrgC7c3aed,font={\fontsize{9}{10}\selectfont},inner sep=0pt] at (1354.000,825.510) {16};
% Q01 -> Q37: 30 maps
\draw[->,>=stealth,draw=qrgC7c3aed,line width=.55pt] (1004.400,674.450) .. controls (1028.330,740.000) and (1066.670,845.000) .. (1096.430,926.520);
\node[anchor=center,text=qrgC7c3aed,font={\fontsize{9}{10}\selectfont},inner sep=0pt] at (1028.500,862.910) {30};
% Q01 -> Q36: 30 maps
\draw[->,>=stealth,draw=qrgC7c3aed,line width=.55pt] (990.000,677.000) .. controls (990.000,740.000) and (990.000,845.000) .. (990.000,925.000);
\node[anchor=center,text=qrgC7c3aed,font={\fontsize{9}{10}\selectfont},inner sep=0pt] at (957.000,768.470) {30};
% Q01 -> Q35: 24 maps
\draw[->,>=stealth,draw=qrgC7c3aed,line width=.55pt] (950.690,649.790) .. controls (808.330,703.330) and (626.670,771.670) .. (468.400,831.200);
\node[anchor=center,text=qrgC7c3aed,font={\fontsize{9}{10}\selectfont},inner sep=0pt] at (680.330,796.360) {24};
% Q01 -> Q34: 24 maps
\draw[->,>=stealth,draw=qrgC7c3aed,line width=.55pt] (948.820,643.240) .. controls (762.880,680.410) and (542.880,748.750) .. (352.980,830.150);
\node[anchor=center,text=qrgC7c3aed,font={\fontsize{9}{10}\selectfont},inner sep=0pt] at (639.840,680.020) {24};
% Q01 -> Q01: 2 maps
\draw[->,>=stealth,draw=qrgC64748b,line width=.4pt] (1007.327,672.159) .. controls (1044.428,705.463) and (1074.822,662.047) .. (1030.844,638.574);
\node[anchor=center,text=qrgC64748b,font={\fontsize{9}{10}\selectfont},inner sep=0pt] at (1088.300,703.830) {2};
% Q04 -> Q04: 8 maps
\draw[->,>=stealth,draw=qrgC64748b,line width=.4pt] (998.208,217.553) .. controls (1020.557,244.138) and (1043.703,216.562) .. (1013.635,199.168);
\node[anchor=center,text=qrgC64748b,font={\fontsize{9}{10}\selectfont},inner sep=0pt] at (1058.940,252.850) {8};
% Q06 -> Q06: 24 maps
\draw[->,>=stealth,draw=qrgC64748b,line width=.4pt] (307.447,283.208) .. controls (280.862,305.557) and (308.438,328.703) .. (325.832,298.635);
\node[anchor=center,text=qrgC64748b,font={\fontsize{9}{10}\selectfont},inner sep=0pt] at (340.020,340.850) {24};
% Q07 -> Q07: 24 maps
\draw[->,>=stealth,draw=qrgC64748b,line width=.4pt] (433.001,254.215) .. controls (406.372,231.909) and (388.368,263.091) .. (421.000,274.999);
\node[anchor=center,text=qrgC64748b,font={\fontsize{9}{10}\selectfont},inner sep=0pt] at (441.300,220.900) {24};
% Q20 -> Q20: 48 maps
\draw[->,>=stealth,draw=qrgC64748b,line width=.4pt] (1669.661,853.763) .. controls (1703.224,862.729) and (1706.356,826.861) .. (1671.750,829.854);
\node[anchor=center,text=qrgC64748b,font={\fontsize{9}{10}\selectfont},inner sep=0pt] at (1739.660,847.840) {48};
% Q30 -> Q30: 16 maps
\draw[->,>=stealth,draw=qrgC64748b,line width=.4pt] (1672.553,266.792) .. controls (1699.138,244.443) and (1671.562,221.297) .. (1654.168,251.365);
\node[anchor=center,text=qrgC64748b,font={\fontsize{9}{10}\selectfont},inner sep=0pt] at (1639.980,209.150) {16};
% Q31 -> Q31: 16 maps
\draw[->,>=stealth,draw=qrgC64748b,line width=.4pt] (1783.386,290.426) .. controls (1816.041,302.278) and (1822.289,266.822) .. (1787.552,266.790);
\node[anchor=center,text=qrgC64748b,font={\fontsize{9}{10}\selectfont},inner sep=0pt] at (1799.650,331.880) {16};
% Q34 -> Q34: 48 maps
\draw[->,>=stealth,draw=qrgC64748b,line width=.4pt] (307.447,848.208) .. controls (280.862,870.557) and (308.438,893.703) .. (325.832,863.635);
\node[anchor=center,text=qrgC64748b,font={\fontsize{9}{10}\selectfont},inner sep=0pt] at (272.150,908.940) {48};
% Q35 -> Q35: 48 maps
\draw[->,>=stealth,draw=qrgC64748b,line width=.4pt] (421.000,840.001) .. controls (388.368,851.909) and (406.372,883.091) .. (433.001,860.785);
\node[anchor=center,text=qrgC64748b,font={\fontsize{9}{10}\selectfont},inner sep=0pt] at (441.300,894.100) {48};
% Q36 -> Q36: 720 maps
\draw[->,>=stealth,draw=qrgC64748b,line width=.4pt] (966.091,947.911) .. controls (932.548,956.930) and (947.762,989.560) .. (976.236,969.661);
\node[anchor=center,text=qrgC64748b,font={\fontsize{9}{10}\selectfont},inner sep=0pt] at (908.430,988.040) {720};
% Q37 -> Q37: 720 maps
\draw[->,>=stealth,draw=qrgC64748b,line width=.4pt] (1107.093,973.909) .. controls (1121.807,1005.379) and (1151.293,984.731) .. (1126.752,960.142);
\node[anchor=center,text=qrgC64748b,font={\fontsize{9}{10}\selectfont},inner sep=0pt] at (1171.190,957.360) {720};
\node[anchor=west,text=qrgC7c3aed,font={\fontsize{9}{10}\selectfont},inner sep=0pt] at (1008.000,525.000) {\(K_2\)};
\draw[draw=qrgC94a3b8,line width=.45pt] (956.000,491.000) -- (956.000,559.000);
\draw[draw=none,fill=qrgC334155,line width=.7pt] (956.000,491.000) circle[radius=0.0320cm];
\draw[draw=none,fill=qrgC334155,line width=.7pt] (956.000,559.000) circle[radius=0.0320cm];
\draw[draw=qrgC16a34a,fill=white,line width=.7pt] (990.000,635.000) circle[radius=0.3200cm];
\node[anchor=center,text=qrgC334155,font={\fontsize{8}{9}\selectfont},inner sep=0pt] at (990.000,635.000) {0};
\node[anchor=west,text=qrgC7c3aed,font={\fontsize{9}{10}\selectfont},inner sep=0pt] at (1008.000,85.000) {\(C_4\)};
\draw[draw=qrgC94a3b8,line width=.45pt] (956.000,51.000) -- (990.000,85.000);
\draw[draw=qrgC94a3b8,line width=.45pt] (990.000,85.000) -- (956.000,119.000);
\draw[draw=qrgC94a3b8,line width=.45pt] (956.000,119.000) -- (922.000,85.000);
\draw[draw=qrgC94a3b8,line width=.45pt] (922.000,85.000) -- (956.000,51.000);
\draw[draw=none,fill=qrgC334155,line width=.7pt] (956.000,51.000) circle[radius=0.0320cm];
\draw[draw=none,fill=qrgC334155,line width=.7pt] (990.000,85.000) circle[radius=0.0320cm];
\draw[draw=none,fill=qrgC334155,line width=.7pt] (956.000,119.000) circle[radius=0.0320cm];
\draw[draw=none,fill=qrgC334155,line width=.7pt] (922.000,85.000) circle[radius=0.0320cm];
\draw[draw=qrgC16a34a,fill=white,line width=.7pt] (990.000,195.000) circle[radius=0.1840cm];
\node[anchor=center,text=qrgC334155,font={\fontsize{8}{9}\selectfont},inner sep=0pt] at (990.000,195.000) {0};
\node[anchor=west,text=qrgC7c3aed,font={\fontsize{9}{10}\selectfont},inner sep=0pt] at (405.500,165.000) {\(K_4\)};
\draw[draw=qrgC94a3b8,line width=.45pt] (353.500,131.000) -- (387.500,165.000);
\draw[draw=qrgC94a3b8,line width=.45pt] (353.500,131.000) -- (353.500,199.000);
\draw[draw=qrgC94a3b8,line width=.45pt] (353.500,131.000) -- (319.500,165.000);
\draw[draw=qrgC94a3b8,line width=.45pt] (387.500,165.000) -- (353.500,199.000);
\draw[draw=qrgC94a3b8,line width=.45pt] (387.500,165.000) -- (319.500,165.000);
\draw[draw=qrgC94a3b8,line width=.45pt] (353.500,199.000) -- (319.500,165.000);
\draw[draw=none,fill=qrgC334155,line width=.7pt] (353.500,131.000) circle[radius=0.0320cm];
\draw[draw=none,fill=qrgC334155,line width=.7pt] (387.500,165.000) circle[radius=0.0320cm];
\draw[draw=none,fill=qrgC334155,line width=.7pt] (353.500,199.000) circle[radius=0.0320cm];
\draw[draw=none,fill=qrgC334155,line width=.7pt] (319.500,165.000) circle[radius=0.0320cm];
\draw[draw=qrgC16a34a,fill=white,line width=.7pt] (330.000,275.000) circle[radius=0.1840cm];
\node[anchor=center,text=qrgC334155,font={\fontsize{8}{9}\selectfont},inner sep=0pt] at (330.000,275.000) {0};
\draw[draw=qrgC64748b,fill=white,line width=.7pt] (445.000,275.000) circle[radius=0.1840cm];
\node[anchor=center,text=qrgC334155,font={\fontsize{8}{9}\selectfont},inner sep=0pt] at (445.000,275.000) {1};
\node[anchor=west,text=qrgC7c3aed,font={\fontsize{9}{10}\selectfont},inner sep=0pt] at (1650.500,970.000) {\(K_{2,4}\)};
\draw[draw=qrgC94a3b8,line width=.45pt] (1598.500,936.000) -- (1569.060,987.000);
\draw[draw=qrgC94a3b8,line width=.45pt] (1598.500,936.000) -- (1569.060,953.000);
\draw[draw=qrgC94a3b8,line width=.45pt] (1627.940,953.000) -- (1569.060,987.000);
\draw[draw=qrgC94a3b8,line width=.45pt] (1627.940,953.000) -- (1569.060,953.000);
\draw[draw=qrgC94a3b8,line width=.45pt] (1627.940,987.000) -- (1569.060,987.000);
\draw[draw=qrgC94a3b8,line width=.45pt] (1627.940,987.000) -- (1569.060,953.000);
\draw[draw=qrgC94a3b8,line width=.45pt] (1598.500,1004.000) -- (1569.060,987.000);
\draw[draw=qrgC94a3b8,line width=.45pt] (1598.500,1004.000) -- (1569.060,953.000);
\draw[draw=none,fill=qrgC334155,line width=.7pt] (1598.500,936.000) circle[radius=0.0320cm];
\draw[draw=none,fill=qrgC334155,line width=.7pt] (1627.940,953.000) circle[radius=0.0320cm];
\draw[draw=none,fill=qrgC334155,line width=.7pt] (1627.940,987.000) circle[radius=0.0320cm];
\draw[draw=none,fill=qrgC334155,line width=.7pt] (1598.500,1004.000) circle[radius=0.0320cm];
\draw[draw=none,fill=qrgC334155,line width=.7pt] (1569.060,987.000) circle[radius=0.0320cm];
\draw[draw=none,fill=qrgC334155,line width=.7pt] (1569.060,953.000) circle[radius=0.0320cm];
\draw[draw=qrgC64748b,fill=white,line width=.7pt] (1650.000,840.000) circle[radius=0.1840cm];
\node[anchor=center,text=qrgC334155,font={\fontsize{8}{9}\selectfont},inner sep=0pt] at (1650.000,840.000) {0};
\node[anchor=west,text=qrgC7c3aed,font={\fontsize{9}{10}\selectfont},inner sep=0pt] at (1675.500,145.000) {\(K_2\vee2K_2\)};
\draw[draw=qrgC94a3b8,line width=.45pt] (1623.500,111.000) -- (1623.500,179.000);
\draw[draw=qrgC94a3b8,line width=.45pt] (1623.500,111.000) -- (1594.060,162.000);
\draw[draw=qrgC94a3b8,line width=.45pt] (1623.500,111.000) -- (1594.060,128.000);
\draw[draw=qrgC94a3b8,line width=.45pt] (1652.940,128.000) -- (1652.940,162.000);
\draw[draw=qrgC94a3b8,line width=.45pt] (1652.940,128.000) -- (1594.060,162.000);
\draw[draw=qrgC94a3b8,line width=.45pt] (1652.940,128.000) -- (1594.060,128.000);
\draw[draw=qrgC94a3b8,line width=.45pt] (1652.940,162.000) -- (1594.060,162.000);
\draw[draw=qrgC94a3b8,line width=.45pt] (1652.940,162.000) -- (1594.060,128.000);
\draw[draw=qrgC94a3b8,line width=.45pt] (1623.500,179.000) -- (1594.060,162.000);
\draw[draw=qrgC94a3b8,line width=.45pt] (1623.500,179.000) -- (1594.060,128.000);
\draw[draw=qrgC94a3b8,line width=.45pt] (1594.060,162.000) -- (1594.060,128.000);
\draw[draw=none,fill=qrgC334155,line width=.7pt] (1623.500,111.000) circle[radius=0.0320cm];
\draw[draw=none,fill=qrgC334155,line width=.7pt] (1652.940,128.000) circle[radius=0.0320cm];
\draw[draw=none,fill=qrgC334155,line width=.7pt] (1652.940,162.000) circle[radius=0.0320cm];
\draw[draw=none,fill=qrgC334155,line width=.7pt] (1623.500,179.000) circle[radius=0.0320cm];
\draw[draw=none,fill=qrgC334155,line width=.7pt] (1594.060,162.000) circle[radius=0.0320cm];
\draw[draw=none,fill=qrgC334155,line width=.7pt] (1594.060,128.000) circle[radius=0.0320cm];
\draw[draw=qrgC64748b,fill=white,line width=.7pt] (1650.000,275.000) circle[radius=0.1840cm];
\node[anchor=center,text=qrgC334155,font={\fontsize{8}{9}\selectfont},inner sep=0pt] at (1650.000,275.000) {0};
\draw[draw=qrgC64748b,fill=white,line width=.7pt] (1765.000,275.000) circle[radius=0.1840cm];
\node[anchor=center,text=qrgC334155,font={\fontsize{8}{9}\selectfont},inner sep=0pt] at (1765.000,275.000) {1};
\node[anchor=west,text=qrgC7c3aed,font={\fontsize{9}{10}\selectfont},inner sep=0pt] at (375.500,970.000) {\(K_{2,2,2}\)};
\draw[draw=qrgC94a3b8,line width=.45pt] (323.500,936.000) -- (352.940,987.000);
\draw[draw=qrgC94a3b8,line width=.45pt] (323.500,936.000) -- (323.500,1004.000);
\draw[draw=qrgC94a3b8,line width=.45pt] (323.500,936.000) -- (294.060,987.000);
\draw[draw=qrgC94a3b8,line width=.45pt] (323.500,936.000) -- (294.060,953.000);
\draw[draw=qrgC94a3b8,line width=.45pt] (352.940,953.000) -- (352.940,987.000);
\draw[draw=qrgC94a3b8,line width=.45pt] (352.940,953.000) -- (323.500,1004.000);
\draw[draw=qrgC94a3b8,line width=.45pt] (352.940,953.000) -- (294.060,987.000);
\draw[draw=qrgC94a3b8,line width=.45pt] (352.940,953.000) -- (294.060,953.000);
\draw[draw=qrgC94a3b8,line width=.45pt] (352.940,987.000) -- (294.060,987.000);
\draw[draw=qrgC94a3b8,line width=.45pt] (352.940,987.000) -- (294.060,953.000);
\draw[draw=qrgC94a3b8,line width=.45pt] (323.500,1004.000) -- (294.060,987.000);
\draw[draw=qrgC94a3b8,line width=.45pt] (323.500,1004.000) -- (294.060,953.000);
\draw[draw=none,fill=qrgC334155,line width=.7pt] (323.500,936.000) circle[radius=0.0320cm];
\draw[draw=none,fill=qrgC334155,line width=.7pt] (352.940,953.000) circle[radius=0.0320cm];
\draw[draw=none,fill=qrgC334155,line width=.7pt] (352.940,987.000) circle[radius=0.0320cm];
\draw[draw=none,fill=qrgC334155,line width=.7pt] (323.500,1004.000) circle[radius=0.0320cm];
\draw[draw=none,fill=qrgC334155,line width=.7pt] (294.060,987.000) circle[radius=0.0320cm];
\draw[draw=none,fill=qrgC334155,line width=.7pt] (294.060,953.000) circle[radius=0.0320cm];
\draw[draw=qrgC64748b,fill=white,line width=.7pt] (330.000,840.000) circle[radius=0.1840cm];
\node[anchor=center,text=qrgC334155,font={\fontsize{8}{9}\selectfont},inner sep=0pt] at (330.000,840.000) {1};
\draw[draw=qrgC64748b,fill=white,line width=.7pt] (445.000,840.000) circle[radius=0.1840cm];
\node[anchor=center,text=qrgC334155,font={\fontsize{8}{9}\selectfont},inner sep=0pt] at (445.000,840.000) {0};
\node[anchor=west,text=qrgC7c3aed,font={\fontsize{9}{10}\selectfont},inner sep=0pt] at (1065.500,1050.000) {\(K_6\)};
\draw[draw=qrgC94a3b8,line width=.45pt] (1013.500,1016.000) -- (1042.940,1033.000);
\draw[draw=qrgC94a3b8,line width=.45pt] (1013.500,1016.000) -- (1042.940,1067.000);
\draw[draw=qrgC94a3b8,line width=.45pt] (1013.500,1016.000) -- (1013.500,1084.000);
\draw[draw=qrgC94a3b8,line width=.45pt] (1013.500,1016.000) -- (984.060,1067.000);
\draw[draw=qrgC94a3b8,line width=.45pt] (1013.500,1016.000) -- (984.060,1033.000);
\draw[draw=qrgC94a3b8,line width=.45pt] (1042.940,1033.000) -- (1042.940,1067.000);
\draw[draw=qrgC94a3b8,line width=.45pt] (1042.940,1033.000) -- (1013.500,1084.000);
\draw[draw=qrgC94a3b8,line width=.45pt] (1042.940,1033.000) -- (984.060,1067.000);
\draw[draw=qrgC94a3b8,line width=.45pt] (1042.940,1033.000) -- (984.060,1033.000);
\draw[draw=qrgC94a3b8,line width=.45pt] (1042.940,1067.000) -- (1013.500,1084.000);
\draw[draw=qrgC94a3b8,line width=.45pt] (1042.940,1067.000) -- (984.060,1067.000);
\draw[draw=qrgC94a3b8,line width=.45pt] (1042.940,1067.000) -- (984.060,1033.000);
\draw[draw=qrgC94a3b8,line width=.45pt] (1013.500,1084.000) -- (984.060,1067.000);
\draw[draw=qrgC94a3b8,line width=.45pt] (1013.500,1084.000) -- (984.060,1033.000);
\draw[draw=qrgC94a3b8,line width=.45pt] (984.060,1067.000) -- (984.060,1033.000);
\draw[draw=none,fill=qrgC334155,line width=.7pt] (1013.500,1016.000) circle[radius=0.0320cm];
\draw[draw=none,fill=qrgC334155,line width=.7pt] (1042.940,1033.000) circle[radius=0.0320cm];
\draw[draw=none,fill=qrgC334155,line width=.7pt] (1042.940,1067.000) circle[radius=0.0320cm];
\draw[draw=none,fill=qrgC334155,line width=.7pt] (1013.500,1084.000) circle[radius=0.0320cm];
\draw[draw=none,fill=qrgC334155,line width=.7pt] (984.060,1067.000) circle[radius=0.0320cm];
\draw[draw=none,fill=qrgC334155,line width=.7pt] (984.060,1033.000) circle[radius=0.0320cm];
\draw[draw=qrgC16a34a,fill=white,line width=.7pt] (990.000,950.000) circle[radius=0.1840cm];
\node[anchor=center,text=qrgC334155,font={\fontsize{8}{9}\selectfont},inner sep=0pt] at (990.000,950.000) {0};
\draw[draw=qrgC64748b,fill=white,line width=.7pt] (1105.000,950.000) circle[radius=0.1840cm];
\node[anchor=center,text=qrgC334155,font={\fontsize{8}{9}\selectfont},inner sep=0pt] at (1105.000,950.000) {1};
\end{tikzpicture}
% END GENERATED QRG CO-EULERIAN
\caption{All co-Eulerian QRGs maps between connected graphs with 2-6 vertices.
Circle labels give \(\alpha\)-parameters \(0,1\), with \(0\) inner.
Arrows are labelled by map counts; grey denotes isomorphisms,
green outlines flatness.}
\label{fig:QRG-quiver-co-eulerian}
\label{fig:QRG-quiver-six}
\end{figure}
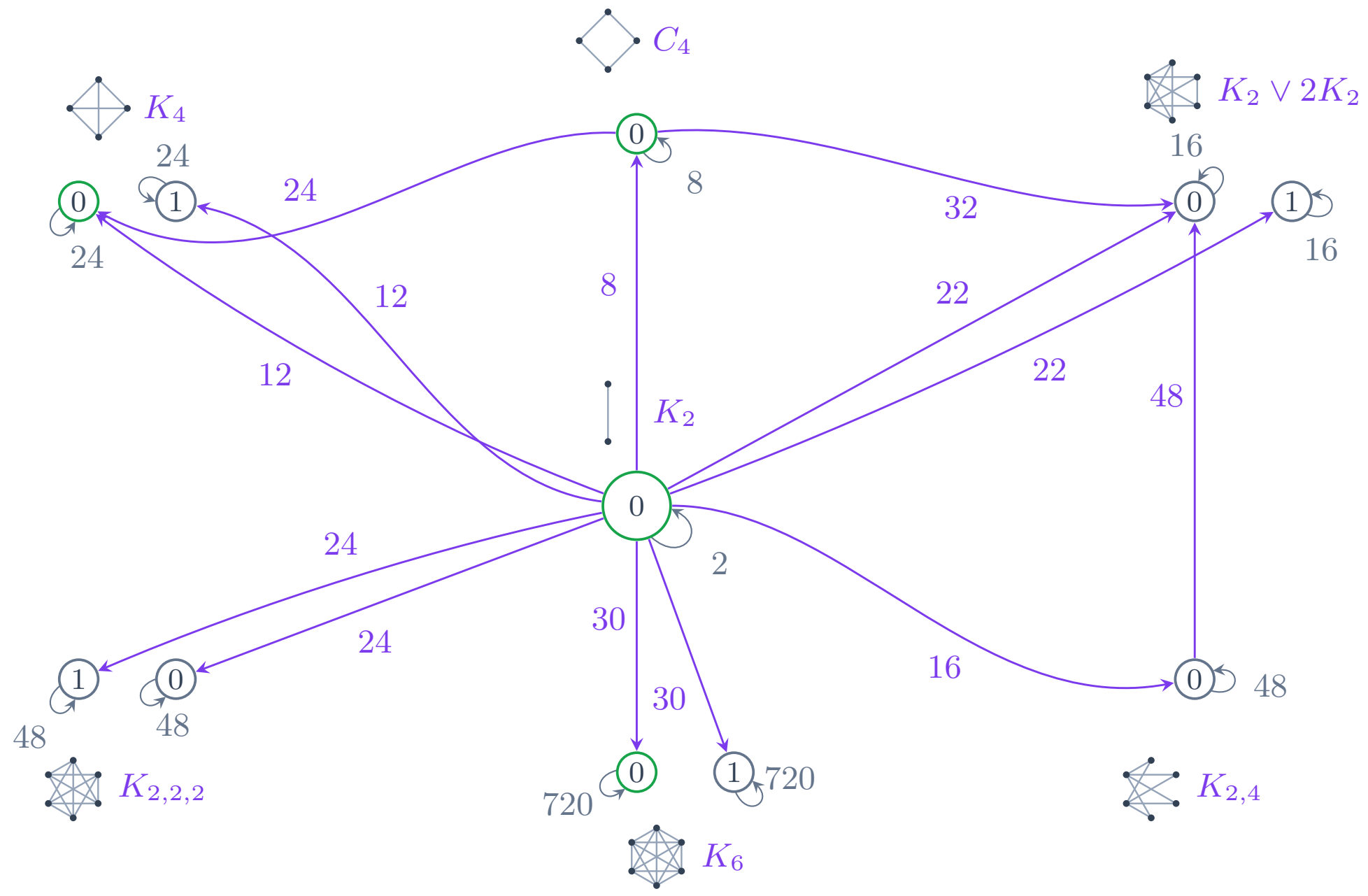

\begin{figure}[tbp]
\centering
% BEGIN GENERATED QRG EULERIAN
\definecolor{qrgC0369a1}{HTML}{0369A1}
\definecolor{qrgC16a34a}{HTML}{16A34A}
\definecolor{qrgC334155}{HTML}{334155}
\definecolor{qrgC64748b}{HTML}{64748B}
\definecolor{qrgC94a3b8}{HTML}{94A3B8}
\begin{tikzpicture}[x=.008cm,y=-.008cm]
\path[use as bounding box] (0,0) rectangle (2000,1150);
% Q17 -> Q02: 6 maps
\draw[->,>=stealth,draw=qrgC0369a1,line width=.55pt] (355.000,180.000) .. controls (623.330,180.000) and (916.670,180.000) .. (1185.000,180.000);
\node[anchor=center,text=qrgC0369a1,font={\fontsize{9}{10}\selectfont},inner sep=0pt] at (854.200,147.000) {6};
% Q02 -> Q02: 6 maps
\draw[->,>=stealth,draw=qrgC64748b,line width=.4pt] (1234.000,179.999) .. controls (1266.632,168.091) and (1248.628,136.909) .. (1221.999,159.215);
\node[anchor=center,text=qrgC64748b,font={\fontsize{9}{10}\selectfont},inner sep=0pt] at (1231.131,130.636) {6};
% Q03 -> Q03: 6 maps
\draw[->,>=stealth,draw=qrgC64748b,line width=.4pt] (1345.785,192.000) .. controls (1380.000,198.000) and (1380.000,162.000) .. (1345.785,168.000);
\node[anchor=center,text=qrgC64748b,font={\fontsize{9}{10}\selectfont},inner sep=0pt] at (1396.446,180.000) {6};
% Q05 -> Q05: 8 maps
\draw[->,>=stealth,draw=qrgC64748b,line width=.4pt] (290.785,462.000) .. controls (325.000,468.000) and (325.000,432.000) .. (290.785,438.000);
\node[anchor=center,text=qrgC64748b,font={\fontsize{9}{10}\selectfont},inner sep=0pt] at (341.446,450.000) {8};
% Q08 -> Q08: 10 maps
\draw[->,>=stealth,draw=qrgC64748b,line width=.4pt] (690.785,462.000) .. controls (725.000,468.000) and (725.000,432.000) .. (690.785,438.000);
\node[anchor=center,text=qrgC64748b,font={\fontsize{9}{10}\selectfont},inner sep=0pt] at (760.000,450.000) {10};
% Q09 -> Q09: 8 maps
\draw[->,>=stealth,draw=qrgC64748b,line width=.4pt] (1043.764,469.661) .. controls (1072.238,489.560) and (1087.452,456.930) .. (1053.909,447.911);
\node[anchor=center,text=qrgC64748b,font={\fontsize{9}{10}\selectfont},inner sep=0pt] at (1060.434,491.745) {8};
% Q10 -> Q10: 8 maps
\draw[->,>=stealth,draw=qrgC64748b,line width=.4pt] (1136.752,439.858) .. controls (1161.293,415.269) and (1131.807,394.621) .. (1117.093,426.091);
\node[anchor=center,text=qrgC64748b,font={\fontsize{9}{10}\selectfont},inner sep=0pt] at (1155.982,391.473) {8};
% Q11 -> Q11: 4 maps
\draw[->,>=stealth,draw=qrgC64748b,line width=.4pt] (1183.027,466.968) .. controls (1164.494,512.959) and (1199.266,522.271) .. (1206.208,473.183);
\node[anchor=center,text=qrgC64748b,font={\fontsize{9}{10}\selectfont},inner sep=0pt] at (1178.593,529.878) {4};
% Q12 -> Q12: 4 maps
\draw[->,>=stealth,draw=qrgC64748b,line width=.4pt] (1381.212,473.182) .. controls (1401.165,501.615) and (1426.615,476.165) .. (1398.182,456.212);
\node[anchor=center,text=qrgC64748b,font={\fontsize{9}{10}\selectfont},inner sep=0pt] at (1425.519,500.519) {4};
% Q13 -> Q13: 12 maps
\draw[->,>=stealth,draw=qrgC64748b,line width=.4pt] (1648.386,465.426) .. controls (1681.041,477.278) and (1687.289,441.822) .. (1652.552,441.790);
\node[anchor=center,text=qrgC64748b,font={\fontsize{9}{10}\selectfont},inner sep=0pt] at (1664.650,506.880) {12};
% Q14 -> Q14: 12 maps
\draw[->,>=stealth,draw=qrgC64748b,line width=.4pt] (1735.785,462.000) .. controls (1770.000,468.000) and (1770.000,432.000) .. (1735.785,438.000);
\node[anchor=center,text=qrgC64748b,font={\fontsize{9}{10}\selectfont},inner sep=0pt] at (1805.000,450.000) {12};
% Q15 -> Q15: 120 maps
\draw[->,>=stealth,draw=qrgC64748b,line width=.4pt] (288.386,745.426) .. controls (321.041,757.278) and (327.289,721.822) .. (292.552,721.790);
\node[anchor=center,text=qrgC64748b,font={\fontsize{9}{10}\selectfont},inner sep=0pt] at (304.650,786.880) {120};
% Q16 -> Q16: 120 maps
\draw[->,>=stealth,draw=qrgC64748b,line width=.4pt] (375.785,742.000) .. controls (410.000,748.000) and (410.000,712.000) .. (375.785,718.000);
\node[anchor=center,text=qrgC64748b,font={\fontsize{9}{10}\selectfont},inner sep=0pt] at (445.000,730.000) {120};
% Q17 -> Q17: 12 maps
\draw[->,>=stealth,draw=qrgC64748b,line width=.4pt] (316.236,160.339) .. controls (287.762,140.440) and (272.548,173.070) .. (306.091,182.089);
\node[anchor=center,text=qrgC64748b,font={\fontsize{9}{10}\selectfont},inner sep=0pt] at (278.000,196.000) {12};
% Q18 -> Q18: 4 maps
\draw[->,>=stealth,draw=qrgC64748b,line width=.4pt] (688.386,745.426) .. controls (721.041,757.278) and (727.289,721.822) .. (692.552,721.790);
\node[anchor=center,text=qrgC64748b,font={\fontsize{9}{10}\selectfont},inner sep=0pt] at (710.202,762.443) {4};
% Q19 -> Q19: 4 maps
\draw[->,>=stealth,draw=qrgC64748b,line width=.4pt] (775.785,742.000) .. controls (810.000,748.000) and (810.000,712.000) .. (775.785,718.000);
\node[anchor=center,text=qrgC64748b,font={\fontsize{9}{10}\selectfont},inner sep=0pt] at (826.446,730.000) {4};
% Q21 -> Q21: 48 maps
\draw[->,>=stealth,draw=qrgC64748b,line width=.4pt] (1090.785,742.000) .. controls (1125.000,748.000) and (1125.000,712.000) .. (1090.785,718.000);
\node[anchor=center,text=qrgC64748b,font={\fontsize{9}{10}\selectfont},inner sep=0pt] at (1160.000,730.000) {48};
% Q22 -> Q22: 4 maps
\draw[->,>=stealth,draw=qrgC64748b,line width=.4pt] (1488.386,745.426) .. controls (1521.041,757.278) and (1527.289,721.822) .. (1492.552,721.790);
\node[anchor=center,text=qrgC64748b,font={\fontsize{9}{10}\selectfont},inner sep=0pt] at (1510.202,762.443) {4};
% Q23 -> Q23: 4 maps
\draw[->,>=stealth,draw=qrgC64748b,line width=.4pt] (1575.785,742.000) .. controls (1610.000,748.000) and (1610.000,712.000) .. (1575.785,718.000);
\node[anchor=center,text=qrgC64748b,font={\fontsize{9}{10}\selectfont},inner sep=0pt] at (1626.446,730.000) {4};
% Q24 -> Q24: 6 maps
\draw[->,>=stealth,draw=qrgC64748b,line width=.4pt] (288.386,1025.426) .. controls (321.041,1037.278) and (327.289,1001.822) .. (292.552,1001.790);
\node[anchor=center,text=qrgC64748b,font={\fontsize{9}{10}\selectfont},inner sep=0pt] at (310.202,1042.443) {6};
% Q25 -> Q25: 6 maps
\draw[->,>=stealth,draw=qrgC64748b,line width=.4pt] (375.785,1022.000) .. controls (410.000,1028.000) and (410.000,992.000) .. (375.785,998.000);
\node[anchor=center,text=qrgC64748b,font={\fontsize{9}{10}\selectfont},inner sep=0pt] at (426.446,1010.000) {6};
% Q26 -> Q26: 8 maps
\draw[->,>=stealth,draw=qrgC64748b,line width=.4pt] (688.386,1025.426) .. controls (721.041,1037.278) and (727.289,1001.822) .. (692.552,1001.790);
\node[anchor=center,text=qrgC64748b,font={\fontsize{9}{10}\selectfont},inner sep=0pt] at (710.202,1042.443) {8};
% Q27 -> Q27: 8 maps
\draw[->,>=stealth,draw=qrgC64748b,line width=.4pt] (775.785,1022.000) .. controls (810.000,1028.000) and (810.000,992.000) .. (775.785,998.000);
\node[anchor=center,text=qrgC64748b,font={\fontsize{9}{10}\selectfont},inner sep=0pt] at (826.446,1010.000) {8};
% Q28 -> Q28: 12 maps
\draw[->,>=stealth,draw=qrgC64748b,line width=.4pt] (1088.386,1025.426) .. controls (1121.041,1037.278) and (1127.289,1001.822) .. (1092.552,1001.790);
\node[anchor=center,text=qrgC64748b,font={\fontsize{9}{10}\selectfont},inner sep=0pt] at (1104.650,1066.880) {12};
% Q29 -> Q29: 12 maps
\draw[->,>=stealth,draw=qrgC64748b,line width=.4pt] (1175.785,1022.000) .. controls (1210.000,1028.000) and (1210.000,992.000) .. (1175.785,998.000);
\node[anchor=center,text=qrgC64748b,font={\fontsize{9}{10}\selectfont},inner sep=0pt] at (1245.000,1010.000) {12};
% Q32 -> Q32: 48 maps
\draw[->,>=stealth,draw=qrgC64748b,line width=.4pt] (1488.386,1025.426) .. controls (1521.041,1037.278) and (1527.289,1001.822) .. (1492.552,1001.790);
\node[anchor=center,text=qrgC64748b,font={\fontsize{9}{10}\selectfont},inner sep=0pt] at (1504.650,1066.880) {48};
% Q33 -> Q33: 48 maps
\draw[->,>=stealth,draw=qrgC64748b,line width=.4pt] (1575.785,1022.000) .. controls (1610.000,1028.000) and (1610.000,992.000) .. (1575.785,998.000);
\node[anchor=center,text=qrgC64748b,font={\fontsize{9}{10}\selectfont},inner sep=0pt] at (1645.000,1010.000) {48};
% Q11 -> Q12: 4 maps
\draw[->,>=stealth,draw=qrgC64748b,line width=.4pt] (1221.000,460.000) .. controls (1258.000,473.000) and (1317.000,473.000) .. (1354.000,460.000);
\node[anchor=center,text=qrgC64748b,font={\fontsize{9}{10}\selectfont},inner sep=0pt] at (1287.500,494.000) {4};
% Q12 -> Q11: 4 maps
\draw[->,>=stealth,draw=qrgC64748b,line width=.4pt] (1354.000,440.000) .. controls (1317.000,427.000) and (1258.000,427.000) .. (1221.000,440.000);
\node[anchor=center,text=qrgC64748b,font={\fontsize{9}{10}\selectfont},inner sep=0pt] at (1287.500,406.000) {4};
\node[anchor=west,text=qrgC0369a1,font={\fontsize{9}{10}\selectfont},inner sep=0pt] at (1285.500,70.000) {\(C_3\)};
\draw[draw=qrgC94a3b8,line width=.45pt] (1233.500,36.000) -- (1262.940,87.000);
\draw[draw=qrgC94a3b8,line width=.45pt] (1233.500,36.000) -- (1204.060,87.000);
\draw[draw=qrgC94a3b8,line width=.45pt] (1262.940,87.000) -- (1204.060,87.000);
\draw[draw=none,fill=qrgC334155,line width=.7pt] (1233.500,36.000) circle[radius=0.0320cm];
\draw[draw=none,fill=qrgC334155,line width=.7pt] (1262.940,87.000) circle[radius=0.0320cm];
\draw[draw=none,fill=qrgC334155,line width=.7pt] (1204.060,87.000) circle[radius=0.0320cm];
\draw[draw=qrgC16a34a,fill=white,line width=.7pt] (1210.000,180.000) circle[radius=0.1840cm];
\node[anchor=center,text=qrgC334155,font={\fontsize{8}{9}\selectfont},inner sep=0pt] at (1210.000,180.000) {0};
\draw[draw=qrgC64748b,fill=white,line width=.7pt] (1325.000,180.000) circle[radius=0.1840cm];
\node[anchor=center,text=qrgC334155,font={\fontsize{8}{9}\selectfont},inner sep=0pt] at (1325.000,180.000) {1};
\node[anchor=west,text=qrgC0369a1,font={\fontsize{9}{10}\selectfont},inner sep=0pt] at (288.000,340.000) {\(C_4\)};
\draw[draw=qrgC94a3b8,line width=.45pt] (236.000,306.000) -- (270.000,340.000);
\draw[draw=qrgC94a3b8,line width=.45pt] (270.000,340.000) -- (236.000,374.000);
\draw[draw=qrgC94a3b8,line width=.45pt] (236.000,374.000) -- (202.000,340.000);
\draw[draw=qrgC94a3b8,line width=.45pt] (202.000,340.000) -- (236.000,306.000);
\draw[draw=none,fill=qrgC334155,line width=.7pt] (236.000,306.000) circle[radius=0.0320cm];
\draw[draw=none,fill=qrgC334155,line width=.7pt] (270.000,340.000) circle[radius=0.0320cm];
\draw[draw=none,fill=qrgC334155,line width=.7pt] (236.000,374.000) circle[radius=0.0320cm];
\draw[draw=none,fill=qrgC334155,line width=.7pt] (202.000,340.000) circle[radius=0.0320cm];
\draw[draw=qrgC16a34a,fill=white,line width=.7pt] (270.000,450.000) circle[radius=0.1840cm];
\node[anchor=center,text=qrgC334155,font={\fontsize{8}{9}\selectfont},inner sep=0pt] at (270.000,450.000) {0};
\node[anchor=west,text=qrgC0369a1,font={\fontsize{9}{10}\selectfont},inner sep=0pt] at (688.000,340.000) {\(C_5\)};
\draw[draw=qrgC94a3b8,line width=.45pt] (636.000,306.000) -- (668.340,329.490);
\draw[draw=qrgC94a3b8,line width=.45pt] (668.340,329.490) -- (655.980,367.510);
\draw[draw=qrgC94a3b8,line width=.45pt] (655.980,367.510) -- (616.020,367.510);
\draw[draw=qrgC94a3b8,line width=.45pt] (616.020,367.510) -- (603.660,329.490);
\draw[draw=qrgC94a3b8,line width=.45pt] (603.660,329.490) -- (636.000,306.000);
\draw[draw=none,fill=qrgC334155,line width=.7pt] (636.000,306.000) circle[radius=0.0320cm];
\draw[draw=none,fill=qrgC334155,line width=.7pt] (668.340,329.490) circle[radius=0.0320cm];
\draw[draw=none,fill=qrgC334155,line width=.7pt] (655.980,367.510) circle[radius=0.0320cm];
\draw[draw=none,fill=qrgC334155,line width=.7pt] (616.020,367.510) circle[radius=0.0320cm];
\draw[draw=none,fill=qrgC334155,line width=.7pt] (603.660,329.490) circle[radius=0.0320cm];
\draw[draw=qrgC16a34a,fill=white,line width=.7pt] (670.000,450.000) circle[radius=0.1840cm];
\node[anchor=center,text=qrgC334155,font={\fontsize{8}{9}\selectfont},inner sep=0pt] at (670.000,450.000) {0};
\draw[draw=qrgC94a3b8,line width=.45pt] (1202.500,306.000) -- (1182.520,367.510);
\draw[draw=qrgC94a3b8,line width=.45pt] (1202.500,306.000) -- (1170.160,329.490);
\draw[draw=qrgC94a3b8,line width=.45pt] (1234.840,329.490) -- (1222.480,367.510);
\draw[draw=qrgC94a3b8,line width=.45pt] (1234.840,329.490) -- (1170.160,329.490);
\draw[draw=qrgC94a3b8,line width=.45pt] (1222.480,367.510) -- (1170.160,329.490);
\draw[draw=qrgC94a3b8,line width=.45pt] (1182.520,367.510) -- (1170.160,329.490);
\draw[draw=none,fill=qrgC334155,line width=.7pt] (1202.500,306.000) circle[radius=0.0320cm];
\draw[draw=none,fill=qrgC334155,line width=.7pt] (1234.840,329.490) circle[radius=0.0320cm];
\draw[draw=none,fill=qrgC334155,line width=.7pt] (1222.480,367.510) circle[radius=0.0320cm];
\draw[draw=none,fill=qrgC334155,line width=.7pt] (1182.520,367.510) circle[radius=0.0320cm];
\draw[draw=none,fill=qrgC334155,line width=.7pt] (1170.160,329.490) circle[radius=0.0320cm];
\draw[draw=qrgC64748b,fill=white,line width=.7pt] (1030.000,450.000) circle[radius=0.1840cm];
\node[anchor=center,text=qrgC334155,font={\fontsize{8}{9}\selectfont},inner sep=0pt] at (1030.000,450.000) {00};
\draw[draw=qrgC64748b,fill=white,line width=.7pt] (1115.000,450.000) circle[radius=0.1840cm];
\node[anchor=center,text=qrgC334155,font={\fontsize{8}{9}\selectfont},inner sep=0pt] at (1115.000,450.000) {11};
\draw[draw=qrgC64748b,fill=white,line width=.7pt] (1200.000,450.000) circle[radius=0.1840cm];
\node[anchor=center,text=qrgC334155,font={\fontsize{8}{9}\selectfont},inner sep=0pt] at (1200.000,450.000) {01};
\draw[draw=qrgC64748b,fill=white,line width=.7pt] (1375.000,450.000) circle[radius=0.1840cm];
\node[anchor=center,text=qrgC334155,font={\fontsize{8}{9}\selectfont},inner sep=0pt] at (1375.000,450.000) {10};
\node[anchor=west,text=qrgC0369a1,font={\fontsize{9}{10}\selectfont},inner sep=0pt] at (1640.500,340.000) {\(K_2\vee3K_1\)};
\draw[draw=qrgC94a3b8,line width=.45pt] (1588.500,306.000) -- (1568.520,367.510);
\draw[draw=qrgC94a3b8,line width=.45pt] (1588.500,306.000) -- (1556.160,329.490);
\draw[draw=qrgC94a3b8,line width=.45pt] (1620.840,329.490) -- (1568.520,367.510);
\draw[draw=qrgC94a3b8,line width=.45pt] (1620.840,329.490) -- (1556.160,329.490);
\draw[draw=qrgC94a3b8,line width=.45pt] (1608.480,367.510) -- (1568.520,367.510);
\draw[draw=qrgC94a3b8,line width=.45pt] (1608.480,367.510) -- (1556.160,329.490);
\draw[draw=qrgC94a3b8,line width=.45pt] (1568.520,367.510) -- (1556.160,329.490);
\draw[draw=none,fill=qrgC334155,line width=.7pt] (1588.500,306.000) circle[radius=0.0320cm];
\draw[draw=none,fill=qrgC334155,line width=.7pt] (1620.840,329.490) circle[radius=0.0320cm];
\draw[draw=none,fill=qrgC334155,line width=.7pt] (1608.480,367.510) circle[radius=0.0320cm];
\draw[draw=none,fill=qrgC334155,line width=.7pt] (1568.520,367.510) circle[radius=0.0320cm];
\draw[draw=none,fill=qrgC334155,line width=.7pt] (1556.160,329.490) circle[radius=0.0320cm];
\draw[draw=qrgC64748b,fill=white,line width=.7pt] (1630.000,450.000) circle[radius=0.1840cm];
\node[anchor=center,text=qrgC334155,font={\fontsize{8}{9}\selectfont},inner sep=0pt] at (1630.000,450.000) {0};
\draw[draw=qrgC64748b,fill=white,line width=.7pt] (1715.000,450.000) circle[radius=0.1840cm];
\node[anchor=center,text=qrgC334155,font={\fontsize{8}{9}\selectfont},inner sep=0pt] at (1715.000,450.000) {1};
\node[anchor=west,text=qrgC0369a1,font={\fontsize{9}{10}\selectfont},inner sep=0pt] at (330.500,640.000) {\(K_5\)};
\draw[draw=qrgC94a3b8,line width=.45pt] (278.500,606.000) -- (310.840,629.490);
\draw[draw=qrgC94a3b8,line width=.45pt] (278.500,606.000) -- (298.480,667.510);
\draw[draw=qrgC94a3b8,line width=.45pt] (278.500,606.000) -- (258.520,667.510);
\draw[draw=qrgC94a3b8,line width=.45pt] (278.500,606.000) -- (246.160,629.490);
\draw[draw=qrgC94a3b8,line width=.45pt] (310.840,629.490) -- (298.480,667.510);
\draw[draw=qrgC94a3b8,line width=.45pt] (310.840,629.490) -- (258.520,667.510);
\draw[draw=qrgC94a3b8,line width=.45pt] (310.840,629.490) -- (246.160,629.490);
\draw[draw=qrgC94a3b8,line width=.45pt] (298.480,667.510) -- (258.520,667.510);
\draw[draw=qrgC94a3b8,line width=.45pt] (298.480,667.510) -- (246.160,629.490);
\draw[draw=qrgC94a3b8,line width=.45pt] (258.520,667.510) -- (246.160,629.490);
\draw[draw=none,fill=qrgC334155,line width=.7pt] (278.500,606.000) circle[radius=0.0320cm];
\draw[draw=none,fill=qrgC334155,line width=.7pt] (310.840,629.490) circle[radius=0.0320cm];
\draw[draw=none,fill=qrgC334155,line width=.7pt] (298.480,667.510) circle[radius=0.0320cm];
\draw[draw=none,fill=qrgC334155,line width=.7pt] (258.520,667.510) circle[radius=0.0320cm];
\draw[draw=none,fill=qrgC334155,line width=.7pt] (246.160,629.490) circle[radius=0.0320cm];
\draw[draw=qrgC64748b,fill=white,line width=.7pt] (270.000,730.000) circle[radius=0.1840cm];
\node[anchor=center,text=qrgC334155,font={\fontsize{8}{9}\selectfont},inner sep=0pt] at (270.000,730.000) {0};
\draw[draw=qrgC64748b,fill=white,line width=.7pt] (355.000,730.000) circle[radius=0.1840cm];
\node[anchor=center,text=qrgC334155,font={\fontsize{8}{9}\selectfont},inner sep=0pt] at (355.000,730.000) {1};
\node[anchor=west,text=qrgC0369a1,font={\fontsize{9}{10}\selectfont},inner sep=0pt] at (348.000,70.000) {\(C_6\)};
\draw[draw=qrgC94a3b8,line width=.45pt] (296.000,36.000) -- (325.440,53.000);
\draw[draw=qrgC94a3b8,line width=.45pt] (325.440,53.000) -- (325.440,87.000);
\draw[draw=qrgC94a3b8,line width=.45pt] (325.440,87.000) -- (296.000,104.000);
\draw[draw=qrgC94a3b8,line width=.45pt] (296.000,104.000) -- (266.560,87.000);
\draw[draw=qrgC94a3b8,line width=.45pt] (266.560,87.000) -- (266.560,53.000);
\draw[draw=qrgC94a3b8,line width=.45pt] (266.560,53.000) -- (296.000,36.000);
\draw[draw=none,fill=qrgC334155,line width=.7pt] (296.000,36.000) circle[radius=0.0320cm];
\draw[draw=none,fill=qrgC334155,line width=.7pt] (325.440,53.000) circle[radius=0.0320cm];
\draw[draw=none,fill=qrgC334155,line width=.7pt] (325.440,87.000) circle[radius=0.0320cm];
\draw[draw=none,fill=qrgC334155,line width=.7pt] (296.000,104.000) circle[radius=0.0320cm];
\draw[draw=none,fill=qrgC334155,line width=.7pt] (266.560,87.000) circle[radius=0.0320cm];
\draw[draw=none,fill=qrgC334155,line width=.7pt] (266.560,53.000) circle[radius=0.0320cm];
\draw[draw=qrgC16a34a,fill=white,line width=.7pt] (330.000,180.000) circle[radius=0.1840cm];
\node[anchor=center,text=qrgC334155,font={\fontsize{8}{9}\selectfont},inner sep=0pt] at (330.000,180.000) {0};
\draw[draw=qrgC94a3b8,line width=.45pt] (712.500,606.000) -- (683.060,657.000);
\draw[draw=qrgC94a3b8,line width=.45pt] (712.500,606.000) -- (683.060,623.000);
\draw[draw=qrgC94a3b8,line width=.45pt] (741.940,623.000) -- (683.060,657.000);
\draw[draw=qrgC94a3b8,line width=.45pt] (741.940,623.000) -- (683.060,623.000);
\draw[draw=qrgC94a3b8,line width=.45pt] (741.940,657.000) -- (712.500,674.000);
\draw[draw=qrgC94a3b8,line width=.45pt] (741.940,657.000) -- (683.060,623.000);
\draw[draw=qrgC94a3b8,line width=.45pt] (712.500,674.000) -- (683.060,623.000);
\draw[draw=none,fill=qrgC334155,line width=.7pt] (712.500,606.000) circle[radius=0.0320cm];
\draw[draw=none,fill=qrgC334155,line width=.7pt] (741.940,623.000) circle[radius=0.0320cm];
\draw[draw=none,fill=qrgC334155,line width=.7pt] (741.940,657.000) circle[radius=0.0320cm];
\draw[draw=none,fill=qrgC334155,line width=.7pt] (712.500,674.000) circle[radius=0.0320cm];
\draw[draw=none,fill=qrgC334155,line width=.7pt] (683.060,657.000) circle[radius=0.0320cm];
\draw[draw=none,fill=qrgC334155,line width=.7pt] (683.060,623.000) circle[radius=0.0320cm];
\draw[draw=qrgC64748b,fill=white,line width=.7pt] (670.000,730.000) circle[radius=0.1840cm];
\node[anchor=center,text=qrgC334155,font={\fontsize{8}{9}\selectfont},inner sep=0pt] at (670.000,730.000) {0};
\draw[draw=qrgC64748b,fill=white,line width=.7pt] (755.000,730.000) circle[radius=0.1840cm];
\node[anchor=center,text=qrgC334155,font={\fontsize{8}{9}\selectfont},inner sep=0pt] at (755.000,730.000) {1};
\node[anchor=west,text=qrgC0369a1,font={\fontsize{9}{10}\selectfont},inner sep=0pt] at (1070.500,640.000) {\(K_{2,4}\)};
\draw[draw=qrgC94a3b8,line width=.45pt] (1018.500,606.000) -- (989.060,657.000);
\draw[draw=qrgC94a3b8,line width=.45pt] (1018.500,606.000) -- (989.060,623.000);
\draw[draw=qrgC94a3b8,line width=.45pt] (1047.940,623.000) -- (989.060,657.000);
\draw[draw=qrgC94a3b8,line width=.45pt] (1047.940,623.000) -- (989.060,623.000);
\draw[draw=qrgC94a3b8,line width=.45pt] (1047.940,657.000) -- (989.060,657.000);
\draw[draw=qrgC94a3b8,line width=.45pt] (1047.940,657.000) -- (989.060,623.000);
\draw[draw=qrgC94a3b8,line width=.45pt] (1018.500,674.000) -- (989.060,657.000);
\draw[draw=qrgC94a3b8,line width=.45pt] (1018.500,674.000) -- (989.060,623.000);
\draw[draw=none,fill=qrgC334155,line width=.7pt] (1018.500,606.000) circle[radius=0.0320cm];
\draw[draw=none,fill=qrgC334155,line width=.7pt] (1047.940,623.000) circle[radius=0.0320cm];
\draw[draw=none,fill=qrgC334155,line width=.7pt] (1047.940,657.000) circle[radius=0.0320cm];
\draw[draw=none,fill=qrgC334155,line width=.7pt] (1018.500,674.000) circle[radius=0.0320cm];
\draw[draw=none,fill=qrgC334155,line width=.7pt] (989.060,657.000) circle[radius=0.0320cm];
\draw[draw=none,fill=qrgC334155,line width=.7pt] (989.060,623.000) circle[radius=0.0320cm];
\draw[draw=qrgC64748b,fill=white,line width=.7pt] (1070.000,730.000) circle[radius=0.1840cm];
\node[anchor=center,text=qrgC334155,font={\fontsize{8}{9}\selectfont},inner sep=0pt] at (1070.000,730.000) {0};
\draw[draw=qrgC94a3b8,line width=.45pt] (1512.500,606.000) -- (1483.060,657.000);
\draw[draw=qrgC94a3b8,line width=.45pt] (1512.500,606.000) -- (1483.060,623.000);
\draw[draw=qrgC94a3b8,line width=.45pt] (1541.940,623.000) -- (1483.060,657.000);
\draw[draw=qrgC94a3b8,line width=.45pt] (1541.940,623.000) -- (1483.060,623.000);
\draw[draw=qrgC94a3b8,line width=.45pt] (1541.940,657.000) -- (1512.500,674.000);
\draw[draw=qrgC94a3b8,line width=.45pt] (1541.940,657.000) -- (1483.060,623.000);
\draw[draw=qrgC94a3b8,line width=.45pt] (1512.500,674.000) -- (1483.060,657.000);
\draw[draw=qrgC94a3b8,line width=.45pt] (1483.060,657.000) -- (1483.060,623.000);
\draw[draw=none,fill=qrgC334155,line width=.7pt] (1512.500,606.000) circle[radius=0.0320cm];
\draw[draw=none,fill=qrgC334155,line width=.7pt] (1541.940,623.000) circle[radius=0.0320cm];
\draw[draw=none,fill=qrgC334155,line width=.7pt] (1541.940,657.000) circle[radius=0.0320cm];
\draw[draw=none,fill=qrgC334155,line width=.7pt] (1512.500,674.000) circle[radius=0.0320cm];
\draw[draw=none,fill=qrgC334155,line width=.7pt] (1483.060,657.000) circle[radius=0.0320cm];
\draw[draw=none,fill=qrgC334155,line width=.7pt] (1483.060,623.000) circle[radius=0.0320cm];
\draw[draw=qrgC64748b,fill=white,line width=.7pt] (1470.000,730.000) circle[radius=0.1840cm];
\node[anchor=center,text=qrgC334155,font={\fontsize{8}{9}\selectfont},inner sep=0pt] at (1470.000,730.000) {0};
\draw[draw=qrgC64748b,fill=white,line width=.7pt] (1555.000,730.000) circle[radius=0.1840cm];
\node[anchor=center,text=qrgC334155,font={\fontsize{8}{9}\selectfont},inner sep=0pt] at (1555.000,730.000) {1};
\draw[draw=qrgC94a3b8,line width=.45pt] (312.500,886.000) -- (283.060,937.000);
\draw[draw=qrgC94a3b8,line width=.45pt] (312.500,886.000) -- (283.060,903.000);
\draw[draw=qrgC94a3b8,line width=.45pt] (341.940,903.000) -- (312.500,954.000);
\draw[draw=qrgC94a3b8,line width=.45pt] (341.940,903.000) -- (283.060,903.000);
\draw[draw=qrgC94a3b8,line width=.45pt] (341.940,937.000) -- (312.500,954.000);
\draw[draw=qrgC94a3b8,line width=.45pt] (341.940,937.000) -- (283.060,937.000);
\draw[draw=qrgC94a3b8,line width=.45pt] (312.500,954.000) -- (283.060,937.000);
\draw[draw=qrgC94a3b8,line width=.45pt] (312.500,954.000) -- (283.060,903.000);
\draw[draw=qrgC94a3b8,line width=.45pt] (283.060,937.000) -- (283.060,903.000);
\draw[draw=none,fill=qrgC334155,line width=.7pt] (312.500,886.000) circle[radius=0.0320cm];
\draw[draw=none,fill=qrgC334155,line width=.7pt] (341.940,903.000) circle[radius=0.0320cm];
\draw[draw=none,fill=qrgC334155,line width=.7pt] (341.940,937.000) circle[radius=0.0320cm];
\draw[draw=none,fill=qrgC334155,line width=.7pt] (312.500,954.000) circle[radius=0.0320cm];
\draw[draw=none,fill=qrgC334155,line width=.7pt] (283.060,937.000) circle[radius=0.0320cm];
\draw[draw=none,fill=qrgC334155,line width=.7pt] (283.060,903.000) circle[radius=0.0320cm];
\draw[draw=qrgC64748b,fill=white,line width=.7pt] (270.000,1010.000) circle[radius=0.1840cm];
\node[anchor=center,text=qrgC334155,font={\fontsize{8}{9}\selectfont},inner sep=0pt] at (270.000,1010.000) {0};
\draw[draw=qrgC64748b,fill=white,line width=.7pt] (355.000,1010.000) circle[radius=0.1840cm];
\node[anchor=center,text=qrgC334155,font={\fontsize{8}{9}\selectfont},inner sep=0pt] at (355.000,1010.000) {1};
\draw[draw=qrgC94a3b8,line width=.45pt] (712.500,886.000) -- (683.060,937.000);
\draw[draw=qrgC94a3b8,line width=.45pt] (712.500,886.000) -- (683.060,903.000);
\draw[draw=qrgC94a3b8,line width=.45pt] (741.940,903.000) -- (741.940,937.000);
\draw[draw=qrgC94a3b8,line width=.45pt] (741.940,903.000) -- (712.500,954.000);
\draw[draw=qrgC94a3b8,line width=.45pt] (741.940,937.000) -- (712.500,954.000);
\draw[draw=qrgC94a3b8,line width=.45pt] (741.940,937.000) -- (683.060,937.000);
\draw[draw=qrgC94a3b8,line width=.45pt] (741.940,937.000) -- (683.060,903.000);
\draw[draw=qrgC94a3b8,line width=.45pt] (712.500,954.000) -- (683.060,937.000);
\draw[draw=qrgC94a3b8,line width=.45pt] (712.500,954.000) -- (683.060,903.000);
\draw[draw=qrgC94a3b8,line width=.45pt] (683.060,937.000) -- (683.060,903.000);
\draw[draw=none,fill=qrgC334155,line width=.7pt] (712.500,886.000) circle[radius=0.0320cm];
\draw[draw=none,fill=qrgC334155,line width=.7pt] (741.940,903.000) circle[radius=0.0320cm];
\draw[draw=none,fill=qrgC334155,line width=.7pt] (741.940,937.000) circle[radius=0.0320cm];
\draw[draw=none,fill=qrgC334155,line width=.7pt] (712.500,954.000) circle[radius=0.0320cm];
\draw[draw=none,fill=qrgC334155,line width=.7pt] (683.060,937.000) circle[radius=0.0320cm];
\draw[draw=none,fill=qrgC334155,line width=.7pt] (683.060,903.000) circle[radius=0.0320cm];
\draw[draw=qrgC64748b,fill=white,line width=.7pt] (670.000,1010.000) circle[radius=0.1840cm];
\node[anchor=center,text=qrgC334155,font={\fontsize{8}{9}\selectfont},inner sep=0pt] at (670.000,1010.000) {0};
\draw[draw=qrgC64748b,fill=white,line width=.7pt] (755.000,1010.000) circle[radius=0.1840cm];
\node[anchor=center,text=qrgC334155,font={\fontsize{8}{9}\selectfont},inner sep=0pt] at (755.000,1010.000) {1};
\draw[draw=qrgC94a3b8,line width=.45pt] (1112.500,886.000) -- (1083.060,937.000);
\draw[draw=qrgC94a3b8,line width=.45pt] (1112.500,886.000) -- (1083.060,903.000);
\draw[draw=qrgC94a3b8,line width=.45pt] (1141.940,903.000) -- (1141.940,937.000);
\draw[draw=qrgC94a3b8,line width=.45pt] (1141.940,903.000) -- (1112.500,954.000);
\draw[draw=qrgC94a3b8,line width=.45pt] (1141.940,903.000) -- (1083.060,937.000);
\draw[draw=qrgC94a3b8,line width=.45pt] (1141.940,903.000) -- (1083.060,903.000);
\draw[draw=qrgC94a3b8,line width=.45pt] (1141.940,937.000) -- (1112.500,954.000);
\draw[draw=qrgC94a3b8,line width=.45pt] (1141.940,937.000) -- (1083.060,937.000);
\draw[draw=qrgC94a3b8,line width=.45pt] (1141.940,937.000) -- (1083.060,903.000);
\draw[draw=qrgC94a3b8,line width=.45pt] (1112.500,954.000) -- (1083.060,937.000);
\draw[draw=qrgC94a3b8,line width=.45pt] (1112.500,954.000) -- (1083.060,903.000);
\draw[draw=none,fill=qrgC334155,line width=.7pt] (1112.500,886.000) circle[radius=0.0320cm];
\draw[draw=none,fill=qrgC334155,line width=.7pt] (1141.940,903.000) circle[radius=0.0320cm];
\draw[draw=none,fill=qrgC334155,line width=.7pt] (1141.940,937.000) circle[radius=0.0320cm];
\draw[draw=none,fill=qrgC334155,line width=.7pt] (1112.500,954.000) circle[radius=0.0320cm];
\draw[draw=none,fill=qrgC334155,line width=.7pt] (1083.060,937.000) circle[radius=0.0320cm];
\draw[draw=none,fill=qrgC334155,line width=.7pt] (1083.060,903.000) circle[radius=0.0320cm];
\draw[draw=qrgC64748b,fill=white,line width=.7pt] (1070.000,1010.000) circle[radius=0.1840cm];
\node[anchor=center,text=qrgC334155,font={\fontsize{8}{9}\selectfont},inner sep=0pt] at (1070.000,1010.000) {0};
\draw[draw=qrgC64748b,fill=white,line width=.7pt] (1155.000,1010.000) circle[radius=0.1840cm];
\node[anchor=center,text=qrgC334155,font={\fontsize{8}{9}\selectfont},inner sep=0pt] at (1155.000,1010.000) {1};
\node[anchor=west,text=qrgC0369a1,font={\fontsize{9}{10}\selectfont},inner sep=0pt] at (1500.500,920.000) {\(K_{2,2,2}\)};
\draw[draw=qrgC94a3b8,line width=.45pt] (1448.500,886.000) -- (1477.940,937.000);
\draw[draw=qrgC94a3b8,line width=.45pt] (1448.500,886.000) -- (1448.500,954.000);
\draw[draw=qrgC94a3b8,line width=.45pt] (1448.500,886.000) -- (1419.060,937.000);
\draw[draw=qrgC94a3b8,line width=.45pt] (1448.500,886.000) -- (1419.060,903.000);
\draw[draw=qrgC94a3b8,line width=.45pt] (1477.940,903.000) -- (1477.940,937.000);
\draw[draw=qrgC94a3b8,line width=.45pt] (1477.940,903.000) -- (1448.500,954.000);
\draw[draw=qrgC94a3b8,line width=.45pt] (1477.940,903.000) -- (1419.060,937.000);
\draw[draw=qrgC94a3b8,line width=.45pt] (1477.940,903.000) -- (1419.060,903.000);
\draw[draw=qrgC94a3b8,line width=.45pt] (1477.940,937.000) -- (1419.060,937.000);
\draw[draw=qrgC94a3b8,line width=.45pt] (1477.940,937.000) -- (1419.060,903.000);
\draw[draw=qrgC94a3b8,line width=.45pt] (1448.500,954.000) -- (1419.060,937.000);
\draw[draw=qrgC94a3b8,line width=.45pt] (1448.500,954.000) -- (1419.060,903.000);
\draw[draw=none,fill=qrgC334155,line width=.7pt] (1448.500,886.000) circle[radius=0.0320cm];
\draw[draw=none,fill=qrgC334155,line width=.7pt] (1477.940,903.000) circle[radius=0.0320cm];
\draw[draw=none,fill=qrgC334155,line width=.7pt] (1477.940,937.000) circle[radius=0.0320cm];
\draw[draw=none,fill=qrgC334155,line width=.7pt] (1448.500,954.000) circle[radius=0.0320cm];
\draw[draw=none,fill=qrgC334155,line width=.7pt] (1419.060,937.000) circle[radius=0.0320cm];
\draw[draw=none,fill=qrgC334155,line width=.7pt] (1419.060,903.000) circle[radius=0.0320cm];
\draw[draw=qrgC64748b,fill=white,line width=.7pt] (1470.000,1010.000) circle[radius=0.1840cm];
\node[anchor=center,text=qrgC334155,font={\fontsize{8}{9}\selectfont},inner sep=0pt] at (1470.000,1010.000) {0};
\draw[draw=qrgC64748b,fill=white,line width=.7pt] (1555.000,1010.000) circle[radius=0.1840cm];
\node[anchor=center,text=qrgC334155,font={\fontsize{8}{9}\selectfont},inner sep=0pt] at (1555.000,1010.000) {1};
\end{tikzpicture}
% END GENERATED QRG EULERIAN
\caption{All Eulerian QRGs maps between connected graphs with 2-6 vertices.
Circle labels give \(\alpha\)-parameters \(0,1\) or \(00,01,10,11\);
zero denotes the inner connection. Other conventions are as in
Figure~\ref{fig:QRG-quiver-co-eulerian}.}
\label{fig:QRG-quiver-eulerian}
\end{figure}
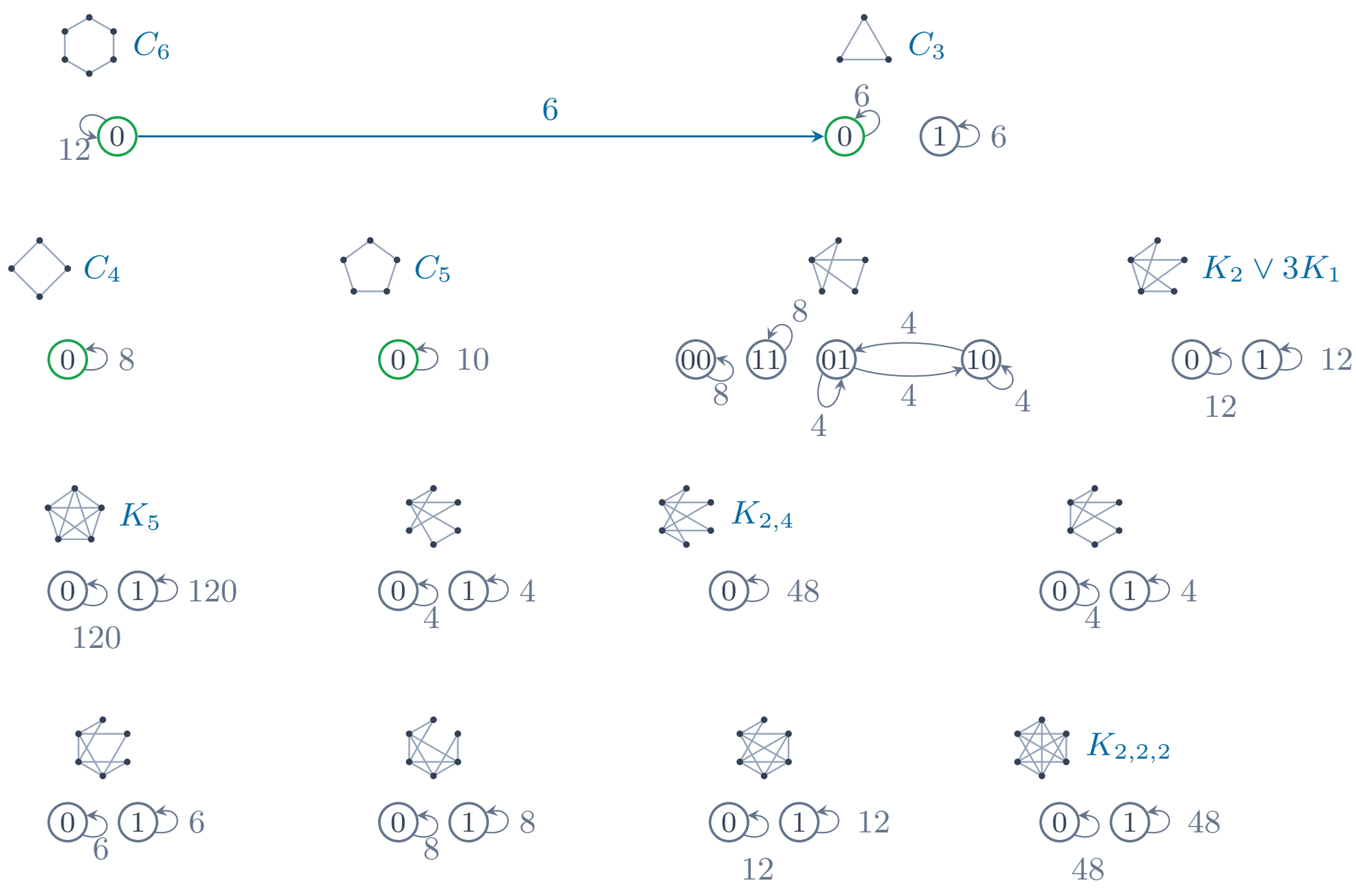

By Theorem~\ref{thm:binary-QRG-map}, this quiver is the disjoint union
of its co-Eulerian and Eulerian parts, shown in
Figures~\ref{fig:QRG-quiver-co-eulerian} and~\ref{fig:QRG-quiver-eulerian}.
They have \(11\) and \(26\) objects, respectively. Circles near a
graph drawing are its QRGs of the indicated type, with their
\(\alpha\)-parameters inside. Each drawn arrow represents all maps
between its ordered pair of endpoints, with the multiplicity beside
it; loops include identity maps. The six non-isomorphisms in the
Eulerian part are the double covers \(C_6\to C_3\), landing on
its inner flat geometry. As suggested by the figure, Eulerian morphisms between different underlying graphs appear to be sparse: through nine vertices the only similar possibilities are the
inner cycle coverings \(C_6\to C_3\), \(C_8\to C_4\) and
\(C_9\to C_3\). We also see that within a cluster, most morphisms are QRG automorphisms for the same
connection. To go between distinct
connections would require a nontrivial action of graph automorphisms on
\(V_*(G)\), which is not possible when \(d_*(G)=1\). The condition \(d_*(G)=1\) holds for almost every Eulerian graph as \(n\to\infty\), by
Theorem~\ref{thm:typical-eulerian}. Corollary~\ref{cor:typical-map-rigidity}
strengthens this for the non-inner geometry: for almost every Eulerian
graph, every incoming map from a connected star-preserving QRG on at
least two vertices is an
isomorphism, even when larger source graphs are allowed.

\appendix 

\section{Towards physics over \(\F_2\)}
\label{app:phys}

Here we look at the results of the paper in the context of developing physics over $\F_2$ using QRG\cite{MaPac:cla,MaPac:dig,Ma:boo}. Although not conclusive, we provide examples of real statistical weights on
the space of binary geometries and of conserved curvature tensors within
the geometry over \(\F_2\).

\subsection{A partition function on Eulerian geometries}
\label{sec:part}

In quantum gravity it is typical to fix a manifold and write down a partition function where one attempts to `integrate' over all metrics with a suitable action. However, one should in principle also `sum' over all manifolds. For QRGs over $\F_2$ on graphs, we have no metric freedom on a given graph $G$ (for a fixed braiding, the star-preserving QLCs are parametrised by $\alpha\in V_*(G)$) but we can sum over all graphs. We have seen in Proposition~\ref{prop:qlc-asymptotic} that Eulerian graphs dominate among large graphs admitting QLCs, so we are led to a partition function of the form
\begin{equation}
\label{eq:binary-partition}
 Z_n(\beta)=
 \sum_{G\ {\rm connected\ Eulerian}\atop 4\le |G|\le n} e^{-\beta\mathcal S(G)} 
 \sum_{\alpha\in V_*(G)}e^{-\beta S_G(\alpha)},
\end{equation}
where we sum over connected Eulerian graph isomorphism classes with up to $n$ vertices, excluding $|G|<4$ as too special. We think of $n$ as a cutoff and would eventually like to take $n\to \infty$. Moreover, $\beta$ is a real coupling constant with $\beta\to 0$ the `strong gravity' limit. The simplest choice for the other action is  $S_G(\alpha)=0$, so the inner sum is $2^{d^*(G)}$. Theorem~\ref{thm:typical-eulerian} shows that almost every connected Eulerian graph has $d^*(G)=1$. Conjecture~\ref{conj:star-asymptotic} predicts that the connection sum contributes an overall factor of two asymptotically under uniform graph weighting, which drops out of normalised graph observables. This motivates omitting the inner sum as an approximation to the physical model, provided the action does not enhance the exceptional graphs enough to offset their rarity.

For $\mathcal S(G)$ we propose to use the Lin--Lu--Yau Ricci curvature, in Smith's formulation~\cite{Smi}. We give each edge of \(G\) unit length and for adjacent vertices $x\to y$
\[
 R_{xy}:=\lim_{\varepsilon\downarrow0}
 \frac{1-W_1(m_x^\varepsilon,m_y^\varepsilon)}{\varepsilon},
 \qquad
 m_x^\varepsilon=(1-\varepsilon)\delta_x+
 \frac{\varepsilon}{\deg(x)}\sum_{x\to z}\delta_z,
\]
Here the functions \(m_x^\varepsilon\) are probability masses on
vertices. With \(d(u,v)\) the shortest-path distance, their transport
distance is the finite linear programme
\[
 W_1(m_x^\varepsilon,m_y^\varepsilon)
 =\min_{\pi\geq0}\sum_{u,v\in V}d(u,v)\pi(u,v),
 \qquad
 \sum_v\pi(u,v)=m_x^\varepsilon(u),\quad
 \sum_u\pi(u,v)=m_y^\varepsilon(v).
\]
Thus \(\pi(u,v)\) is the mass transported from \(u\) to \(v\).
Only sources in \(\{x\}\cup N(x)\) and destinations in
\(\{y\}\cup N(y)\) occur. These vertices are all within two steps
of \(x\), but a source and destination can be three steps apart:
for \(u\in N(x)\), \(v\in N(y)\), the path
\(u\to x\to y\to v\) has length three. The relevant costs are
therefore just \(0,1,2,3\): zero at coincident vertices, one for
an edge, two for a nonadjacent pair with a common neighbour, and
three otherwise. Such a common neighbour need not itself carry
mass, but lies within two steps of \(x\). Hence the induced
radius-two neighbourhood of \(x\) suffices for the calculation.

For example, on a cycle let \(u\) and \(v\) be the other neighbours
of adjacent vertices \(x\) and \(y\). For small \(\varepsilon\),
an optimal plan keeps the common masses in place, transports
\(1-3\varepsilon/2\) from \(x\) to \(y\), and
\(\varepsilon/2\) from \(u\) to \(v\). The latter transport is
unnecessary when \(u=v\). Consequently
\[
 W_1=1-\frac{\varepsilon}{2}(3-d(u,v)),
 \qquad R_{xy}=\frac{3-d(u,v)}2.
\]
The distances are \(0,1,2,3\) for \(C_3,C_4,C_5,C_n\)
with \(n\geq6\), respectively.

We view this curvature as a central tensor ${\rm Ricci}:=\sum_{x\to y}R_{xy}\omega_{xyx}$ on a graph differential calculus over $\R$. The scalar curvature at a vertex $x$ is  $\sum_{x\to y}R_{xy}$ and summing over $x$ (i.e. with counting measure on vertices) suggests a graph Einstein-Hilbert action
\begin{equation}
\label{eq:smith-action}
 \mathcal S(G)=\sum_{x\to y}R_{xy}.
\end{equation}
For example, we can show that 
\(\mathcal S(C_3)=9\), \(\mathcal S(C_4)=8\),
\(\mathcal S(C_5)=5\), and \(\mathcal S(C_n)=0\) for \(n\geq6\), while \(\mathcal S(G)=-4|G|\) for any four-regular
graph of girth at least six. Note that if we had varying lengths then we could apply the inverse metric in computing the scalar curvature as well as use the metric in the transport problem. We are also free to take some other measure on the vertices. It is not known, however, how such an Einstein-Hilbert action could arise from the curvature within QRG over $\R$. 

As a first look at quantum gravity, we can consider as an observable the averaged scalar curvature itself, i.e. $\mathcal S(G)/|G|$. We take $S_G(\alpha)=0$, with results shown on the left in Figure~\ref{fig:binary-scalar-moments}.  Because the action can be negative, it is not clear that we can have stable values as $n\to \infty$ even with renormalisation. It could also be of interest to look at the `Lorentzian' theory with $\imath=\sqrt{-1}$ in front of $\beta$ in the weight. The standard deviation in the curvature has been proposed\cite{Whe,BliMa:cur} as a measure of vacuum energy from spacetime fluctuations, albeit arising in our case from fluctuation of the underlying graph rather than of the metric.

\subsection{Conserved curvature tensors over \(\F_2\)}

A separate question is to construct an Einstein tensor
\(\mathsf G\in\Omega^1\ten\Omega^1\) from curvature working in our case of a QRG over $\F_2$. We require this to be built from the curvature and to be conserved, for example in the sense
\[
 \operatorname{Div}\mathsf G:=
 ((\ ,\ )\ten\id)\nabla\mathsf G=0.
\]
For a central tensor \(H=\sum_{x\to y}h_{xy}\omega_{xyx}\),
expansion of the induced connection gives
\begin{equation}
\label{eq:central-divergence}
 (\operatorname{Div}H)_{xz}
 =h_{xz}+h_{zx}+(q+1)\sum_{x\to y}h_{xy}.
\end{equation}
where we include the last term in the Eulerian $q=0$ case. The \(\alpha\)-terms do not survive the contraction. Thus an
edge-symmetric tensor is conserved in an Eulerian geometry precisely
when each row sum vanishes; in a co-Eulerian geometry every such tensor
is conserved.

For an Eulerian QLC a curvature-derived example is
\begin{equation}
\label{eq:conserved-curvature-tensor}
 H_\alpha=\sum_{x\to y}
 \operatorname{codeg}(x,y)a_{xy}a_{yx}\,\omega_{xyx},
 \qquad \operatorname{Div}H_\alpha=0,
\end{equation}
where \(\alpha(\omega_{xy})=a_{xy}P^{\mathrm{sum}}_{xy}\), with
\(a_{xy}=0\) on edges in no triangle. The triangle equations make
\(a_{xy}a_{yx}\) constant on the three edges of each triangle, so
each incident triangle contributes twice to a row sum. To see the
curvature origin, let \(r(v)\) represent the return block of
\(R_\nabla(\omega_{xy})\) with final arrow \(x\to y\). The return-block
formula in the proof of Theorem~\ref{thm:flat-alpha} gives
\[
 \sum_{v\in N(x)\cap N(y)}(r(v)+r(y))
 =\operatorname{codeg}(x,y)a_{xy}a_{yx}.
\]
The expression on the left is independent of the representative.
This tensor vanishes for every inner connection, including curved ones.

Allowing noncentral tensors retains more curvature information. For a
bimodule lift \(i:\Omega^2_{\min}\to\Omega^1\ten\Omega^1\), use
the usual contraction~\cite{BegMa}
\[
 \Ric_i=((\ ,\ )\ten\id\ten\id)
 (\id\ten i\ten\id)(\id\ten R_\nabla)g,
 \qquad \Scal_i=(\ ,\ )\Ric_i.
\]
Every lift has, on a nonempty endpoint block, the form
\[
 i(\omega_{xv}\wedge\omega_{vb})
 =\omega_{xvb}+\lambda_v^{xb}P^{\mathrm{sum}}_{xb},
 \qquad \sum_{v\in N(x)\cap N(b)}\lambda_v^{xb}=1.
\]
The usual averaging lift over \(\C\) divides by the block size
\(m_{xb}=|N(x)\cap N(b)|\). Over \(\F_2\), it therefore fails
on every block of even size. In particular, every nonempty return
block of an Eulerian graph has even size \(m_{xx}=\deg(x)\).
A lift respects \(*\) precisely when
\[
 \lambda_v^{xb}=\lambda_v^{bx}.
\]
Such lifts always exist: choose a coefficient vector of sum one for
each return block and each pair of reversed non-return blocks.
Write \(\mathcal I_*(G)\) for this finite set of lifts.

For the inner Eulerian connection $\nabla_0$, the curvature coefficient in
\eqref{eq:flat-curvature-coeff} reduces to
\(\delta_{vy}+\delta_{vc}+\delta_{xb}\). Substitution gives
\[
 \Ric_i=\sum_{t\to b\to c}r_{tbc}\omega_{tbc},
 \qquad
 r_{tbc}=\sum_{x\to t}(\lambda_t^{xb}+\lambda_c^{xb}),
\]
where \(\lambda\) is extended by zero outside its block.
Consequently \(r_{tbt}=0\) and, for fixed \(t,b\),
\(\sum_{c:b\to c}r_{tbc}=0\).
The first identity implies \(\Scal_i=0\) for every lift, without
forcing \(\Ric_i=0\). The divergence reduces to
\( (\operatorname{Div}\Ric_i)_{st}=\sum_{u\to s}r_{ust}\);
hence \(\Ric_i^\dagger=\Ric_i\) implies conservation in this case.

For example, on the `butterfly' consisting of two triangles \(012\) and \(034\) sharing
vertex \(0\), choose \(\lambda^{xb}\) supported at the least
intermediate vertex in each block. Then the inner connection has
\[
 \Ric_i=\sum_{\substack{u,v\in\{2,3,4\}\\u\ne v}}
              \omega_{u0v}\ne0,
 \qquad \Ric_i^\dagger=\Ric_i,
 \qquad \operatorname{Div}\Ric_i=0.
\]
In contrast, for the triangle \(012\) and square \(03450\) sharing
\(0\), substitution into the same formulas gives
\(\operatorname{Div}\Ric_i=\omega_{01}+\omega_{02}\) for every
lift. Since \(\Scal_i=0\), a scalar-curvature multiple of the metric
cannot correct this defect.

More generally, for a star-preserving Eulerian geometry and a star-compatible lift,
\begin{equation}
 \Scal_i(y)=\sum_{x\to y}a_{xy}
       \sum_{v\in N(x)\cap N(y)}\lambda_v^{xx}.
 \label{eq:binary-scalar}
\end{equation}
This is because the scalar contraction selects the coefficient of
\(\omega_{xv}\wedge\omega_{vb}\ten\omega_{by}\) in
\(R_\nabla(\omega_{xy})\). Equation~\eqref{eq:flat-curvature-coeff}
gives \(\delta_{xb}+a_{xy}A_{vy}\). The constant part cancels
under lifting, while the terms with \(x\ne b\) cancel in reversed
pairs: star compatibility gives \(\lambda_v^{xb}=\lambda_v^{bx}\),
and the triangles \(xyv\), \(byv\) give \(a_{xy}=a_{by}\).
Only the return-block parameters remain.

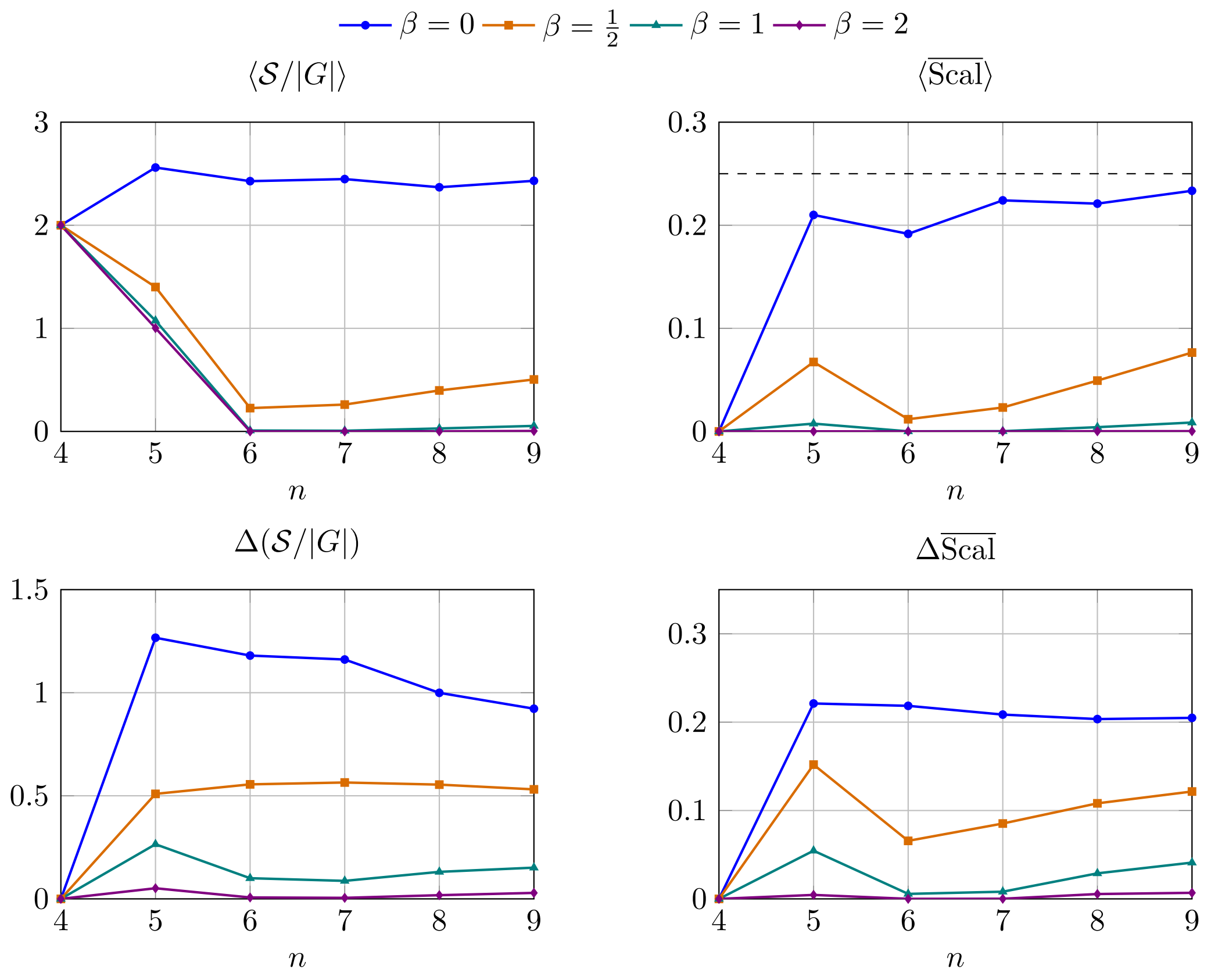
\begin{figure}[tbp]
\centering
\begin{tikzpicture}
\begin{groupplot}[
 group style={group size=2 by 2,horizontal sep=.13\textwidth,vertical sep=1.8cm},
 width=.43\textwidth,height=5.1cm,
 xmin=4,xmax=9,xtick={4,5,6,7,8,9},xlabel={$n$},
 ymin=0,grid=major,
 tick label style={font=\small},label style={font=\small},
 title style={font=\small},
 legend style={font=\small,draw=none,fill=white,cells={anchor=west}},
 legend columns=4]
\nextgroupplot[title={$\langle\mathcal S/|G|\rangle$},ymax=3,ytick={0,1,2,3},legend to name=leg:scalar-beta]
\addplot[blue,thick,mark=*,mark size=1pt] coordinates {
(4,2) (5,2.56) (6,2.42777777777778) (7,2.44801920768307) (8,2.36847218747394) (9,2.4312210741739)};
\addlegendentry{$\beta=0$}
\addplot[orange!85!black,thick,mark=square*,mark size=1pt] coordinates {
(4,2) (5,1.4006153953787) (6,0.226207712202379) (7,0.260351073592056) (8,0.397427961025474) (9,0.503965759818578)};
\addlegendentry{$\beta=\frac12$}
\addplot[teal,thick,mark=triangle*,mark size=1pt] coordinates {
(4,2) (5,1.07469148694154) (6,0.00863694247299414) (7,0.00659290682410054) (8,0.0287466249418923) (9,0.0535478672339211)};
\addlegendentry{$\beta=1$}
\addplot[violet,thick,mark=diamond*,mark size=1pt] coordinates {
(4,2) (5,1.00265914247637) (6,4.58366368511463e-05) (7,2.42969615359033e-05) (8,0.000769728806814082) (9,0.00417908397575291)};
\addlegendentry{$\beta=2$}
\nextgroupplot[title={$\langle\overline{\Scal}\rangle$},ymax=.3,ytick={0,.1,.2,.3}]
\addplot[black,dashed,forget plot] coordinates {(4,.25) (9,.25)};
\addplot[blue,thick,mark=*,mark size=1pt] coordinates {
(4,0) (5,0.21) (6,0.191666666666667) (7,0.224129651860744) (8,0.22096572429322) (9,0.233449535789158)};
\addplot[orange!85!black,thick,mark=square*,mark size=1pt] coordinates {
(4,0) (5,0.0673394336660071) (6,0.0118206165202227) (7,0.0232201022968132) (8,0.0492919605766773) (9,0.076508478177125)};
\addplot[teal,thick,mark=triangle*,mark size=1pt] coordinates {
(4,0) (5,0.0075605413683636) (6,9.15406520297538e-05) (7,0.000236550894341695) (8,0.00412370869487246) (9,0.00861875354024648)};
\addplot[violet,thick,mark=diamond*,mark size=1pt] coordinates {
(4,0) (5,5.09343538804872e-05) (6,6.21975811942441e-09) (7,9.9999158337401e-08) (8,0.000158857351534052) (9,0.000250113995681185)};
\nextgroupplot[title={$\Delta(\mathcal S/|G|)$},ymax=1.5,ytick={0,.5,1,1.5}]
\addplot[blue,thick,mark=*,mark size=1pt] coordinates {
(4,0) (5,1.26664912268552) (6,1.18016660259327) (7,1.16096254080977) (8,0.99935199155967) (9,0.922415414538178)};
\addplot[orange!85!black,thick,mark=square*,mark size=1pt] coordinates {
(4,0) (5,0.509270265821066) (6,0.555434825463093) (7,0.564472382472594) (8,0.554161831495851) (9,0.531147521312062)};
\addplot[teal,thick,mark=triangle*,mark size=1pt] coordinates {
(4,0) (5,0.26495053768577) (6,0.100202975691303) (7,0.087269451523653) (8,0.131030084496404) (9,0.15143077760552)};
\addplot[violet,thick,mark=diamond*,mark size=1pt] coordinates {
(4,0) (5,0.0515153356658357) (6,0.00679340995984568) (7,0.00494785093714181) (8,0.0178606859186971) (9,0.0289438608442598)};
\nextgroupplot[title={$\Delta\overline{\Scal}$},ymax=.35,ytick={0,.1,.2,.3}]
\addplot[blue,thick,mark=*,mark size=1pt] coordinates {
(4,0) (5,0.22113344387496) (6,0.218475356691508) (7,0.208574901248895) (8,0.203482808403184) (9,0.204855985748687)};
\addplot[orange!85!black,thick,mark=square*,mark size=1pt] coordinates {
(4,0) (5,0.151949815380121) (6,0.0656602360772818) (7,0.0852817480143121) (8,0.108128602740194) (9,0.121511419008245)};
\addplot[teal,thick,mark=triangle*,mark size=1pt] coordinates {
(4,0) (5,0.0545368849887813) (6,0.00563766844403949) (7,0.00811199595543648) (8,0.0288953227648609) (9,0.0410591035290805)};
\addplot[violet,thick,mark=diamond*,mark size=1pt] coordinates {
(4,0) (5,0.00447560841419608) (6,4.35341851174874e-05) (7,0.000149169863427612) (8,0.00547818378083984) (9,0.0068453881374381)};
\end{groupplot}
\node[anchor=south] at
 ($(group c1r1.north west)!0.5!(group c2r1.north east)+(0,.65cm)$)
 {\ref{leg:scalar-beta}};
\end{tikzpicture}
\caption{Scalar curvature means (top) and standard deviations (bottom) for quantum gravity partition function with \(S_G(\alpha)=0\) and truncated to \(4\le|G|\le n\). The dashed line is the  \(\beta=0\) asymptote for \(\langle\overline{\Scal}\rangle\) implied by Conjecture~\ref{conj:star-asymptotic}.}
\label{fig:binary-scalar-moments}
\end{figure}

We can still use this $\F_2$-scalar curvature as a real observable on the collection of QRGs. Let
\(\iota:\F_2\to\{0,1\}\subset\R\) take the usual integer
representatives. As there is no preferred lift, we average over
\(\mathcal I_*(G)\) and put
\[
 \overline{\Scal}(G,\alpha)=
 \frac1{|G|\,|\mathcal I_*(G)|}
 \sum_{i\in\mathcal I_*(G)}\sum_{x\in V}\iota(\Scal_i(x))
\]
as a non-negative curvature measure from the $\F_2$ geometry.
The averages are over \(\R\). Averaging over vertices and lifts makes
this observable invariant under relabelling, without needing to match
vertices of different graphs. For the inner Eulerian connection,
\(\Scal_i=0\) for every lift, so \(\overline{\Scal}=0\).

Results for this observable with the same partition function
\(Z_n(\beta)\) as in Section~\ref{sec:part}, using the Smith action
\(\mathcal S(G)\) and \(S_G(\alpha)=0\), are shown on the right in
Figure~\ref{fig:binary-scalar-moments}. All moments are taken over
graphs and connections after the lift average, with
\(\Delta\overline{\Scal}=\sqrt{\langle\overline{\Scal}^2\rangle
-\langle\overline{\Scal}\rangle^2}\).
It is remarkable that the behaviour is broadly comparable with the
Smith observable, although the plots use different vertical scales.
Alternatively, we could include the lift as geometric data in the
partition sum, with the normalised average
\(|\mathcal I_*(G)|^{-1}\sum_i\), and take moments of the scalar curvature 
averaged only over vertices. This gives the same mean and includes
additional lift fluctuations in the variance. The standard-deviation
plots are similar but slightly higher.

Moreover, Conjecture~\ref{conj:star-asymptotic} predicts for $\beta=0$ that
\(\langle\overline{\Scal}\rangle \to1/4\) as $n\to\infty$. This is because, for the typical Eulerian graphs of
Theorem~\ref{thm:typical-eulerian}, every vertex lies in a triangle
and \(d^*(G)=1\). The inner connection has zero scalar curvature.
For the other connection, \eqref{eq:binary-scalar} is a nonconstant
affine function of the lift parameters at each vertex, giving mean
\(1/2\) under the uniform lift average. Averaging the two connections
gives \(1/4\), and the conjecture makes the exceptional connections
negligible in the total count. For $\beta>0$, convergence as $n\to\infty$ is a separate question, as is the behaviour of the Lorentzian version.

\subsection{Geometric and bare Laplacians}
\label{sec:bare}

The geometric Laplacian for any QRG and a `bare' Laplacian that directly depends on the metric are respectively\cite{BegMa}
\[
 \Delta_\nabla=(\ ,\ )\nabla\extd,
 \qquad Lf=(\ ,\ )(\theta\ten\extd f).
\]
The latter is known \cite{Ma:gra} for any field to be closely related to the usual graph Laplacian in discrete geometry. We look at these constructions for QRG over $\F_2$. 

\begin{proposition}
For every QLC classified in Theorem~\ref{thm:general-alpha-kernel}, the geometric Laplacian vanishes,
\(\Delta_\nabla=0\). The bare Laplacian is
\[
 (Lf)(x)=\sum_{x\to y}(f(y)+f(x)),
 \qquad L=D+A,
\]
where \(A\) is the adjacency matrix over \(\F_2\) and
\(D=\operatorname{diag}(\deg(x)\bmod2)\).
For Eulerian graphs \(L=A\); for connected co-Eulerian graphs
\(L^2=0\).
\end{proposition}

\begin{proof}
Contracting the connection on a basis arrow gives
\[
 (\ ,\ )\nabla\omega_{xy}
 =\delta_y+\bigl(1+(q+1)\deg(x)\bigr)\delta_x
 =\delta_y+\delta_x.
\]
The \(\alpha\)-term has distinct endpoints and contributes nothing;
the last equality follows from the degree parity for \(q=0\),
and directly for \(q=1\). Therefore
\[
 \Delta_\nabla f
 =\sum_{x\to y}(f(y)+f(x))(\delta_x+\delta_y)=0,
\]
as the two orientations of each edge cancel. Contracting
\(\theta\ten\extd f\) gives the stated formula for \(L\).
For a connected co-Eulerian graph, the common degree parity
\(\delta\in\F_2\) gives \(D=\delta I\) and \(A^2=\delta I\).
Thus \(L^2=(A+\delta I)^2=0\).
\end{proof}

Although all eigenvalues vanish in the co-Eulerian case, the rank
\(r=\rank L\) determines its Jordan form:
\[
 L\sim J_2(0)^{\oplus r}\oplus 0^{\oplus(n-2r)},
 \qquad J_2(0)=\begin{pmatrix}0&1\\0&0\end{pmatrix},
 \qquad n=|V|.
\]
Indeed, \(L^2=0\) permits only Jordan blocks of size one or two,
and each size-two block contributes one to the rank. The quotient
\[
 \ker L/\im L,
 \qquad \dim_{\F_2}(\ker L/\im L)=n-2r
\]
therefore counts the size-one blocks. For example, \(C_4\) and
\(K_4\) have ranks two and one, respectively, and hence different
Jordan forms despite having the same characteristic polynomial \(t^4\).

The kernel has a graph-theoretic interpretation without assuming
co-Eulerianity. Let \(B:\F_2^{E}\to\F_2^{V}\) be the incidence matrix for a graph $G=(V,E)$,  with one column for each undirected edge, containing a one
at each endpoint. A binary edge vector specifies an edge subset,
with addition given by symmetric difference. The \emph{cycle space}
\(\ker B\) consists of subsets having even degree at every vertex,
equivalently sums of cycles. The \emph{cut space} \(\im B^T\)
consists of all edges crossing from a vertex subset \(S\) to its
complement: this cut is \(B^T\chi_S\) where $\chi_S$ is the characteristic function of $S$. Their intersection
\(\mathcal B(G)=\ker B\cap\im B^T\) is the \emph{bicycle space},
whose elements are cuts which also have even degree at every vertex.
For example, all four edges of \(C_4\) form both a cycle and the
cut separating its two pairs of opposite vertices.

\begin{proposition}
\label{prop:laplacian-bicycle}
For a connected graph, \(L=BB^T\), and \(B^T\) induces an isomorphism
\[
 \ker L/\Span\{1\}\cong\mathcal B(G),
 \qquad \dim\ker L=1+\dim\mathcal B(G),
\]
where \(1\) denotes the constant function.
\end{proposition}

\begin{proof}
The diagonal entries of \(BB^T\) are the degree parities and its
off-diagonal entries are those of \(A\). Thus \(Lf=0\) precisely
when the cut \(B^Tf\) lies in \(\ker B\). Every bicycle is such a
cut, and connectedness gives \(\ker B^T=\Span\{1\}\).
\end{proof}

For general Eulerian graphs, \(L=A\) need not be nilpotent.
Its characteristic and minimal polynomials over \(\F_2\), together
with \(\dim\ker L^k\), give spectral invariants; the latter determine
the Jordan block sizes at eigenvalue zero.

For real statistical weights, one can instead use the ordinary
Laplacian \(L_{\R}=D_{\R}-A_{\R}\), where
\(D_{\R}=\operatorname{diag}(d_x)\), \(d_x=\deg(x)\), and
\(A_{\R}\) is the adjacency matrix over \(\R\). The degrees here
are their full integer values, not their parities.
Its reduction modulo two is \(L\), but it
retains information lost in that reduction. On a connected graph,
write its real eigenvalues as \(0=\lambda_1<\lambda_2\leq\cdots\leq\lambda_n\).
For example, the heat trace and the determinant with the constant
mode removed are
\[
 K_G(t)=\sum_{j=2}^n e^{-t\lambda_j}\quad(t>0),
 \qquad \det{}'L_{\R}=\prod_{j=2}^n\lambda_j.
\]
Either \(K_G(t)\) at a chosen scale or \(\log\det{}'L_{\R}\)
could supply further spectral data that depend on the graph alone and could be used in the construction of other candidates for $S[G]$.

\section{Digital realisation of QRGs}
\label{app:digital}

We realise the smallest curved star-preserving QRG over \(\F_2\),
on \(C_3\), as a classical digital process. Each basis vector labels
one wire, carrying its coefficient \(0\) or \(1\). A linear map is
therefore a binary circuit built from XOR gates and copies of wires.
We first give the individual gates, then their tensor extensions and
finally the curvature process.

\subsection{The individual gates}
\label{ssec:digital-gates}

\begin{figure}[tbp]
\centering
\begin{tikzpicture}[x=1cm,y=1cm,line width=.65pt,
  font=\small,
  xor gate/.style={circle,draw=black,fill=white,minimum size=4.8mm,
                   inner sep=0pt,font=\large}]
\foreach \j/\edge/\colour in
 {0/01/blue!75!black,1/02/orange!80!black,2/10/green!45!black,
  3/12/violet,4/20/red!75!black,5/21/cyan!60!black}{
  \coordinate (bus\edge) at ({.55*\j},.6);
  \draw[\colour] ({.55*\j},.6) -- ({.55*\j},-8.6);
  \node[above] at (bus\edge) {$\omega_{\edge}$};
}
\foreach \row/\u/\v/\w/\ret/\nonret/\cu/\cv/\cw in {
  0/1/2/3/010/102/orange!80!black/green!45!black/violet,
  1/0/3/1/101/012/blue!75!black/violet/orange!80!black,
  2/0/4/5/020/201/blue!75!black/red!75!black/cyan!60!black,
  3/1/5/0/202/021/orange!80!black/cyan!60!black/blue!75!black,
  4/3/4/2/212/120/violet/red!75!black/green!45!black,
  5/2/5/4/121/210/green!45!black/cyan!60!black/red!75!black}{
  \begin{scope}[yshift={-1.5*\row cm}]
  \node[xor gate] (first) at (4.8,0) {$+$};
  \node[xor gate] (second) at (7.5,-.55) {$+$};
  \draw[\cu] ({.55*\u},.17) -- (first.north west);
  \fill[\cu] ({.55*\u},.17) circle[radius=1.2pt];
  \draw[\cv] ({.55*\v},-.17) -- (first.south west);
  \fill[\cv] ({.55*\v},-.17) circle[radius=1.2pt];
  \draw (first.east) -- (5.5,0) -- (5.8,0) -- (5.8,.20)
        -- (9.3,.20) node[right] {$\omega_{\ret}$};
  \fill (5.5,0) circle[radius=1.2pt];
  \draw (5.5,0) |- (second.north west);
  \draw[\cw] ({.55*\w},-.88) -- (7.05,-.88)
        -- (second.south west);
  \fill[\cw] ({.55*\w},-.88) circle[radius=1.2pt];
  \draw (second.east) -- (9.3,-.55) node[right] {$\omega_{\nonret}$};
  \end{scope}
}
\end{tikzpicture}
\caption{The curved star-preserving \(C_3\) connection. Circled
pluses are XOR gates; dots mark junctions, while unmarked crossings
are not connected.}
\label{fig:K3-digital}
\end{figure}
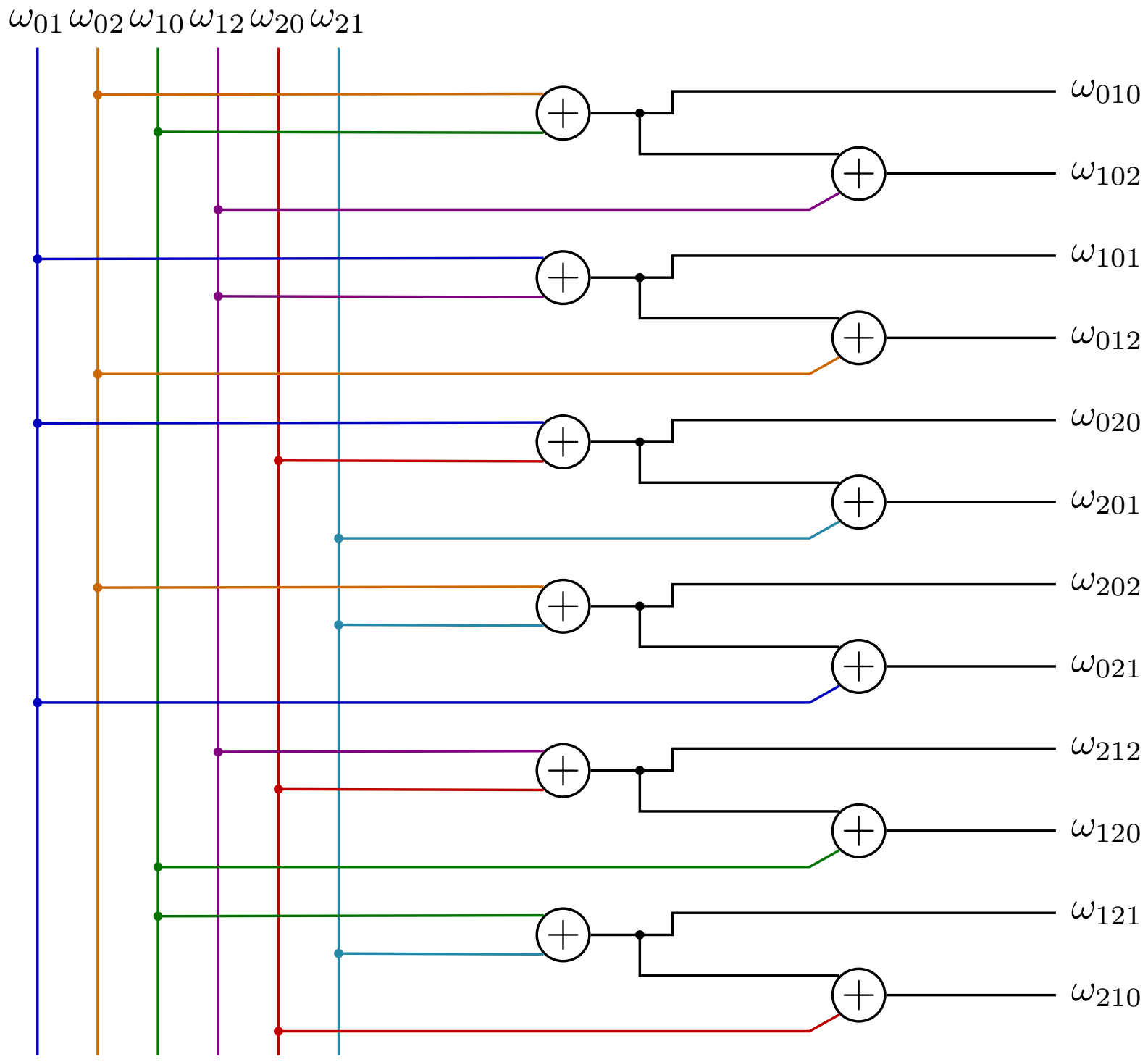

Take vertices \(0,1,2\) and use the ordered basis of \(\Omega^1\)
\[
 \mathcal B=(\omega_{01},\omega_{02},\omega_{10},\omega_{12},\omega_{20},\omega_{21}).
\]
The input bundle has six wires. With
\(\omega_{abc}=\omega_{ab}\ten\omega_{bc}\), order the twelve-wire
basis of \(\Omega^1\ten\Omega^1\) as
\[
 \bigl(\omega_{010},\omega_{102},\omega_{101},\omega_{012},
       \omega_{020},\omega_{201},
       \omega_{202},\omega_{021},\omega_{212},\omega_{120},
       \omega_{121},\omega_{210}\bigr).
\]
The curved connection in Remark~\ref{rem:K3twoparam}, with
\(t_1=t_2=1\), is
\begin{align*}
 \nabla\omega_{01}&=\omega_{012}+\omega_{020}+\omega_{021}+\omega_{101}+\omega_{201},\\
 \nabla\omega_{02}&=\omega_{010}+\omega_{012}+\omega_{021}+\omega_{102}+\omega_{202},\\
 \nabla\omega_{10}&=\omega_{010}+\omega_{102}+\omega_{120}+\omega_{121}+\omega_{210},\\
 \nabla\omega_{12}&=\omega_{012}+\omega_{101}+\omega_{102}+\omega_{120}+\omega_{212},\\
 \nabla\omega_{20}&=\omega_{020}+\omega_{120}+\omega_{201}+\omega_{210}+\omega_{212},\\
 \nabla\omega_{21}&=\omega_{021}+\omega_{121}+\omega_{201}+\omega_{202}+\omega_{210}.
\end{align*}
Its twelve XOR gates are shown in Figure~\ref{fig:K3-digital}, with
output wires in the order above.

\begin{figure}[tbp]
\centering
\begin{tikzpicture}[x=1cm,y=1cm,font=\small,line width=.55pt,
 xor gate/.style={circle,draw,fill=white,minimum size=4.8mm,
                  inner sep=0pt,font=\large}]
\node at (2.5,.7) {\(\extd\)};
% The six input wires are in the stated Omega^1 basis order.
\foreach \i/\edge in {0/01,1/02,2/10,3/12,4/20,5/21}{
 \coordinate (in\edge) at (0,{-.6*\i});
 \node[left] at (in\edge) {$\omega_{\edge}$};
}
\foreach \i/\a/\b in {0/01/10,1/02/20,2/12/21}{
 \node[xor gate] (s\i) at (1.5,{-1.5*\i}) {$+$};
 \draw (in\a) -- (s\i.north west);
 \draw (in\b) -- (s\i.south west);
 \draw (s\i.east) -- ({2.1+.25*\i},{-1.5*\i});
 \fill ({2.1+.25*\i},{-1.5*\i}) circle[radius=1pt];
 \node[xor gate] (z\i) at (4.2,{-1.5*\i}) {$+$};
 \draw (z\i.east) -- (5,{-1.5*\i}) node[right] {$\eta_\i$};
}
\draw (2.1,0) -- (3.4,0) -- (z0.north west);
\draw (2.1,0) -- (2.1,-.8) -- (z1.north west);
\draw (2.35,-1.5) -- (3.2,-.7) -- (z0.south west);
\draw (2.35,-1.5) -- (2.35,-2.3) -- (z2.north west);
\draw (2.6,-3) -- (2.6,-2.3) -- (z1.south west);
\draw (2.6,-3) -- (3.4,-3) -- (z2.south west);
\begin{scope}[xshift=7.2cm]
\node at (1.3,.7) {\(\wedge\)};
\foreach \i/\pathlabel in {0/010,1/102,2/101,3/012,4/020,5/201,
 6/202,7/021,8/212,9/120,10/121,11/210}{
 \coordinate (p\pathlabel) at (0,{-.36*\i});
 \node[left] at (p\pathlabel) {$\omega_{\pathlabel}$};
}
\foreach \i/\a/\b in {0/010/020,1/101/121,2/202/212}{
 \node[xor gate] (w\i) at (1.9,{-1.5*\i}) {$+$};
 \draw (p\a) -- (w\i.north west);
 \draw (p\b) -- (w\i.south west);
 \draw (w\i.east) -- (2.8,{-1.5*\i}) node[right] {$\eta_\i$};
}
\foreach \i/\pathlabel in {1/102,3/012,5/201,7/021,9/120,11/210}{
 \draw (p\pathlabel) -- (.45,{-.36*\i});
 \draw (.45,{-.36*\i-.07}) -- (.45,{-.36*\i+.07});
}
\end{scope}
\begin{scope}[xshift=.8cm,yshift=-6.6cm]
\node at (1.6,1.8) {copy};
\foreach \i/\colour in {1/blue!75!black,2/orange!80!black,
 3/green!45!black,4/violet}{
 \pgfmathsetmacro{\yin}{.3-.65*(\i-1)}
 \pgfmathsetmacro{\yu}{1.2-.45*(\i-1)}
 \pgfmathsetmacro{\yl}{-2.2-.45*(\i-1)}
 \draw[\colour] (0,\yin) node[left,text=black] {$v_\i$} -- (1,\yin);
 \fill[\colour] (1,\yin) circle[radius=1.2pt];
 \draw[\colour] (1,\yin) .. controls (1.8,\yin) and (2,\yu) .. (3,\yu)
   node[right,text=black] {$v_\i$};
 \draw[\colour] (1,\yin) .. controls (1.8,\yin) and (2,\yl) .. (3,\yl)
   node[right,text=black] {$v_\i$};
}
\begin{scope}[xshift=6.4cm]
\node at (1.4,1.8) {addition \(\oplus_4\)};
\foreach \i/\colour in {1/blue!75!black,2/orange!80!black,
 3/green!45!black,4/violet}{
 \pgfmathsetmacro{\level}{.3-.65*(\i-1)}
 \pgfmathsetmacro{\yu}{1.2-.45*(\i-1)}
 \pgfmathsetmacro{\yl}{-2.2-.45*(\i-1)}
 \node[xor gate] (plus\i) at (2,\level) {$+$};
 \draw[\colour] (0,\yu) node[left,text=black] {$v_\i$}
   .. controls (1,\yu) and (1.2,{\level+.17}) .. (plus\i.north west);
 \draw[\colour] (0,\yl) node[left,text=black] {$v_\i$}
   .. controls (1,\yl) and (1.2,{\level-.17}) .. (plus\i.south west);
 \draw[\colour] (plus\i.east) -- (3,\level)
   node[right,text=black] {$v_\i$};
}
\end{scope}
\end{scope}
\end{tikzpicture}
\caption{The \(\extd\) and \(\wedge\) gates (top); copy and addition banks (bottom). Bars terminate unused wires.}
\label{fig:K3-form-gates}
\label{fig:digital-banks}
\end{figure}
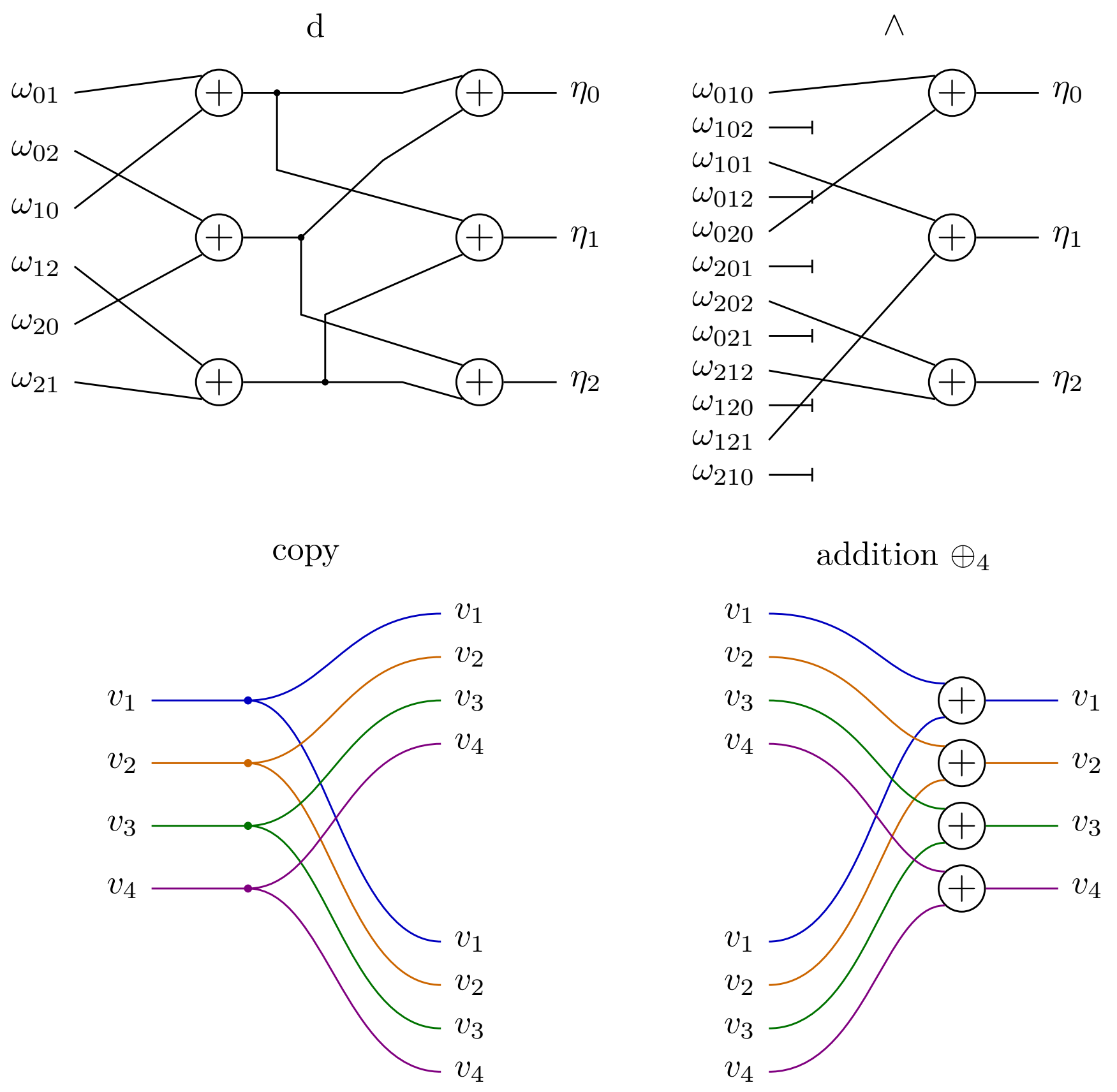

The ordered basis for the three-wire bundle \(\Omega^2_{\min}\) is
\[
 (\omega_{01}\wedge\omega_{10},\,
  \omega_{12}\wedge\omega_{21},\,
  \omega_{20}\wedge\omega_{02})=(\eta_0,\eta_1,\eta_2).
\]
The relations identify the two return paths at each vertex:
\[
 \omega_{02}\wedge\omega_{20}=\eta_0,\qquad
 \omega_{10}\wedge\omega_{01}=\eta_1,\qquad
 \omega_{21}\wedge\omega_{12}=\eta_2,
\]
while every non-returning two-path vanishes. Thus
\(\extd\omega_{ab}=\eta_a+\eta_b\), and \(\wedge\) adds the two
return-path coefficients at each vertex. Their internal circuits are
shown in Figure~\ref{fig:K3-form-gates} (top).

Copy and addition will also be needed to xerox or add elements of a vector space $V$ with basis $v_1,\cdots,v_n$, say. This is illustrated in Figure~\ref{fig:digital-banks} (bottom) for $n=4$. A vector is carried in a bundle of 4 wires with each wire labelled by its basis vector (that element is included if the wire is switched on). 

\subsection{Tensoring a gate with the identity}
\label{ssec:digital-tensor}

Let \(V,W,Z\) be finite-dimensional \(A\)-bimodules with ordered
bases in their endpoint components. A right \(A\)-linear map
\(\phi:V\to W\) induces
\(\phi\otimes_A\id:V\ten Z\to W\ten Z\).
For each basis vector \(z\in\delta_aZ\), use the restricted gate
\(\phi:V\delta_a\to W\delta_a\), attaching the unchanged factor
\(z\) to every input and output label. No other wires occur in this
block. Thus tensoring with the identity repeats the appropriate
endpoint block of the gate for each basis vector of \(Z\).

\begin{figure}[tbp]
\centering
\begin{tikzpicture}[x=1cm,y=1cm,font=\small,line width=.55pt,
 xor gate/.style={circle,draw,fill=white,minimum size=4.8mm,
                  inner sep=0pt,font=\large}]
\node at (1.5,.9) {$\wedge$};
\node[xor gate] (plain) at (1.5,0) {$+$};
\draw (0,.3) node[left] {$\omega_{010}$} -- (plain.north west);
\draw (0,-.3) node[left] {$\omega_{020}$} -- (plain.south west);
\draw (plain.east) -- (3,0) node[right] {$\eta_0$};
\foreach \level/\pathlabel in {-.9/120,-1.5/210}{
 \draw (0,\level) node[left] {$\omega_{\pathlabel}$} -- (.55,\level);
 \draw (.55,{\level-.07}) -- (.55,{\level+.07});
}
\begin{scope}[xshift=7cm]
\node at (1.5,.9) {$\wedge\otimes_A\id$};
\node[xor gate] (tensor) at (1.5,0) {$+$};
\draw (0,.3) node[left] {$\omega_{010}\ten\omega_{01}$} -- (tensor.north west);
\draw (0,-.3) node[left] {$\omega_{020}\ten\omega_{01}$} -- (tensor.south west);
\draw (tensor.east) -- (3,0) node[right] {$\eta_0\ten\omega_{01}$};
\foreach \level/\pathlabel in {-.9/120,-1.5/210}{
 \draw (0,\level) node[left] {$\omega_{\pathlabel}\ten\omega_{01}$} -- (.55,\level);
 \draw (.55,{\level-.07}) -- (.55,{\level+.07});
}
\end{scope}
\begin{scope}[yshift=-3.3cm]
\node at (1.3,.7) {$J$};
\draw (0,-1.1) node[left] {$\omega_{010}$} -- (2.6,-1.1)
 node[right] {$\omega_{01}\otimes_{\F_2}\omega_{10}$};
\draw (0,-1.65) node[left] {$\omega_{012}$} -- (2.6,-1.65)
 node[right] {$\omega_{01}\otimes_{\F_2}\omega_{12}$};
\foreach \level/\edge in {0/01,-.55/02,-2.2/20,-2.75/21}{
 \draw (.7,\level) node[left] {$0$} -- (2.6,\level)
  node[right] {$\omega_{01}\otimes_{\F_2}\omega_{\edge}$};
}
\begin{scope}[xshift=7.8cm]
\node at (.9,.7) {$P_2$};
\foreach \level/\edge in {0/01,-.55/02}{
 \draw (0,\level) node[left] {$\eta_0\otimes_{\F_2}\omega_{\edge}$}
  -- (1.8,\level) node[right] {$\eta_0\ten\omega_{\edge}$};
}
\foreach \level/\edge in {-1.1/10,-1.65/12,-2.2/20,-2.75/21}{
 \draw (0,\level) node[left] {$\eta_0\otimes_{\F_2}\omega_{\edge}$} -- (.55,\level);
 \draw (.55,{\level-.07}) -- (.55,{\level+.07});
}
\end{scope}
\end{scope}
\end{tikzpicture}
\caption{Endpoint blocks of \(\wedge\), \(\wedge\otimes_A\id\) (top), and \(J,P_2\) (bottom). Zeros are constant inputs.}
\label{fig:digital-tensor-wedge}
\label{fig:digital-lift-project}
\end{figure}
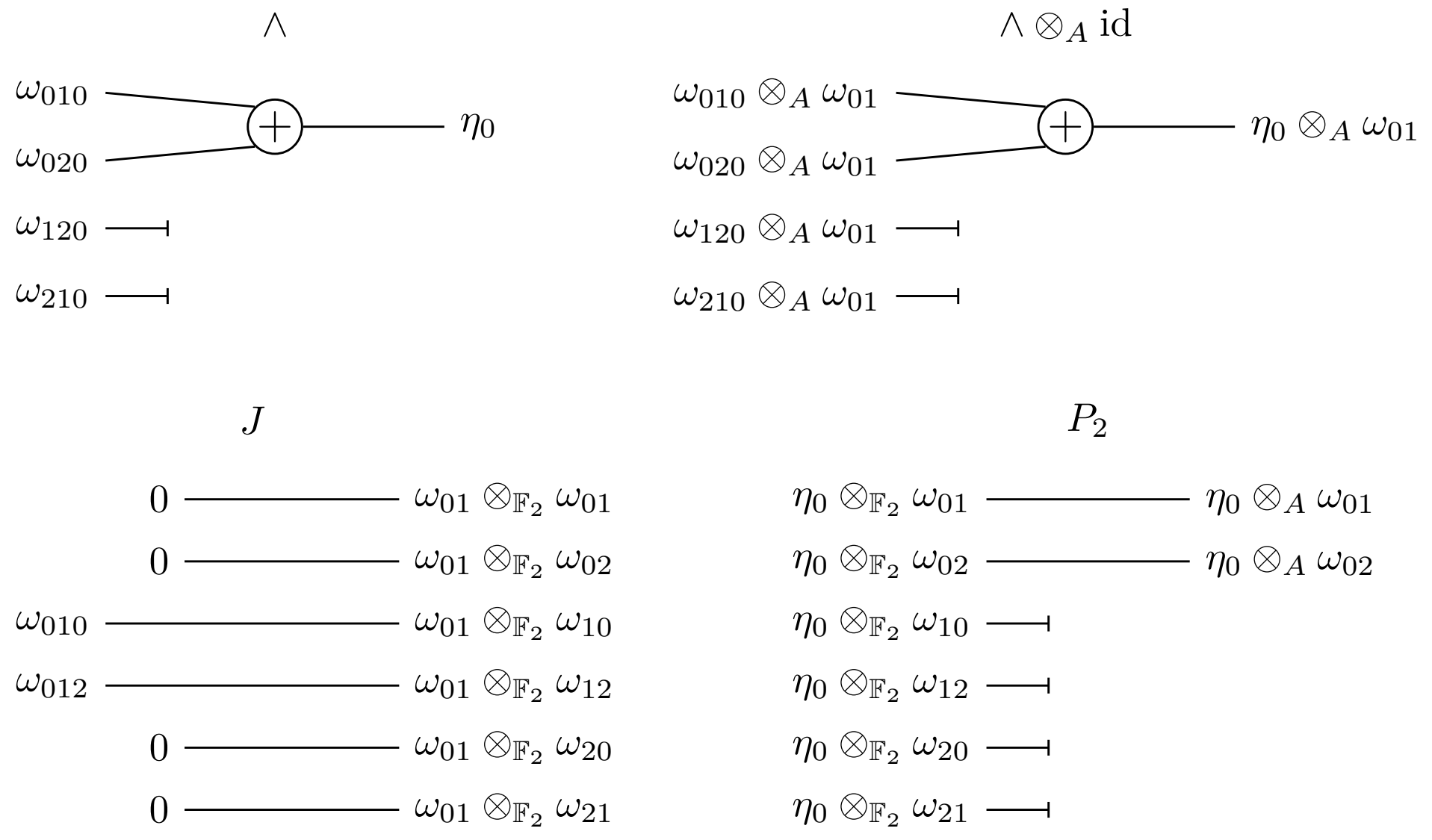

For example, \(\wedge\) is a bimodule map, so
\[
 \wedge\otimes_A\id:
 \Omega^1\ten\Omega^1\ten\Omega^1\longrightarrow
 \Omega^2_{\min}\ten\Omega^1
\]
needs no lifting or projection. Fix the last factor to be
\(\omega_{01}\). The possible first two factors are
\[
 (\omega_{010},\omega_{020},\omega_{120},\omega_{210}).
\]
The first two map to \(\eta_0\), while the last two map to zero.
Figure~\ref{fig:digital-tensor-wedge} (top) shows how the same XOR gate
acts with \(\omega_{01}\) attached to every label. Repeating this
for \(01,02,10,12,20,21\) gives the full 24-input, six-output gate:
there are four inputs and one output for each fixed last arrow.

When \(\phi\) is only \(\F_2\)-linear, it may be one part of a
process which descends to \(\otimes_A\), although
\(\phi\otimes\id\) does not descend separately. Write
\[
 \pi_{V,Z}:V\otimes_{\F_2}Z\longrightarrow V\ten Z
\]
for the quotient by \(vf\otimes z=v\otimes fz\).
Choose a representative for each basis vector of \(V\ten Z\) and
extend linearly to a section
\[
 J:V\ten Z\longrightarrow V\otimes_{\F_2}Z,\qquad
 \pi_{V,Z}J=\id.
\]
The choice is fixed for the whole process. If
\(P:W\otimes_{\F_2}Z\to W\ten Z\) is the output quotient, the
corresponding piece of the circuit is \(P(\phi\otimes_{\F_2}\id)J\).
Unlike \(P\), this piece can depend on the choice of \(J\).

For \(V=Z=\Omega^1(C_3)\), use \(\mathcal B\) from
Section~\ref{ssec:digital-gates} and order
\(\mathcal B\otimes_{\F_2}\mathcal B\) by the first factor,
then the second. These are 36 basis vectors. The quotient has the
12 two-path basis vectors of Section~\ref{ssec:digital-gates}, and
\[
 \pi_{\Omega^1,\Omega^1}(\omega_{ab}\otimes_{\F_2}\omega_{cd})
 =\begin{cases}\omega_{abd},&b=c,\\0,&b\ne c.\end{cases}
\]
We choose
\[
 J(\omega_{abc})=\omega_{ab}\otimes_{\F_2}\omega_{bc}.
\]
Thus \(J\) routes each path coefficient to its composable slot and
inserts 24 zero coefficients. The relations do not require those coefficients
to be zero in every representative: noncomposable tensors lie in the
kernel of \(\pi\) and could be added without changing its image.
Zero insertion is our fixed choice. Figure~\ref{fig:digital-lift-project}
(lower left) shows the block with first factor \(\omega_{01}\); the other five
blocks have the same form, with the two live slots determined by the
end of the first arrow.

The output quotient needed for \(\extd\otimes\id\) is
\[
 P_2:\Omega^2_{\min}\otimes_{\F_2}\Omega^1
       \longrightarrow\Omega^2_{\min}\ten\Omega^1,
 \qquad
 P_2(\eta_a\otimes_{\F_2}\omega_{bc})
 =\begin{cases}\eta_a\ten\omega_{ac},&a=b,\\0,&a\ne b.\end{cases}
\]
Order its 18 input basis vectors by \(\eta_0,\eta_1,\eta_2\), with
\(\mathcal B\) in each block. Its six outputs are
\[
 (\eta_0\ten\omega_{01},\eta_0\ten\omega_{02},
  \eta_1\ten\omega_{10},\eta_1\ten\omega_{12},
  \eta_2\ten\omega_{20},\eta_2\ten\omega_{21}).
\]
It retains these six wires and discards the other twelve.

For example, keeping the second factor \(\omega_{01}\) fixed gives
\begin{align*}
 P_2(\extd\otimes_{\F_2}\id)J(\omega_{101})
   &=P_2\bigl((\eta_1+\eta_0)\otimes_{\F_2}\omega_{01}\bigr)
     =\eta_0\ten\omega_{01},\\
 P_2(\extd\otimes_{\F_2}\id)J(\omega_{201})
   &=P_2\bigl((\eta_2+\eta_0)\otimes_{\F_2}\omega_{01}\bigr)
     =\eta_0\ten\omega_{01}.
\end{align*}
Hence the two input coefficients are XORed into the wire labelled
\(\eta_0\ten\omega_{01}\). Equivalently, pad their first-factor
inputs to six wires, apply the \(\extd\)-gate of
Section~\ref{ssec:digital-gates}, and retain its \(\eta_0\) output.

For the other branch of curvature we use the quotient
\[
 P_3:\Omega^1\otimes_{\F_2}(\Omega^1\ten\Omega^1)
       \longrightarrow\Omega^1\ten\Omega^1\ten\Omega^1.
\]
On \(\omega_{ab}\otimes_{\F_2}\omega_{cde}\), it keeps the tensor
when \(b=c\) and sends it to zero otherwise. This selects the 24
composable three-paths from 72 inputs. The quotient from
\((\Omega^1\ten\Omega^1)\otimes_{\F_2}\Omega^1\) is defined
in the same way, matching the last vertex of the two-path to the
first vertex of the arrow.

\subsection{The curvature process}
\label{ssec:digital-curv}

\begin{figure}[tbp]
\centering
\begin{tikzpicture}[x=1cm,y=1cm,font=\small,line width=.35pt,
 gate/.style={draw,line width=.6pt,minimum width=.8cm,
              minimum height=.7cm,inner sep=2pt},
 route/.style={draw,line width=.6pt,minimum width=.6cm,
              minimum height=.7cm,inner sep=2pt}]
\node[gate] (nabla) at (.8,0) {$\nabla$};
\node[route] (lift) at (2,0) {$J$};
\node[gate,minimum width=1cm,minimum height=1.1cm] (copy) at (3.3,0) {copy};
\node[route] (pu) at (7,1.6) {$P_2$};
\node[route,minimum width=1cm] (routing) at (7.5,-1.6) {$P_3$};
\node[gate,minimum height=1.1cm] (xor) at (13.2,0) {$\oplus_6$};
% Each small rectangle is one copy of a previously defined gate.
\foreach \x/\y/\op in {5.5/1.6/\extd,5.5/-1.6/\nabla,9.5/-1.6/\wedge}{
 \foreach \j in {0,...,5}{
  \pgfmathsetmacro{\level}{\y+.32*(2.5-\j)}
  \node[draw,minimum width=.65cm,minimum height=.25cm,
        inner sep=.2pt,font=\scriptsize] at (\x,\level) {$\op$};
  \draw ({\x-.55},\level) -- ({\x-.325},\level);
  \draw ({\x+.325},\level) -- ({\x+.55},\level);
 }
}
% Double rails denote coefficient bundles; adjacent numbers give their widths.
\foreach \o in {-.045,.045}{
 \draw (0,\o) -- (.4,\o);
 \draw (1.2,\o) -- (1.7,\o);
 \draw (2.3,\o) -- (2.8,\o);
 \draw (3.8,{.25+\o}) .. controls (4.15,{.25+\o}) and
  (4.15,{1.6+\o}) .. (4.65,{1.6+\o});
 \draw (3.8,{-.25+\o}) .. controls (4.15,{-.25+\o}) and
  (4.15,{-1.6+\o}) .. (4.65,{-1.6+\o});
 \draw (6.35,{1.6+\o}) -- (6.7,{1.6+\o});
 \draw (6.35,{-1.6+\o}) -- (7,{-1.6+\o});
 \draw (8,{-1.6+\o}) -- (8.65,{-1.6+\o});
 \draw (7.3,{1.6+\o}) -- (12,{1.6+\o})
  .. controls (12.4,{1.6+\o}) and (12.4,{.25+\o}) .. (12.8,{.25+\o});
 \draw (10.35,{-1.6+\o}) -- (12,{-1.6+\o})
  .. controls (12.4,{-1.6+\o}) and (12.4,{-.25+\o}) .. (12.8,{-.25+\o});
 \draw (13.6,\o) -- (14.6,\o);
}
% Fan the indexed slices into and out of the repeated-gate banks.
\foreach \x/\y in {5.5/1.6,5.5/-1.6,9.5/-1.6}{
 \foreach \j in {0,...,5}{
  \pgfmathsetmacro{\level}{\y+.32*(2.5-\j)}
  \pgfmathsetmacro{\port}{\y+.018*(2.5-\j)}
  \draw ({\x-.85},\port) -- ({\x-.55},\level);
  \draw ({\x+.55},\level) -- ({\x+.85},\port);
 }
}
\node[above] at (.2,.15) {$\omega$};
\node[below] at (.2,-.1) {6};
\node[above] at (1.45,.1) {12};
\node[above] at (2.55,.1) {36};
\node[above] at (4.3,1.75) {36};
\node[below] at (4.3,-1.75) {36};
\node[above] at (6.5,1.75) {18};
\node[below] at (6.65,-1.75) {72};
\node[below] at (8.3,-1.75) {24};
\node[above] at (9,1.75) {6};
\node[below] at (11.4,-1.75) {6};
\node[above] at (14.1,.15) {$R_\nabla(\omega)$};
\node[below] at (14.1,-.1) {6};
\node[above] at (5.5,2.65) {six copies};
\node[below] at (5.5,-2.65) {six copies};
\node[below] at (9.5,-2.65) {six wedge blocks};
\end{tikzpicture}
\caption{The curvature process. Numbers give bundle widths;
\(J\) inserts representatives; \(P_2,P_3\) select composable outputs.}
\label{fig:K3-curvature-circuit}
\end{figure}
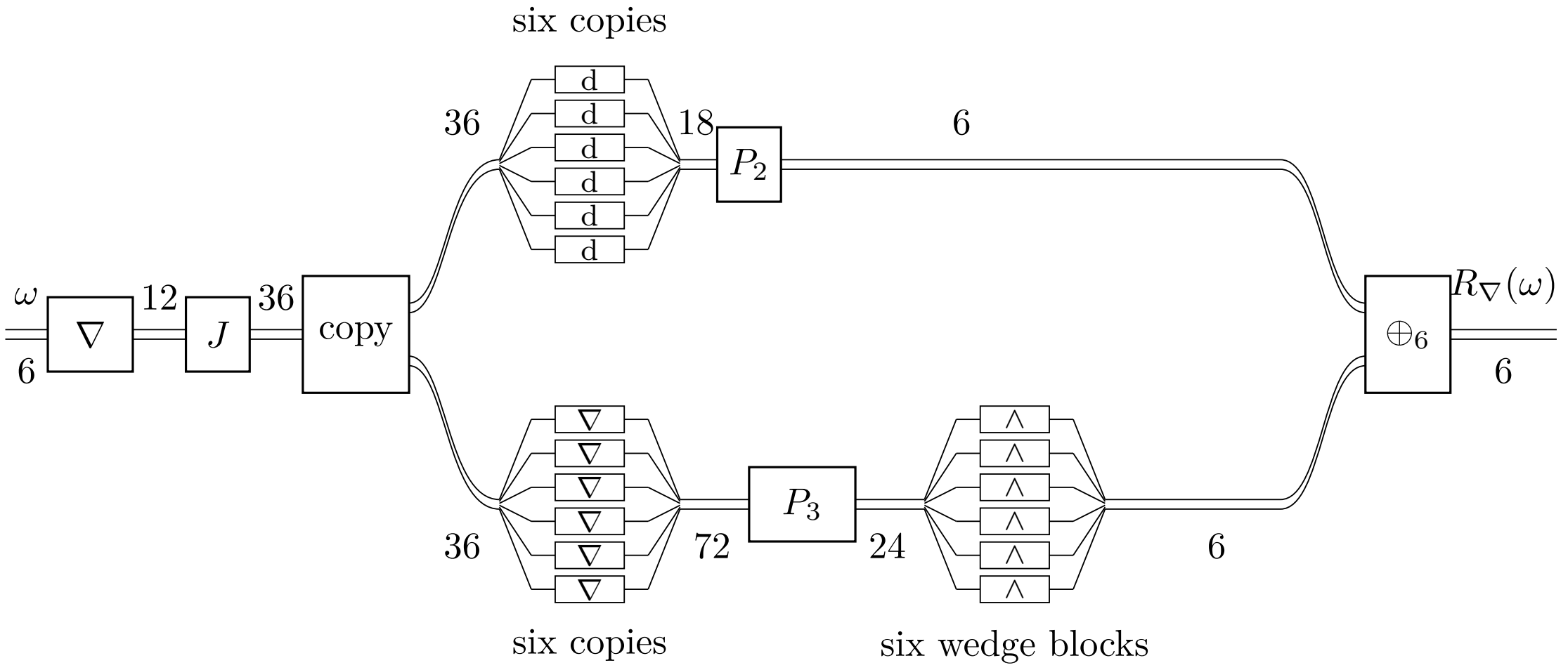

Apply \(\nabla\), lift its twelve-wire output by \(J\), and copy
the resulting 36-wire representative into two branches. The first
applies \(\extd\) to the first factor and then \(P_2\). The second
applies \(\nabla\) to the second factor, then \(P_3\) and
\(\wedge\otimes_A\id\). The two six-wire outputs, in the order
specified above, are added by \(\oplus_6\):
\[
 R_\nabla=\bigl[P_2(\extd\otimes_{\F_2}\id)
   +(\wedge\otimes_A\id)P_3
       (\id\otimes_{\F_2}\nabla)\bigr]J\nabla.
\]
This is the process in Figure~\ref{fig:K3-curvature-circuit}.

The two terms in brackets agree on the kernel of the input quotient,
so their sum is independent of the representative. Using the same
\(J\) after copying would be equivalent, but unrelated choices on
the two branches need not cancel. For example,
\(\omega_{01}\otimes_{\F_2}\omega_{02}\) represents zero, yet
each branch sends it to \(\eta_0\ten\omega_{02}\). Putting \(J\)
before the copy bank ensures that both branches use the same lift.

For this particular connection and the chosen bases,
\[
 R_\nabla(\omega_{ab})=\eta_a\ten\omega_{ab}.
\]
Hence the entire process simplifies to six straight wires in the
stated input and output bases. Such simplification can be done using a digital emulator to simplify the circuit developed geometrically as above, and similarly for other more complex processes. Over $\C$ one can similarly build gates and compose them but now more like quantum computer circuits. 

\section{Python scripts and catalogue}
\label{app:computational-classification}
The accompanying Python~3 programs provide solution archives,
statistical counts and human-readable catalogues\cite{github}. The file \texttt{Python\_doc.pdf} provided there documents the command lines,
JSON Lines and CSV formats, restart procedure and quiver identifiers, while here we provide a brief overview. 

\subsection{Direct enumeration: \texttt{qlc-enumerator.py}}
\label{app:enumerator}
The enumerator constructs all QLCs on connected graph representatives,
recording each \(\sigma\), \(\alpha\), curvature and star condition in
a machine-readable archive. The original polynomial solver supplied
the independent enumeration through \(n=9\), agreeing with
Tables~\ref{tab:graph-classes} and~\ref{tab:QLC-statistics}.
The present default uses the classification: it tests the Eulerian
and co-Eulerian conditions, then solves the linear \(\alpha\)-equations.
The original solver remains available for independent checks.
The machine-readable archive through nine vertices contains 19,323 QLCs,
of which 5,341 are star-preserving.

\subsection{Statistics and solution reports: \texttt{qlc-report.py}}
The reporter converts a solution archive into graph-level CSV statistics
and a \LaTeX\ catalogue of graphs, connections and curvatures. Reports
can be selected by vertex range or prescribed degree sequence, with
independent \(\alpha\)-parameters and, where available, a cyclic
Cayley presentation. For example, the degree sequence
\((4,4,4,4,4,4)\) selects the octahedral graph \(K_{2,2,2}\),
with four QLCs and two star-preserving QLCs for each braiding.
The full catalogue \texttt{solutions8.pdf} through eight vertices has 3,833 pages; restricted reports
are more convenient for inspecting individual families.

\subsection{Graph counts: \texttt{qlc-graph-counts.py}}
This program produced Tables~\ref{tab:graph-classes}
and~\ref{tab:QLC-statistics} through \(n=12\), without generating
individual connections. It streams connected graphs, applies the
parity tests and counts triangle components and flat connections.
Parallel workers save resumable counts rather than a solution database.
For example, at \(n=12\) it finds 86,539,128 Eulerian graphs,
325 co-Eulerian graphs and an intersection of size 159.
The Cayley column is computed separately by
\texttt{qlc-cayley-counts.py}, which identifies graph-isomorphic
presentations across all groups of the relevant orders.

\subsection{Eulerian counts: \texttt{eulerian-counts.py}}
For larger orders this program generates only connected Eulerian
graphs, by parity completion of graphs on one fewer vertex and
canonical identification. Its \texttt{--mode star} option gives the
histogram of \(d^*(G)\) and hence \(\sum_G2^{d^*(G)}\), together with
the co-Eulerian intersection. The \texttt{--mode full} option counts
\(d(G)\) and tests flatness, while \texttt{--mode all} computes both
in one pass. The star and full modes supplied the completed \(n=13\) results in two passes.
At odd order the Eulerian search covers every graph admitting QLCs,
including both braidings where they occur. The full mode therefore
also tests Conjecture~\ref{conj:odd-flat} by checking the inner braidings
on every connected Eulerian graph of the chosen order.

\subsection{Maps between QRGs: \texttt{qlc-qrg-quiver.py}}
\label{app:quiver-identifiers}
The quiver program reads a solution archive from
Appendix~\ref{app:enumerator} and tests vertex maps
for metric and connection preservation. Through six vertices it gives the results reported in Figures~\ref{fig:QRG-quiver-co-eulerian}
and~\ref{fig:QRG-quiver-eulerian}. Its JSON output retains every vertex map and
the independent \(\alpha\)-parameters at each node.
The companion \texttt{qlc-present-quiver.py} produces an interactive
HTML catalogue. Both programs accept \texttt{--type eulerian} and
\texttt{--type co-eulerian} to select the two cases separately. 
Graph drawings, graph archive codes and parameter supports are matched in the supplementary documentation.

\subsection{Cyclic Cayley families: \texttt{qlc-cyclic-catalogue.py}}
\label{app:cyclic-catalogue}
The cyclic generator constructs connected Cayley graphs of \(\Z_n\),
removes graph-isomorphic duplicates and computes the full and
star-preserving \(\alpha\)-kernels directly. Its machine-readable
catalogue through nine vertices retains the supports of every kernel
basis vector. Table~\ref{tab:cyclic-QLC-catalogue} illustrates the
output through six vertices as a small sample. Here
\[
 \sigma_q(e_i\ten e_j)=e_i\ten e_j+(q+\delta_{i+j,0})D_{i+j},
 \qquad D_s=\sum_{\substack{a,b\in \CC\\a+b=s}}e_a\ten e_b.
\]
The symbol \(0\) denotes \(\alpha=0\), and \(T\) denotes
\(\alpha(e_i)=tD_i\) for all \(i\in \CC\), with \(t\in\F_2\).
For \(P_J\), where \(\CC=J\sqcup(-J)\), set
\[
 \alpha(e_i)=
 \begin{cases}t_+D_i&i\in J,\\t_-D_i&i\in-J,\end{cases}
 \qquad t_+,t_-\in\F_2.
\]
The \(T\)-family is star-preserving; \(P_J\) is star-preserving
exactly when \(t_+=t_-\). A dash denotes absence of a QLC.

\begin{table}[tbp]
\centering
\begin{tabular}{r l c c c l}
\toprule
\(n\)&Generating set \(\CC\)&\(q\)&\(d\)&\(d_*\)&\(\alpha\)-family\\
\midrule
2&\(\{1\}\)&\(1\)&0&0&\(0\)\\
3&\(\{\pm1\}\)&\(0\)&2&1&\(P_{\{1\}}\)\\
4&\(\{\pm1\}\)&\(0,1\)&0&0&\(0\)\\
4&\(\{\pm1,2\}\)&\(1\)&1&1&\(T\)\\
5&\(\{\pm1\}\)&\(0\)&0&0&\(0\)\\
5&\(\{\pm1,\pm2\}\)&\(0\)&1&1&\(T\)\\
6&\(\{\pm1\}\)&\(0\)&0&0&\(0\)\\
6&\(\{\pm1,3\}\)&\(-\)&-&-&\(-\)\\
6&\(\{\pm2,3\}\)&\(-\)&-&-&\(-\)\\
6&\(\{\pm1,\pm2\}\)&\(0,1\)&2&1&\(P_{\{1,-2\}}\)\\
6&\(\{\pm1,\pm2,3\}\)&\(1\)&1&1&\(T\)\\
\bottomrule
\end{tabular}
\caption{Cyclic Cayley graphs through six vertices.}
\label{tab:cyclic-QLC-catalogue}
\end{table}

\end{document}